\input amstex
\magnification=1200
\loadmsbm
\loadeufm
\loadeusm
\UseAMSsymbols
\input amssym.def

\font\BIGtitle=cmr10 scaled\magstep3
\font\bigtitle=cmr10 scaled\magstep1
\font\boldsectionfont=cmb10 scaled\magstep1
\font\section=cmsy10 scaled\magstep1

\def\scr#1{{\fam\eusmfam\relax#1}}
\def\scrA{{\scr A}}
\def\scrB{{\scr B}}
\def\scrC{{\scr C}}
\def\scrD{{\scr D}}
\def\scrE{{\scr E}}
\def\scrF{{\scr F}}
\def\scrG{{\scr G}}
\def\scrH{{\scr H}}
\def\scrI{{\scr I}}
\def\scrK{{\scr K}}
\def\scrJ{{\scr J}}
\def\scrL{{\scr L}}
\def\scrM{{\scr M}}
\def\scrN{{\scr N}}

\def\scrP{{\scr P}}

\def\scrQ{{\scr Q}}
\def\scrS{{\scr S}}
\def\scrU{{\scr U}}
\def\scrV{{\scr V}}
\def\scrZ{{\scr Z}}
\def\scrW{{\scr W}}
\def\scrR{{\scr R}}
\def\scrT{{\scr T}}
\def\scrX{{\scr X}}
\def\scrY{{\scr Y}}
\def\gr#1{{\fam\eufmfam\relax#1}}

\def\grA{{\gr A}}	
	
\def\grC{{\gr C}}	
\def\grD{{\gr D}}	
	\def\gre{{\gr e}}
\def\grF{{\gr F}}	
\def\grG{{\gr G}}	
\def\grH{{\gr H}}	
\def\grI{{\gr I}}	
\def\grJ{{\gr J}}	
	
\def\grL{{\gr L}}	

\def\grM{{\gr M}}	\def\grm{{\gr m}}
\def\grN{{\gr N}}	
\def\grO{{\gr O}}	
\def\grP{{\gr P}}	 
\def\grQ{{\gr Q}}	
\def\grR{{\gr R}}	
\def\grS{{\gr S}}	\def\grs{{\gr s}}
\def\grT{{\gr T}}	
\def\grU{{\gr U}}	
	
\def\grW{{\gr W}}

	\def\grz{{\gr z}}

\def\db#1{{\fam\msbfam\relax#1}}

\def\dbA{{\db A}} 
\def\dbC{{\db C}} 
 \def\dbF{{\db F}}
\def\dbG{{\db G}} \def\dbH{{\db H}}

 \def\dbN{{\db N}}
 
\def\dbQ{{\db Q}} \def\dbR{{\db R}}

 \def\dbZ{{\db Z}}

\def\Rtil{\widetilde{R}}

\def\Ker{\text{Ker}}
\def\der{\text{der}}
\def\Sh{\hbox{\rm Sh}}

\def\sc{\text{sc}}
\def\Res{\text{Res}}
\def\ab{\text{ab}}
\def\ad{\text{ad}}

\def\Gal{\text{Gal}}
\def\GL{\text{GL}}
\def\Hom{\text{Hom}}
\def\End{\text{End}}

\def\Spec{\text{Spec}\,}
\def\Spf{\text{Spf}\,}

\def\Lie{\text{Lie}}

\def\leaderfill{\leaders\hbox to 1em
     {\hss.\hss}\hfill}
\def\nspace{\lineskip=1pt\baselineskip=12pt\lineskiplimit=0pt}

\def\endproof{$\hfill \square$}
\def\finishproclaim{\par\rm
     \ifdim\lastskip<\medskipamount\removelastskip
     \penalty55\medskip\fi}
\def\proof{\par\noindent {\it Proof:}\enspace}
\def\references#1{\par
  \centerline{\boldsectionfont References}\smallskip
     \parindent=#1pt\nspace}
\def\Ref[#1]{\par\hang\indent\llap{\hbox to\parindent
     {[#1]\hfil\enspace}}\ignorespaces}
\def\Item#1{\par\smallskip\hang\indent\llap{\hbox to\parindent
     {#1\hfill$\,\,$}}\ignorespaces}
\def\ItemItem#1{\par\indent\hangindent2\parindent
     \hbox to \parindent{#1\hfill\enspace}\ignorespaces}

\def\Le{{\mathchoice{\,{\scriptstyle\le}\,}
  {\,{\scriptstyle\le}\,}
  {\,{\scriptscriptstyle\le}\,}{\,{\scriptscriptstyle\le}\,}}}
\def\Ge{{\mathchoice{\,{\scriptstyle\ge}\,}
  {\,{\scriptstyle\ge}\,}
  {\,{\scriptscriptstyle\ge}\,}{\,{\scriptscriptstyle\ge}\,}}}

\def\arrowsim{\,\smash{\mathop{\to}\limits^{\lower1.5pt
  \hbox{$\scriptstyle\sim$}}}\,}

\def\doublemaprights#1#2#3#4{\raise3pt\hbox{$\mathop{\,\,\hbox to
     #1pt{\rightarrowfill}\kern-30pt\lower3.95pt\hbox to
     #2pt{\rightarrowfill}\,\,}\limits_{#3}^{#4}$}}

\def\rightcapdownarrow{\raise9pt\hbox{$\ssize\cap$}\kern-7.75pt
     \Big\downarrow}

\def\rcapmapdown#1{\rightcapdownarrow\kern-1.0pt\vcenter{
     \hbox{$\scriptstyle#1$}}}

\def\rmapdown#1{\Big\downarrow\kern-1.0pt\vcenter{
     \hbox{$\scriptstyle#1$}}}
\def\rightsubsetarrow#1{{\ssize\subset}\kern-4.5pt\lower2.85pt
     \hbox to #1pt{\rightarrowfill}}
\def\longtwoheadedrightarrow#1{\raise2.2pt\hbox to #1pt{\hrulefill}
     \!\!\!\twoheadrightarrow}

\def\Gal{\operatorname{\hbox{Gal}}}
\def\Hom{\operatorname{\hbox{Hom}}}

\def\Im{\hbox{Im}}

\NoBlackBoxes
\parindent=25pt
\document
\footline={\hfil}

\null
\centerline{\BIGtitle CM lifts for Isogeny Classes of}

\centerline{\BIGtitle Shimura F-crystals over Finite Fields}

\vskip 0.1in 
\centerline{\bigtitle Adrian Vasiu}

\vskip 0.1in
\centerline{June 19, 2012}
\footline={\hfill}
\null

\noindent
{\bf ABSTRACT}. We extend to large contexts pertaining to Shimura varieties of Hodge type a result of Zink on the existence of lifts to characteristic $0$ of suitable representatives of certain isogeny classes of abelian varieties endowed with Frobenius and other endomorphisms over $\dbF_{p^q}$, whose $p$-divisible groups in mixed characteristic $(0,p)$ are with complex multiplication. 



\footline={\hss\tenrm \folio\hss}
\pageno=1

\bigskip
\noindent
{\boldsectionfont 1. Introduction}

\bigskip
Let $p\in\dbN^*$ be a prime. Let $r\in\dbN^*$. Let $q:=p^r$. Let $k:=\dbF_q$ be the field with $q$ elements. Let $W(k)$ be the ring of Witt vectors  with coefficients in $k$. Let $B(k):=W(k)[{1\over p}]$ be the field of fractions of $W(k)$. Let $\sigma:=\sigma_k$ be the Frobenius automorphism of $k$, $W(k)$, and $B(k)$. The {\it Honda--Serre--Tate theory} classified the {\it isogeny classes} of {\it abelian varieties} over $k$ (see [Ta2, Thm. 1]) and in particular it proved that each abelian variety over $k$, up to an extension to a finite field extension of $k$ and up to an isogeny, lifts to an abelian scheme with {\it complex multiplication} over a discrete valuation ring of mixed characteristic $(0,p)$ (see [Ta2, Thm. 2]). We recall that an abelian scheme of relative dimension $d$ over an integral scheme is with complex multiplication if its ring of endomorphisms has a commutative $\dbZ$-subalgebra of rank $2d$. Zink generalized [Ta2, Thms. 1 and 2] to contexts that pertain to suitable abelian varieties endowed with endomorphisms, cf. [Zi1, Thms. 4.4 and 4.7]. Special cases of loc. cit. were obtained or announced previously, cf. [Ii1-3], [Lan], and [Mi1]. To detail these contexts and to prepare the background for the present work,  we will use the language of {\it reductive group schemes} and of {\it crystalline cohomology}. 

We recall that a group scheme $F$ over an affine scheme $\Spec R$ is called reductive if it is smooth and affine and its fibres are connected and have trivial unipotent radicals. We denote by $F^{\der}$ and $F^{\ad}$ the derived group scheme and the adjoint group scheme (respectively) of $F$. If $S$ is a closed subgroup scheme of $F$, let $\Lie(S)$ be its {\it Lie algebra} over $R$. If $R$ is a $Q$-algebra, then the pull-back of an object or a morphism $\dag$ or $\dag_{Q}$ (resp. $\dag_*$ with $*$ an index) of the category of $\Spec Q$-schemes to $\Spec R$ is denoted by $\dag_R$ (resp. $\dag_{*,R}$). If $Q\to R$ is a finite, flat monomorphism, let $\Res_{R/Q} S$ be the group scheme over $Q$ obtained from $S$ through the Weil restriction of scalars (see [BT, Subsect. 1.5] and [BLR, Ch. 7, Sect. 7.6]). If moreover $R$ is an \'etale $Q$-algebra, then  $\Res_{R/Q} F$ is a reductive group scheme over $Q$. For a free $R$-module $O$ of finite rank let $\pmb{\GL}_O$  (resp. $\pmb{\text{SL}}_O$) be the reductive group scheme over $R$ of linear automorphisms (resp. of linear automorphisms of determinant $1$) of $O$. If $f_1$ and $f_2$ are two $\dbZ$-endomorphisms of $O$, let $f_1f_2:=f_1\circ f_2$. Let $\bar K$ be an algebraic closure of a field $K$.

\bigskip\noindent 
{\bf 1.1. Isogeny classes.} Let $D$ be a $p$-divisible group over $k$. Let $(M,\phi)$ be the (contravariant) {\it Dieudonn\'e module} of $D$. Thus $M$ is a free $W(k)$-module of finite rank and $\phi:M\to M$ is a $\sigma$-linear endomorphism such that we have an inclusion $pM\subseteq\phi(M)$. We denote also by $\phi$ the $\sigma$-linear automorphism of $\End(M[{1\over p}])$ that maps $e\in\End(M[{1\over p}])$ to $\phi(e):=\phi\circ e\circ\phi^{-1}\in\End(M[{1\over p}])$. 
Let $\scrG$ be a reductive, closed subgroup scheme of $\pmb{\GL}_M$. We recall from [Va7,8] that the triple 
$$\scrC:=(M,\phi,\scrG)$$ 
is called a {\it Shimura $F$-crystal} over $k$ if there exists a direct sum decomposition $M=F^1\oplus F^0$ such that the following two axioms hold:

\medskip
{\bf (i)} we have identities $\phi(M+{1\over p}F^1)=M$ and $\phi(\Lie(\scrG_{B(k)}))=\Lie(\scrG_{B(k)})$, and 

\smallskip
{\bf (ii)} the {\it cocharacter} $\mu:\dbG_m\to \pmb{\GL}_M$ that acts trivially on $F^0$ and as the inverse of the identical character of $\dbG_m$ on $F^1$ (i.e., with weight $-1$ on $F^1$), factors through $\scrG$. 

\medskip
Until the end we will assume that $\scrC$ is a Shimura $F$-crystal over $k$ and that $M=F^1\oplus F^0$ is a direct sum decomposition for which the axioms (i) and (ii) hold. 

The quadruple $(M,F^1,\phi,\scrG)$ is called a {\it Shimura filtered $F$-crystal} over $k$. Either $(M,F^1,\phi,\scrG)$ or $F^1$ is called a {\it lift} of $\scrC$ (to $W(k)$). By an {\it endomorphism} of $\scrC$ (resp. of $(M,F^1,\phi,\scrG)$) we mean an element $e\in\Lie(\scrG)$ fixed by $\phi$ (resp. fixed by $\phi$ and such that we have an inclusion $e(F^1)\subseteq F^1$). We emphasize that the set of endomorphisms of $\scrC$ (resp. of $(M,F^1,\phi,\scrG)$) is in general only a Lie algebra over $\dbZ_p$ (and not a $\dbZ_p$-algebra).   

\medskip
Let $\grP(\scrC)$ be the set of elements $h\in \pmb{\GL}_M(B(k))$ for which the triple
$$(h(M),\phi,\scrG(h))\leqno (1)$$
is a Shimura $F$-crystal over $k$ that can be extended to a Shimura filtered $F$-crystal $(h(M),\tilde h(F^1[{1\over p}])\cap h(M),\phi,\scrG(h))$ over $k$, where $\tilde h\in \scrG(B(k))$ and where $\scrG(h)$ is the schematic closure of $\scrG_{B(k)}$ in $\pmb{\GL}_{h(M)}$. Let 
$\grI(\scrC):=\grP(\scrC)\cap \scrG(B(k))$. It is easy to see that we have:
$$\grI(\scrC)=\{h\in \scrG(B(k))|\exists u\in \scrG(W(k))\,\,\text{such}\,\,\text{that}\,\, u^{-1}h^{-1}\phi hu\phi^{-1}\in \scrG(W(k))\}.$$ 
\noindent
The reductive group scheme $\scrG(h)$ is isomorphic to $\scrG$ (if $h\notin \grI(\scrC)$, then this follows from [Ti2]). For $i\in\{1,2\}$ let $h_i\in\grI(\scrC)$ and $g_i\in \scrG(h_i)(W(k))$. By an {\it inner isomorphism} between $(h_1(M),g_1\phi,\scrG(h_1))$ and $(h_2(M),g_2\phi,\scrG(h_2))$ we mean an element $g\in \scrG(B(k))$ such that we have $g(h_1(M))=h_2(M)$ and $gg_1\phi=g_2\phi g$ (and therefore $g\scrG(h_1)g^{-1}=\scrG(h_2)$). 

By the isogeny class of $\scrC$ we mean the set $\scrI(\scrC)$ of inner isomorphism classes of Shimura $F$-crystals over $k$ that are of the form $(h(M),\phi,\scrG(h))$ with $h\in \grI(\scrC)$.  Ideally, one would like to describe the set $\scrI(\scrC)$ in a way which allows ``the reading" of different Lie algebras of endomorphisms of ({\it ramified}) {\it lifts} of its representatives. Abstract ramified lifts of $\scrC$ (or of $D$ with respect to $\scrG$) are formalized in Subsection 3.2. In this introduction we will only mention the abelian schemes counterpart of ramified lifts. 

\medskip\noindent
{\bf 1.1.1. Two geometric operations.} Until Subsubsection 1.4.1 we will assume that $D$ is the $p$-divisible group of an abelian variety $A$ over $k$.

By a $\dbZ[{1\over p}]$-isogeny between two abelian schemes $A_1$ and $A_2$ over a given scheme we mean a $\dbQ$--isomorphism between $A_1$ and $A_2$ that induces an isomorphism $A_1[N]\arrowsim A_2[N]$ for all natural numbers $N$ that are relatively prime to $p$. For each $h\in\grP(\scrC)$ there exists a unique abelian variety $A(h)$ over $k$ which is $\dbZ[{1\over p}]$-isogenous to $A$ and such that under this $\dbZ[{1\over p}]$-isogeny the Dieudonn\'e module of $A(h)$ is identified with $(h(M),\phi)$. 

If $h\in\grI(\scrC)$, then we say $A(h)$ is {\it $\scrG$-isogenous} to $A$. In all that follows we study the pair $(A,\scrG)$ only up to the following two operations.

\medskip\noindent
$\grO_1:$ {\it the pull-back of $A$ to a finite field extension of $k$.}

\smallskip\noindent
$\grO_2:$ {\it the replacement of $A$ by an abelian variety $A(h)$ over $k$ that is $\scrG$-isogenous to it.}

\bigskip\noindent
{\bf 1.2. Main Problem.} {\it Up to operations $\grO_1$ and $\grO_2$, find conditions which guarantee that there exists a triple $(V,A_V,\scrG^\prime_V)$, where $V$ is a finite, discrete valuation ring extension of $W(k)$ of residue field $k$, $A_V$ is an abelian scheme over $V$ that lifts $A$, and $\scrG^\prime_V$ is a reductive, closed  subgroup scheme of $\pmb{\GL}_{H^1_{\text{dR}}(A_V/V)}$, such that the following four conditions hold:

\medskip
{\bf (a)} the abelian scheme $A_V$ is with complex multiplication;

\smallskip
{\bf (b)} under the canonical identifications $M/pM=H_{\text{dR}}^1(A/k)=H^1_{\text{dR}}(A/V)\otimes_V k$, the group scheme $\scrG^\prime_V$ lifts $\scrG_k$;

\smallskip
{\bf (c)} under the canonical identification $H^1_{\text{dR}}(A/V)[{1\over p}]=M\otimes_{W(k)} V[{1\over p}]$ (see [BO, Thm. 1.3]), the generic fibre of $\scrG^\prime_V$ is the pull-back to $\Spec V[{1\over p}]$ of $\scrG_{B(k)}$;

\smallskip
{\bf (d)} there exists a cocharacter $\dbG_m\to \scrG^\prime_V$ that acts on $F^1_V$ via the inverse of the identical character of $\dbG_m$ and that fixes $H^1_{\text{dR}}(A_V/V)/F^1_V$, where $F^1_V$ is the direct summand of $H^1_{\text{dR}}(A_V/V)$ which is the Hodge filtration of $A_V$.}
\medskip

If (c) holds, then the group schemes $\scrG^\prime_V$ and $\scrG_V$ are isomorphic (cf. [Ti2]). If only (b) to (d) hold and $V=W(k)$ (resp. and $V\neq W(k)$), then we refer to $A_V$ as a lift of $A$ (resp. as a ramified lift of $A$ to $V$) with respect to $\scrG$.

Let $\gre$ be the $B(k)$-span inside $\End(M[{1\over p}])$ of those endomorphisms of $(M,\phi,\scrG)$ which are crystalline realizations of endomorphisms of $A$. It is the Lie algebra of a unique connected subgroup $\scrE$ of $\scrG_{B(k)}$. The uniqueness of $\scrE$ follows from [Bo, Ch. II, Subsect. 7.1] and the existence of $\scrE$ is a standard application of the fact that the $\dbQ$--algebra of $\dbQ$--endomorphisms of $A$ is semisimple. The triple $(V,A_V,\scrG^\prime_V)$ does not always exist (simple examples can be constructed with $\scrG$ a {\it torus}). The reason for this is: in general the ranks of $\scrE$ and $\scrG_{B(k)}$ are not equal. Thus in order to motivate the Main Problem and to list accurately conditions under which one expects that such triples exist, next we will recall some basic terminology. 

\medskip\noindent
{\bf 1.2.1. Review.} We use the terminology of [De3, Sect. 2] for {\it Hodge cycles} on an abelian scheme $B$ over a reduced $\Spec \dbQ$-scheme $Z$. Thus we write a Hodge cycle $v$ on $B$ as a pair $(v_{\text{dR}},v_{\acute et})$, where $v_{\text{dR}}$ and $v_{\acute et}$ are the {\it de Rham} component and the {\it \'etale} component (respectively) of $v$. As its turn, $v_{\acute et}$ has an $l$-component $v_{\acute et}^l$ for each prime $l\in\dbN^*$. 

For instance, if $Z$ is the spectrum of a subfield $E$ of $\bar{\dbQ}\subseteq\dbC$, then $v_{\acute et}^p$ is a suitable $\Gal(E)$-invariant tensor of the tensor algebra of $H^1_{\acute et}(B_{\bar{\dbQ}},\dbQ_p)\oplus (H^1_{\acute et}(B_{\bar{\dbQ}},\dbQ_p))^{\vee}\oplus\dbQ_p(1)$, where $(H^1_{\acute et}(B_{\bar{\dbQ}},\dbQ_p))^{\vee}$ is the dual vector space of $H^1_{\acute et}(B_{\bar{\dbQ}},\dbQ_p)$ (i.e., it is the tensorization with $\dbQ_p$ of the Tate module of $B_{\bar{\dbQ}}$) and where $\dbQ_p(1)$ is the usual Tate twist. The Betti realization $v_B$ of $v$ corresponds to $v_{\text{dR}}$ (resp. to $v_{\acute et}^l$) via the standard isomorphism that relates the de Rham (resp. the $\dbQ_l$ \'etale) cohomology of $B_{\dbC}$  with the Betti cohomology of the complex manifold $B(\dbC)$ with $\dbQ$--coefficients (see [De3, Sects. 1 and 2]).

A {\it Shimura pair} $(G,\scrX)$ consists of a reductive group $G$ over $\dbQ$ and a $G(\dbR)$-conjugacy class $\scrX$ of homomorphisms $\Res_{\dbC/\dbR} \dbG_m\to G_{\dbR}$ that satisfy Deligne's axioms of [De2, Subsubsect. 2.1.1]: the Hodge $\dbQ$--structure of $\Lie(G)$ defined by any $x\in \scrX$ is of type $\{(-1,1),(0,0),(1,-1)\}$, $\text{Ad}(x(i))$ defines a Cartan involution of $\Lie(G^{\ad}_{\dbR})$, and no simple factor of $G^{\ad}$ becomes compact over $\dbR$. Here $\text{Ad}:G_{\dbR}\to\pmb{\GL}_{\Lie(G_{\dbR}^{\ad})}$ is the adjoint representation. The axioms imply that $\scrX$ has a natural structure of a hermitian symmetric domain, cf. [De2, Cor. 1.1.17]. The most studied Shimura pairs are constructed as follows. Let $W$ be a vector space over $\dbQ$ of even dimension $2d$. Let $\psi$ be a non-degenerate alternative form on $W$. Let $\scrS$ be the set of all monomorphisms $\Res_{\dbC/\dbR} \dbG_m\hookrightarrow \pmb{\text{GSp}}(W\otimes_{\dbQ} {\dbR},\psi)$ that define Hodge $\dbQ$--structures on $W$ of type $\{(-1,0),(0,-1)\}$ and that have either $2\pi i\psi$ or $-2\pi i\psi$ as polarizations. The pair $(\pmb{\text{GSp}}(W,\psi),\scrS)$ is a Shimura pair that defines a {\it Siegel modular variety}, cf. [Mi3, p. 161]. See [De1,2], [Mi4], and [Va1, Subsect. 2.5] for different types of Shimura pairs and for their attached {\it Shimura varieties}. We recall that $(G,\scrX)$ is called of {\it Hodge type}, if it can be embedded into a Shimura pair of the form $(\pmb{\text{GSp}}(W,\psi),\scrS)$. We recall that Shimura varieties of Hodge type are moduli spaces of polarized abelian schemes endowed with Hodge cycles, cf. [De1,2], [Mi4], or [Va1, Subsect. 4.1]. 

In this paragraph we will assume that the adjoint group $G^{\ad}$ is $\dbQ$--simple. Let $\theta$ be the {\it Lie type} of each simple factor of $G^{\ad}_{\dbC}$. If $\theta\in\{A_n,B_n,C_n|n\in\dbN^*\}$, then $(G,\scrX)$ is said to be of $\theta$ type. If $\theta=D_n$ with $n\Ge 4$, then $(G,\scrX)$ is of one of the following three disjoint types: $D_n^{\dbH}$, $D_n^{\dbR}$, and $D_n^{\text{mixed}}$ (cf. [De2] and [Mi4]). If $(G,\scrX)$ is of $D_n^{\dbR}$ (resp. of $D_n^{\dbH}$) type, then all simple, non-compact factors of $G_{\dbR}^{\ad}$ are isomorphic to $\pmb{\text{SO}}(2,2n-2)^{\ad}_{\dbR}$ (resp. to $\pmb{\text{SO}}^*(2n)_{\dbR}^{\ad}$) and the converse of this statement holds for $n\Ge 5$ (see [He, p. 445] for the classical groups $\pmb{\text{SO}}(2,2n-2)^{\ad}_{\dbR}$ and $\pmb{\text{SO}}^*(2n)_{\dbR}^{\ad}$). If moreover $(G,\scrX)$ is of Hodge type, then $(G,\scrX)$ is of one of the following five possible types: $A_n$, $B_n$, $C_n$, $D_n^{\dbH}$, and $D_n^{\dbR}$ (see [Sa1,2] and [De2, Table 2.3.8]).

\medskip\noindent
{\bf 1.2.2. Conjecture.} {\it We assume that one of the following two conditions holds:

\medskip
{\bf (i)} the subgroup $\scrE$ of $\scrG_{B(k)}$ has the same rank as $\scrG_{B(k)}$;

\smallskip
{\bf (ii)} there exists an abelian scheme $A_{W(k)}$ over $W(k)$ which lifts $A$ and for which there exists a family $(w_{\alpha})_{\alpha\in\scrJ}$ of Hodge cycles on its generic fibre $A_{B(k)}$ such that $\scrG_{B(k)}$ is the subgroup of $\pmb{\GL}_{M[{1\over p}]}$ that fixes the crystalline realization $t_{\alpha}$ of $w_{\alpha}$ for all $\alpha\in\scrJ$ (under the canonical identification $M=H^1_{\text{dR}}(A_{W(k)}/W(k))$, the crystalline and the de Rham realizations of $w_{\alpha}$ coincide).

\medskip
Then up to the operations $\grO_1$ and $\grO_2$, there exists a triple $(V,A_V,\scrG^\prime_V)$ such that all conditions 1.2 (a) to (d) hold.}

\medskip
If (i) (resp. (ii)) holds, then we refer to Conjecture 1.2.2 as Conjecture 1.2.2 (i) (resp. Conjecture 1.2.2 (ii)). Conjecture 1.2.2 stems from the {\it Langlands--Rapoport conjecture} (see [LR], [Mi2,3], [Pf], and [Re2]) on $\bar k$-valued points of special fibres of (see [Va1] for precise definitions) {\it integral canonical models} of Shimura varieties of Hodge type in mixed characteristic $(0,p)$. This motivic conjecture of combinatorial nature is a key ingredient in the understanding of zeta functions of Shimura varieties of Hodge type and of different trace functions that pertain to $\dbQ_l$-local systems on quotients of finite type of such integral canonical models (for instance, see [LR], [Ko2], and [Mi3]; here $l$ is a prime different from $p$). Conjecture 1.2.2 (ii) is in fact only a slight refinement of an adequate translation of a part of the Langlands--Rapoport conjecture. The Langlands--Rapoport conjecture is known to be true for Siegel modular varieties (see [Mi2]) and for certain Shimura varieties of $A_1$ type (see [Ii2,3] and [Re1]). 

We added Conjecture 1.2.2 (i) due to the following two reasons. First, if one assumes the standard Hodge conjecture for complex abelian varieties (see [Le, Ch. 7]), then (ii)$\Rightarrow$(i). 

Second, often due to technical reasons one assumes that $\scrG^{\der}$ is simply connected and this excludes the cases related to Shimura pairs of $D_n^{\dbH}$ type (see [De2, Rm. 1.3.10 (ii)]). Thus to handle Conjecture 1.2.2 (ii) in cases related to Shimura pairs of $D_n^{\dbH}$ type, one can proceed in two steps as follows: (a) first solve Conjecture 1.2.2 (ii) in cases related to Shimura pairs of $C_{2n}$ type and then (b) appeal to {\it relative PEL situations} defined in [Va1, Subsubsect. 4.3.16] and constructed as in [Va4, Rm. 4.8.2 (b)], in order to reduce Conjecture 1.2.2 (ii) to Conjecture 1.2.2 (i) for these cases related to the $D_n^{\dbH}$ type. 

In what follows we will also refer to the following Subproblem of the Main Problem.

\medskip\noindent
{\bf 1.2.3. Subproblem.} {\it Same as Main Problem but with condition 1.2 (a) replaced by the weaker condition that the $p$-divisible group of $A_V$ is with complex multiplication (i.e., the image of the $p$-adic Galois representation associated to the Tate module $T_p(A_{V[{1\over p}]})$ of $A_{V[{1\over p}]}$, is formed by semisimple elements that commute).} 

\medskip\noindent
{\bf 1.2.4. Definitions.} {\bf (a)} Let $n\in\dbN^*$. We say $\scrC$ is of $B_n$ and $D_n^{\dbR}$ type if the following three conditions hold:

\medskip\noindent
{\bf (a.i)} each simple factor $\scrV$ of $\scrG^{\ad}_{W(\bar k)}$ is of either $B_n$ or $D_n$ Lie type;

\smallskip\noindent
{\bf (a.ii)} if a simple factor $\scrV$ is of $D_n$ Lie type, then the centralizer in $\scrV$ of the image of $\mu_{W(\bar k)}$ in $\scrV$ is either $\scrV$ itself or is a reductive group scheme whose adjoint is of $D_{n-1}$ Lie type;

\smallskip\noindent
{\bf (a.iii)} if $n=4$ and if a simple factor $\scrV$ is of $D_4$ Lie type, then the images in $\scrV$ of iterates of $\mu_{W(\bar k)}$ under integral powers of $\phi\otimes\sigma_{\bar k}$, are all $\scrV(W(\bar k))$-conjugate.

\medskip
{\bf (b)} By the standard axioms for $\scrC$ we mean the following two axioms:

\medskip\noindent
{\bf (b.i)} there exists a family $(t_{\alpha})_{\alpha\in\scrJ}$ of tensors of the tensor algebra of $M\oplus \Hom_{W(k)}(M,W(k))$ fixed by $\phi$ and $\scrG$ and such that $\scrG$ is the schematic closure in $\pmb{\GL}_M$ of the subgroup of $\pmb{\GL}_{M[{1\over p}]}$ that fixes $t_{\alpha}$ for all $\alpha\in\scrJ$;

\smallskip\noindent
{\bf (b.ii)} there exists a set of cocharacters of $\scrG_{W(\bar k)}$ which act on $M\otimes_{W(k)} W(\bar k)$ via the weights $0$ and $-1$ and whose images in $\scrG^{\ad}_{W(\bar k)}$ generate $\scrG^{\ad}_{W(\bar k)}$.

\medskip
Part (a) (resp. (b)) of the following theorem solves an unramified version of (resp. solves a refined version of) the Subproblem 1.2.3 for the case when $\scrC$ is of $B_n$ and $D_n^{\dbR}$ type. 

\medskip\noindent
{\bf 1.2.5. Theorem.} {\it We assume that the standard axioms hold for $\scrC$ and that there exists a polarization $\lambda_A$ of $A$ whose crystalline realization is an alternating form $M\times M\to W(k)$ whose $W(k)$-span is normalized by $\scrG$. We also assume that $\scrC$ is of $B_n$ and $D_n^{\dbR}$ type. 

\medskip
{\bf (a)} Then up to operations $\grO_1$ and $\grO_2$, there exists an abelian scheme $A_{W(k)}$ over $W(k)$ which is a lift of $A$ with respect to $\scrG$ (i.e., conditions 1.2 (b) to 1.2 (d) hold with $(V,\scrG_V^\prime)=(W(k),\scrG)$) and  for which the $p$-divisible group of $A_{W(k)}$ is with complex multiplication and to which the Frobenius endomorphism of $A$ lifts.

\smallskip
{\bf (b)} We also assume that $p>2$. Let $\scrT_{B(k)}$ be an arbitrary maximal torus of $\scrG_{B(k)}$ whose Lie algebra $\Lie(\scrT_{B(k)})$ is generated by crystalline realizations of $\dbQ_p$--endomorphisms of $D$. Then up to operations $\grO_1$ and $\grO_2$, there exists a triple $(V,A_V,\scrG^\prime_V)$ as in the Main Problem 1.2 for which the conditions 1.2 (b) to 1.2 (d) hold and for which the $p$-divisible group of $A_V$ is with complex multiplication and in fact each element of $\Lie(\scrT_{B(k)})$ fixed by $\phi$ is the crystalline realization of a $\dbQ_p$--endomorphism of $A_V$.}

\medskip
Theorem 1.2.5 is proved in Subsection 7.4.

\bigskip\noindent
{\bf 1.3. The classical PEL context.} This is the context in which there exists a principal polarization $\lambda_A$ of $A$ and there exists a $\dbZ_{(p)}$-subalgebra $\Theta$ of $\End(M)$ formed by crystalline realizations of $\dbZ_{(p)}$-endomorphisms of $A$, such that the following two conditions hold:

\medskip
{\bf (i)} the $W(k)$-algebra $\Theta\otimes_{\dbZ_{(p)}} W(k)$ is semisimple, is self dual with respect to the perfect alternating form $\lambda_A:M\times M\to W(k)$ which is the crystalline realization of $\lambda_A$ (and which is denoted in the same way), and is equal to the following $W(k)$-algebra $\{e\in\End(M)|e\,\,\text{fixed}\,\,\text{by}\,\,\scrG\}$;

\smallskip
{\bf (ii)}  the group $\scrG_{B(k)}$ is the identity component of the subgroup $\scrD(\lambda_A)_{B(k)}$ of $\break\pmb{\text{GSp}}(M[{1\over p}],\lambda_A)$ that fixes each element of $\Theta[{1\over p}]$.

\medskip
As $A$ is with complex multiplication and the $\dbQ$--algebra $\End(A)\otimes_{\dbZ} \dbQ$ is semisimple, in this context the group $\scrE$ is reductive and has the same rank as $\scrG_{B(k)}$; thus the condition 1.2.2 (i) holds. The existence up to operations $\grO_1$ and $\grO_2$ of a triple $(V,A_V,\scrG^\prime_V)$ such that all the conditions 1.2 (a) to (d) hold was proved (using a slightly different language) in  [Zi1, Thm. 4.4] for the cases when $\scrG_{B(k)}=\scrD(\lambda_A)_{B(k)}$ (strictly speaking, loc. cit. assumes that $\Theta[{1\over p}]$ is a $\dbQ$--simple algebra; but the case when $\Theta[{1\over p}]$ is not $\dbQ$--simple gets easily reduced to the case when it is so). For the mentioned cases, loc. cit. also shows that (even if $p=2$) we can choose the triple $(V,A_V,\scrG^\prime_V)$ such that $\Theta$ lifts to a family of $\dbZ_{(p)}$-endomorphisms of $A_V$ and that $\lambda_A$ is the crystalline realization of a principal polarization of $A_V$. Some refinements of loc. cit., which are still weaker than the Langlands--Rapoport conjecture for the corresponding Shimura varieties of PEL type and which also consider for $p>2$ the case when $\scrG_{B(k)}$ is the identity component of the subgroup $\scrD(\lambda_A)_{B(k)}$, were obtained in [ReZ] and [Ko2].

\bigskip\noindent
{\bf 1.4. On results and tools.} The goal of the paper is to solve Conjecture 1.2.2 (i) and Subproblem 1.2.3 in contexts general enough (see Theorem 1.2.5, Corollary 7.3, Variant 7.5, and Subsections 8.2 to 8.4) which suffice for the following main two purposes:

\medskip
{\bf (a)} to classify (at least for $p>2$) all maximal tori of $\scrG_{B(k)}$ whose Lie algebras are generated by crystalline realizations of $\dbQ$--endomorphisms of $A$ that lift to some $A_V$, where $A_V$ is part of a triple $(V,A_V,\scrG_V^\prime)$ for which the conditions 1.2 (a) to (d) hold;

\smallskip
{\bf (b)} that works in progress (like [Mi5] and [Va10]) can be plugged in to result in complete proofs of {\it refined forms} of Conjecture 1.2.2 (ii) and of the Langlands--Rapoport conjecture for Shimura varieties of Hodge type which involve simply connected derived groups (here the word refined refers to solutions that: (b.i) also accomplish (a) and (b.ii) allow us to take $V=W(k)$ in Conjecture 1.2.2 (ii)). 

\medskip
The passage from the mentioned solutions to a refined solution of Conjecture 1.2.2 (ii) for the case when $p\ge 5$ and $\scrG^{\der}$ is simply connected, is controlled by [Va1,2,3,6,7] and by the following two extra things (see already [Va10] which in fact works for all $p\ge 2$):

\medskip
{\bf (i)} the {\it weak isogeny property} which says that each {\it rational stratification} as in [Va7, Subsect. 5.3] has only one stratum that has a closed connected component;

\smallskip
{\bf (ii)} recent work of Milne (see [Mi5]).

\medskip
It is well known that the weak isogeny property holds for Siegel modular varieties. For instance, this can be easily proved based on the following fact (see [Oo]): the {\it Newton polygon} stratification of the {\it Mumford moduli scheme} $\scrA_{d,1,N}$ over $\bar k$ has only one closed stratum (the supersingular one); here $d,N\in\dbN^*$, $N\Ge 3$, and $g.c.d.(N,p)=1$. The weak isogeny property for arbitrary rational stratifications requires methods different from the ones of this paper and thus it will be proved in future works (see already [Va10] for Shimura varieties which have compact factors in the sense of [Va5, Subsect. 2.2]). 

\smallskip
The main tools we use in this paper are the following seven.

\medskip
{\bf T1.} The {\it rational classification} of Shimura $F$-crystals over $\bar k$ accomplished in [Va7]. 

\smallskip
{\bf T2.} Approximations of tori of reductive groups over $\dbQ$, cf. [Ha, Lem. 5.5.3]. 

\smallskip
{\bf T3.} A new theory of {\it admissible cocharacters} of extensions of maximal tori of $\scrG_{B(k)}$ contained in tori of $\pmb{\GL}_{M[{1\over p}]}$ whose Lie algebras are $B(k)$-generated by crystalline realizations of $\dbQ_p$--endomorphisms of $A$. In its abstract form, the theory refines [RaZ, Subsects. 1.21 to 1.25] for Shimura $F$-crystals in two ways. First, it is over $k$ and not only over $\bar k$. Second, in many cases it works without assuming that all Newton polygon slopes of $(\Lie(\scrG_{B(k)}),\phi)$ are $0$ and moreover it applies to {\it all} such maximal tori of $\scrG_{B(k)}$. We mention that in connection to either this theory or loc. cit., [FR] does not provide new tools. 

\smallskip
{\bf T4.} In some cases we rely on [Zi1, Thm. 4.4] (see proof of Theorem 8.4, etc.).

\smallskip
{\bf T5.} The classification of isogeny classes of $p$-divisible groups over $p$-adic fields obtained in [Br, Subsect. 5.3] for $p\Ge 3$ and in [Ki1] for $p=2$. 

\smallskip
{\bf T6.} The natural $\dbZ_p$ structure $\scrG_{\dbZ_p}$ of $\scrG$ defined by $\scrC$ via the theory of {\it canonical split cocharacters} of [Wi] and the vanishing of certain classes in the pointed set $H^1(\dbQ_p,\scrG_{\dbQ_p})$. 

\smallskip
{\bf T7.} The recent result [Ki2, Cor. (1.4.3)] which can substitute our theory of {\it well positioned families of tensors} developed in [Va1, Subsect. 4.3] (here one has to assume that either $p\ge 3$ or $p=2$ and the $p$-divisible group of $A$ is connected). 

\medskip
See [Fo] for (weakly admissible or admissible) filtered modules over $p$-adic fields. Next we exemplify how the tools T1 to T7 work under some conditions. We emphasize that often we do have to perform either the operation $\grO_1$ or the operation $\grO_2$ but this will not be repeated in this paragraph. Based on [Va7, Thm. 3.1.2 (b) and (c)], in connection to Conjecture 1.2.2 (i) and to Subproblem 1.2.3 it suffices to refer to the case when all Newton polygon slopes of $(\Lie(\scrG_{B(k)}),\phi)$ are $0$ (see Subsection 2.6). Assuming that the condition 1.2.2 (i) holds, we show based on [Ha, Lem. 5.5.3] that there exist maximal tori of $\scrG_{B(k)}$ as mentioned in the tool T3 but with $\dbQ_p$ replaced by $\dbQ$. The essence of T3 can be described as follows. For an {\it arbitrary} maximal torus of $\scrG_{B(k)}$ as in the tool $T3$, we show the existence of suitable cocharacters of its extension to a finite field extension $V[{1\over p}]$ of $B(k)$ such that the resulting filtered modules over $V[{1\over p}]$ are weakly admissible. For this, in some cases related to Shimura pairs of $A_n$, $C_n$, and $D_n^{\dbH}$ types we rely as well on [Zi1, Thm. 4.4] and accordingly some extra assumptions are imposed (roughly speaking we deal with abelian varieties associated to Shimura varieties of Hodge type constructed as in [De2, Prop. 2.3.10] but in the integral contexts of [Va1, Sects. 5 and 6]). Using the tool T5, we get an isogeny class of $p$-divisible groups over $V$. Using the tool T6 we get natural choices of representatives of this isogeny class so that we end up in the \'etale context with a  reductive group scheme over $\dbZ_p$ whose generic fibre corresponds via Fontaine comparison theory to $\scrG_{V[{1\over p}]}$. Using the tool T7, if $p\ge 3$ we ``transfer backwards" (as in [Va1, Subsects. 5.2 and 5.3] and [Ki2, Cor. (1.4.3) (2)]) the mentioned reductive group scheme over $\dbZ_p$ in order to end up again with a reductive group scheme $\scrG^{\prime}_V$ in the de Rham context over $V$. 
 
\medskip\noindent
{\bf 1.4.1. On contents.} Motivated by general applications, in Sections 2 to 6 we work abstractly. Thus we work with an arbitrary Shimura $F$-crystal $\scrC$ over $k$ and, even if $D$ is the $p$-divisible group of some abelian variety $A$ over $k$, most often we do not impose any geometric condition on the group scheme $\scrG$ over $W(k)$ (of the type of conditions 1.2.2 (i) and (ii) or 1.3 (i) and (ii)). In Section 2 we develop a minute language that pertains to Subsection 1.1 and to the tool T3 which allows us to solve in many cases refined versions of Conjecture 1.2.2 (i) and Subproblem 1.2.3. Different abstract $\text{CM}$-isogeny classifications are formalized in Section 3. In particular, Corollary 3.6.3 shows that if $p\Ge 3$, then the ramified lifts of $D$ with respect to $\scrG$ (see Definition 3.6.2) are in bijection to the ramified lifts of $\scrC$ (see Definitions 3.2.1). This is a refined version of the classification of $p$-divisible groups over $V$ achieved for $p\ge 3$ by Faltings, Breuil, and Zink (see [Fa], [Br], and [Zi2]). 

In Section 4 we state in the abstract context two basic results that pertain to the tool T3 (see Basic Theorems 4.1 and 4.2) and a Corollary 4.3 which is the very essence behind Theorem 1.2.5 (b). Corollary 4.3 presents the very first situations of general nature where complete ramified $\text{CM}$-classifications as defined in Subsubsection 3.2.3 are accomplished; to ``balance" the focus of [Zi1] on Shimura varieties of PEL type (and thus of either $A_n$ or totally non-compact $C_n$ or $D_n^{\dbH}$  type), they involve cases of $B_n$ and $D_n^{\dbR}$ type. Sections 5 and 6 prove the Basic Theorems 4.1 and 4.2 (respectively).

First applications to abelian varieties are included in Section 7 (see Subsections 7.3 to 7.5 for partial solutions to Conjecture 1.2.2 (i) and Subproblem 1.2.3). Section 8 introduces the integral context of moduli spaces of polarized abelian varieties endowed with (specializations of) Hodge cycles. See Subsections 8.2 to 8.4 for different properties and how they lead to generalizations of the results of Zink recalled in Subsection 1.3.

\bigskip\smallskip
\noindent
{\boldsectionfont 2. Preliminaries}
\bigskip

Subsection 2.1 lists notations and conventions. Subsection 2.2 recalls descent properties of connected, affine, algebraic groups in characteristic $0$. Lemma 2.3 pertains to reductive group schemes over $\dbZ_p$. Subsections 2.4 and 2.5 mainly introduce a language. Subsection 2.6 shows that often one can assume that all Newton polygon slopes of $(\Lie(\scrG_{B(k)}),\phi)$ are $0$. Subsection 2.7 proves a general variant of Theorem 1.2.5 (a) in the so called case with many endomorphisms. Subsection 2.8 introduces $W(k)$-algebras that are required for the ramified contexts of Sections 3 to 6.

\bigskip\noindent
{\bf 2.1. Notations and conventions.} Let $R$, $F$, and $O$  be as before Subsection 1.1. We refer to [Va7, Subsect. 2.2] for quasi-cocharacters of $F$. Let $Z(F)$ be the center of $F$; we have $F^{\ad}=F/Z(F)$. Let $Z^0(F)$ be the maximal torus of $Z(F)$; the quotient group scheme $Z(F)/Z^0(F)$ is a finite, flat group scheme over $R$ of multiplicative type. Let $F^{\ab}:=F/F^{\der}$; it is the maximal abelian quotient of $F$. Let $F^{\sc}$ be the simply connected semisimple group scheme cover of $F^{\der}$. If $S$ is a reductive, closed subgroup scheme of $F$, let $C_F(S)$ (resp. $N_F(S)$) be the centralizer (resp. the normalizer) of $S$ in $F$. Thus $C_F(S)$ (resp. $N_F(S)$) is a closed subgroup scheme of $F$, cf. [DG, Vol. II, Exp. XI, Cor. 6.11]. If $S$ is a torus, then $C_F(S)$ is a reductive group scheme (cf. [DG, Vol. III, Exp. XIX, Subsect. 2.8 and Prop. 6.3]) and $N_F(S)$ is a smooth closed subgroup scheme of $F$ whose identity component is $C_F(S)$ (cf. [DG, Vol. III, Exp. XXII, Cors. 5.3.10 and 5.3.18 (ii)]). If $R$ is a finite, discrete valuation ring extension of $W(k)$, then $F(R)$ is called a {\it hyperspecial} subgroup of $F(R[{1\over p}])$ (see [Ti2]). Let $O^{\vee}:=\Hom_{R}(O,R)$. A bilinear form on $O$ is called perfect if it induces naturally an isomorphism
$O\arrowsim O^{\vee}$. We will often use the following free $O$-module
$$
\scrT(O):=\oplus_{s,t\in\dbN} O^{\otimes s}\otimes_{R} O^{{\vee}\otimes t}.
$$ 
We use the same notation for two perfect bilinear forms or tensors of two tensor algebras if they are obtained one from another via either a reduction modulo some ideal or a scalar extension. If $F^1(O)$ is a direct summand of $O$, then $F^0(O^{\vee}):=(O/F^1(O))^{\vee}$ is a direct summand of $O^{\vee}$. By the $F^0$-filtration of $\scrT(O)$ defined by $F^1(O)$ we mean the direct summand of $\scrT(O)$ whose elements have filtration degrees at most $0$, where $\scrT(O)$ is equipped with the tensor product filtration defined by the decreasing, exhaustive, and separated filtrations $(F^i(O))_{i\in\dbZ}$ and $(F^i(O^{\vee}))_{i\in\dbZ}$ of $O$ and $O^{\vee}$ (respectively). Here $F^0(O):=O$, $F^2(O):=0$, $F^{-1}(O^{\vee}):=O^{\vee}$, $F^0(O^{\vee}):=\{x\in O^{\vee}|x(F^1(O))=0\}$, and $F^1(O^{\vee}):=0$. We always identify $\End(O)$ with $O\otimes_{R} O^{\vee}$. Thus $\End(\End(O))=\End(O\otimes_{R} O^{\vee})=O\otimes_{R} O^{\vee}\otimes_{R} O^{\vee}\otimes_{R} O$ is always identified by interchanging the second and fourth factors  with the direct summand $O^{\otimes 2}\otimes_{R} O^{{\vee}\otimes 2}$ of $\scrT(O)$.

Let $x\in R$ be a non-divisor of $0$. A family of tensors of $\scrT(O[{1\over x}])=\scrT(O)[{1\over x}]$ is denoted $(u_{\alpha})_{\alpha\in\scrJ}$, with $\scrJ$ as the set of indexes. Let $O_1$ be another free $O$-module of finite rank. Let $(u_{1,\alpha})_{\alpha\in\scrJ}$ be a family of tensors of $\scrT(O_1[{1\over x}])$ indexed also by the set $\scrJ$. By an isomorphism $(O,(u_{\alpha})_{\alpha\in\scrJ})\arrowsim (O_1,(u_{1,\alpha})_{\alpha\in\scrJ})$ we mean an $R$-linear isomorphism $O\arrowsim O_1$ that extends naturally to an $R[{1\over x}]$-linear isomorphism $\scrT(O[{1\over x}])\arrowsim\scrT(O_1[{1\over x}])$ which takes $u_{\alpha}$ to $u_{1,\alpha}$ for all $\alpha\in\scrJ$.

If $K$ is a $p$-adic field, see [Fo] for de the Rham ring $B_{\text{dR}}(K)$ and for admissible Galois representations of the Galois group $\Gal(K):=\Gal(\bar K/K)$. For the classification of Lie and Dynkin types we refer to [Bou1] and [DG, Vol. III, Exp. XXII and XXIII]. Whenever we use a $D_n$ type, we assume $n\Ge 3$. Let $\dbZ_{(p)}$ be the localization of $\dbZ$ at its prime ideal $(p)$. 

By a Frobenius lift of a flat $\dbZ_{(p)}$-algebra $R$ we mean an endomorphism $\Phi_R:R\to R$ which modulo $p$ is the usual Frobenius endomorphism of $R/pR$. If $\phi_O:O\to O$ is a $\Phi_R$-linear endomorphism such that $O[{1\over p}]$ is $R[{1\over p}]$-generated by $\phi_O(O)$, then we denote also by $\phi_O$ the $\Phi_R$-linear endomorphism of each $R$-submodule of $\scrT(O)[{1\over p}]$ left invariant by $\phi_O$. We recall that $\phi_O$ acts on $O^{\vee}[{1\over p}]$ via the rule: if $f\in O^{\vee}[{1\over p}]$ and $e\in O[{1\over p}]$, then $\phi_O(f)(\phi_O(e))=\Phi_R(f(e))\in R[{1\over p}]$. If $\phi_O$ becomes an isomorphism after inverting $p$ and if $\mu_O$ is a cocharacter of $\pmb{\GL}_{O[{1\over p}]}$, then $\phi_O(\mu_R):=\phi_O\mu_O\phi_O^{-1}$ is a cocharacter of $\pmb{\GL}_{O[{1\over p}]}$. 

Always $\scrC:=(M,\phi,\scrG)$ (resp. $(M,F^1,\phi,\scrG)$) is a Shimura (resp. Shimura filtered) $F$-crystal over $k=\dbF_q$ and $D$ is a $p$-divisible group over $k$ whose Dieudonn\'e module is $(M,\phi)$. We fix a cocharacter $\mu:\dbG_m\to\scrG$ of $\scrC$ as in Subsection 1.1 (thus it normalizes $F^1$); we call it a {\it Hodge cocharacter} of $\scrC$ and we say that it defines $F^1$. Let $\scrP$ be the parabolic subgroup scheme of $\scrG$ which is the normalizer of $F^1$ in $\scrG$, cf. [Va7, paragraph before Cor. 2.3.2]. Let the sets $\grI(\scrC)$, $\grP(\scrC)$, and $\scrI(\scrC)$ be as in Subsection 1.1. Let
$$C:=C_{\pmb{\GL}_M}(\scrG).$$
If $C$ is a reductive group scheme over $W(k)$, then let $C_1:=C_{\pmb{\GL}_M}(C)$. 

See [Va7, Subsubsects. 2.2.1 and 2.2.3] for the Newton quasi-cocharacter of $\scrC$. Let $P_{\scrG}^+(\phi)$, $P_{\scrG}^-(\phi)$, and $L_{\scrG}^0(\phi)_{B(k)}$ be the non-negative parabolic subgroup scheme, the non-positive parabolic subgroup scheme, and the Levi subgroup (respectively) of $\scrC$ we defined in [Va7, Lem. 2.3.1 and Def. 2.3.3]. Thus $P_{\scrG}^+(\phi)$ is the parabolic subgroup scheme of $\scrG$ which is maximal subject to the property that $\Lie(P_{\scrG}^+(\phi)_{B(k)})$ is normalized by $\phi$ and all Newton polygon slopes of $(\Lie(P_{\scrG}^+(\phi)_{B(k)}),\phi)$ are non-negative, $P_{\scrG}^-(\phi)$ is defined similarly but by replacing non-negative with non-positive, and $L_{\scrG}^0(\phi)_{B(k)}$ is the unique Levi subgroup of either $P_{\scrG}^+(\phi)_{B(k)}$ or $P_{\scrG}^-(\phi)_{B(k)}$ with the property that $\Lie(L_{\scrG}^0(\phi)_{B(k)})$ is normalized by $\phi$ and all Newton polygon slopes of $(\Lie(L_{\scrG}^0(\phi)_{B(k)}),\phi)$ are $0$. We have $P_{\scrG}^+(\phi)_{B(k)}\cap P_{\scrG}^-(\phi)_{B(k)}=L_{\scrG}^0(\phi)_{B(k)}$. Let $U_{\scrG}^+(\phi)$ be the unipotent radical of $P_{\scrG}^+(\phi)$. Let $L_{\scrG}^0(\phi)$ be the schematic closure of $L_{\scrG}^0(\phi)_{B(k)}$ in $\scrG$ (or in $P_{\scrG}^+(\phi)$); we emphasize that it is not always a Levi subgroup scheme of $P_{\scrG}^+(\phi)$. We say $\scrC$ is {\it basic} if all Newton polygon slopes of $(\Lie(\scrG_{B(k)}),\phi)$ are $0$ (i.e., if $P_{\scrG}^+(\phi)=P_{\scrG}^-(\phi)=L_{\scrG}^0(\phi)=\scrG$). 

Always $k_1$ is a finite field extension of $k$. For a perfect field $l$ containing $k$, let $W(l)$, $B(l)$, and $\sigma_l$ be the analogues of $W(k)$, $B(k)$, and $\sigma_k$ but for $l$ instead of $k$. Let 
$$\scrC\otimes l:=(M\otimes_{W(k)} W(l),\phi\otimes\sigma_l,\scrG_{W(l)})$$
be the extension of $\scrC$ to $l$. We also refer as $\grO_1$ to the operation of replacing $k$ by $k_1$ and $\scrC$ by $\scrC\otimes k_1$, and as $\grO_2$ to the operation of replacing $\scrC$ by $(h(M),\phi,\scrG(h))$, where $h\in \grI(\scrC)$ and $\scrG(h)$ are as in Subsection 1.1. For $g\in\scrG(W(k))$ let 
$$\scrC_g:=(M,g\phi,\scrG).$$ 
We have $\scrC=\scrC_{1_M}$. Let $\scrF:=\{\scrC_g|g\in G(W(k))\}$ be the {\it family} of Shimura $F$-crystals over $k$ associated naturally to $\scrC$. Let $\scrY(\scrF):=\cup_{g\in \scrG(W(k))} \scrI(\scrC_g)$. The (inner) isomorphism class of some object $\natural$ will be denoted as $[\natural]$.

 Though in this paper we deal only with Shimura $F$-crystals, Sections 2 to 4 are organized in such a way that the interested reader can extend their notions to the context of $p$-divisible objects with a reductive group over $k$ introduced in [Va7, Def. 1.2.1] (often even over an arbitrary perfect field of characteristic $p$). 

\bigskip\noindent
{\bf 2.2. Lemma.} {\it Let $\eta\subseteq\eta_1$ be an extension of fields of characteristic $0$. Let $\grG$ be a connected, affine, algebraic group over $\eta$. Let $\grL$ be a Lie subalgebra of $\Lie(\grG)$. We assume that there exists a connected (resp. reductive) subgroup $\grS_{\eta_1}$ of $\grG_{\eta_1}$ whose Lie algebra is $\grL\otimes_{\eta} \eta_1$. Then the following two properties hold:

\medskip
{\bf (a)} there exists a unique connected (resp. reductive)  subgroup $\grS$ of $\grG$ whose Lie algebra is $\grL$ (the notations match i.e., the extension of $\grS$ to $\eta_1$ is $\grS_{\eta_1}$);

\smallskip
{\bf (b)} if $\grS$ is a reductive group and if $\grG$ is the general linear group $\pmb{\GL}_W$ of  a finite dimensional $\eta$-vector space $W$, then the restriction of the trace form on $\End(W)$ to $\grL$ is non-degenerate.}

\medskip
\proof 
We prove (a). The uniqueness part is implied by [Bo, Ch. I, Subsect. 7.1]. Loc cit. also implies that if $\grS$ exists, then its extension to $\eta_1$ is indeed $\grS_{\eta_1}$. It suffices to prove (a) for the case when $\grS$ is connected. We consider commutative $\eta$-algebras $\kappa$ for which there exists a closed subgroup scheme $\grS_{\kappa}$ of $\grG_{\kappa}$ whose Lie algebra is $\grL\otimes_{\eta} \kappa$. Our hypotheses imply that as $\kappa$ we can take $\eta_1$. Thus as $\kappa$ we can also take a finitely generated $\eta$-subalgebra of $\eta_1$. By considering the reduction modulo a maximal ideal of this last $\eta$-algebra, we can assume that $\kappa$ is a finite field extension of $\eta$. Even more, (as $\eta$ has characteristic $0$) we can assume that $\kappa$ is a finite Galois extension of $\eta$. By replacing $\grS_{\kappa}$ with its identity component, we can assume that $\grS_{\kappa}$ is connected. Due to the mentioned uniqueness part, the Galois group $\Gal(\kappa/\eta)$ acts naturally on the connected subgroup $\grS_{\kappa}$ of $\grG_{\kappa}$. As $\grS_{\kappa}$ is an affine scheme, the resulting Galois descent datum on $\grS_{\kappa}$ with respect to $\Gal(\kappa/\eta)$ is effective (cf. [BLR, Ch. 6, Sect. 6.1, Thm. 5]). This implies the existence of a subgroup $\grS$ of $\grG$ whose extension to $\kappa$ is $\grS_{\kappa}$. As $\Lie(\grS)\otimes_{\eta} \kappa=\Lie(\grS_{\kappa})=\grL\otimes_{\eta} \kappa$, we have $\Lie(\grS)=\grL$. The group $\grS$ is connected as $\grS_{\kappa}$ is so. Thus $\grS$ exists i.e., (a) holds. See [Va6, Prop. 3.2] for another approach to prove (a). 

To check (b) we can assume that $\eta$ is algebraically closed. Using isogenies, it suffices to prove (b) in the case when $\grS$ is either $\dbG_m$ or a semisimple group whose adjoint is simple. If $\grS$ is $\dbG_m$, then the $\grS$-module $W$ is a direct sum of one dimensional $\grS$-modules. We easily get that there exists an element $x\in\grL\setminus\{0\}$ which is a semisimple element of $\End(W)$ whose eigenvalues are integers. The trace of $x^2$ is a non-trivial sum of squares of non-zero integers and thus it is non-zero. Thus (b) holds if $\grS$ is $\dbG_m$. If $\grS$ is a semisimple group whose adjoint is simple, then $\grL$ is a simple Lie algebra over $\eta$. From Cartan solvability criterion we get that the restriction of the trace form on $\End(W)$ to $\grL$ is non-zero and therefore (as $\grL$ is a simple Lie algebra) it is non-degenerate. Thus (b) holds.\endproof

\medskip\noindent
{\bf 2.2.1. Example.} We take $\eta=\dbQ_p$ and $\eta_1=B(\star)$, where $\star$ is a perfect field of characteristic $p$. Let $(\scrW,\varphi)$ be an $F$-crystal over $\star$; thus $\scrW$ is a free $W(\star)$-module of finite rank and $\varphi$ is a $\sigma_{\star}$-linear monomorphism of $\scrW$. Let $\grG$ be the group over $\dbQ_p$ which is the group scheme of invertible elements of the $\dbQ_p$-algebra $\{e\in\End(\scrW)|\varphi(e)=e\}$. Let $\diamond$ be a connected subgroup of $\grG_{\eta_1}$ whose Lie algebra is $\eta_1$-generated by elements fixed by $\varphi$. From Lemma 2.2 (a) we get that $\diamond$ is the extension to $\eta_1$ of the unique connected subgroup $\diamond_{\dbQ_p}$ of $\grG$ whose Lie algebra is $\{e\in\Lie(\diamond)|\varphi(e)=e\}$. We refer to $\diamond_{\dbQ_p}$ as the {\it $\dbQ_p$-form of $\diamond$ with respect to $(\scrW,\varphi)$}. This example generalizes [Va7, part of proof of Claim 2.2.2].

\medskip
We recall with details the following well known general result. 
 
\bigskip\noindent
{\bf 2.3. Lemma.} {\it Let $\scrH$ be a reductive group scheme over $\dbZ_p$. Let $\scrT$ be a torus over $\dbZ_p$ whose generic fibre $\scrT_{\dbQ_p}$ is a maximal torus of $\scrH_{\dbQ_p}$. Then there exists an element $h\in \scrH(\dbQ_p)$ such that $\scrT$ is naturally a maximal torus of $h\scrH h^{-1}$.}

\medskip
\proof
We can assume that $\scrT_{W(k)}$ is split. The Galois group $\Gamma:=\Gal(k/\dbF_p)$ acts on the set of reductive group schemes over $W(k)$ whose generic fibres are $\scrH_{B(k)}$ and which have $\scrT_{W(k)}$ as a maximal torus (this set is non-empty, cf. [Ti2]). Let $o=\{\scrH_1,\ldots,\scrH_s\}$ be an orbit of this action and let $\scrH_o:=\scrH_1\times_{W(k)}\cdots\times_{W(k)}\scrH_s$. We have a canonical action of $\Gamma$ on $\scrH_0$ which is compatible with the action of $\Gamma$ on the $\dbZ_p$-scheme $\Spec W(k)$ and which for $\gamma\in\Gamma$ maps the factor $\scrH_i$ to the factor $\gamma(\scrH_i)$. As $\scrH_o$ is an affine group scheme, the resulting Galois descent datum on it with respect to $\Gamma$ is effective (cf. [BLR, Ch. 6, Sect. 6.1, Thm. 5]) and thus there exists a reductive group scheme $\scrH^\prime$ over $\dbZ_p$ such that we have an identification $\scrH_o=\scrH^\prime_{W(k)}$ compatible with the natural actions of $\Gamma$. The generic fibre of $\scrH^\prime$ is $\scrH_{\dbQ_p}$. The closed subgroup scheme $\prod_{i=1}^s \scrT_{W(k)}$ is naturally a maximal torus of $\scrH_o$ invariant under the action of $\Gamma$ and therefore $\scrT$ is a maximal torus of $\scrH^\prime$. Let $g\in\scrH^{\ad}(\dbQ_p)$ be such that $\scrH^\prime(\dbZ_p)=g\scrH(\dbZ_p) g^{-1}$, cf. [Ti2, Subsects. 2.5 and 3.2]. Thus we have $\scrH^\prime=g\scrH g^{-1}$, cf. [Va3, Lem. 2.3 (a)]. We can replace $g$ by $tg$, where $t\in\scrT_0(\dbQ_p)$ with $\scrT_0:=\text{Im}(\scrT_{\dbQ_p}\to \scrH^{\ad}_{\dbQ_p})$. 

Let $\scrM$ be a free $\dbZ_p$-module of finite rank such that $\scrH$ is a closed subgroup scheme of $\pmb{GL}_{\scrM}$. We recall the standard argument that allows us to consider the $\dbZ_p$-lattice $g(\scrM)$  of $\scrM[{1\over p}]$. Let $\ast$ be a Galois extension of $\dbQ_p$ such that there exists an element $h_0\in\scrH^{\der}(\ast)$ that maps to $g$. As the torus $\scrT_{0,B(k)}$ is split, by replacing $g$ with $tg$ for some $t\in\scrT_0(\dbQ_p)$, we can assume that $\ast=B(k)$. For $\gamma\in\Gal(B(k)/\dbQ_p)=\Gal(k/\dbF_p)$ we have $\gamma(h_0)=h_0z_{\gamma}$ with $z_{\gamma}\in Z(\scrH^{\der})(B(k))=Z(\scrH^{\der})(W(k))$. We easily get that $\gamma(h_0)(\scrM\otimes_{\dbZ_p} W(k))=h_0(\scrM\otimes_{\dbZ_p} W(k))$. From this and a standard Galois descent for modules we get that there exists a unique $\dbZ_p$-lattice $g(\scrM)$ of $\scrM[{1\over p}]$ with the property that $g(\scrM\otimes_{\dbZ_p} W(k))=h_0(\scrM\otimes_{\dbZ_p} W(k))$; it does not depend on the choice of $h_0\in \scrH^{\der}(B(k))$ mapping to $g$.

Let $(w_{\alpha})_{\alpha\in\scrJ_{\scrH}}$ be a family of tensors of $\scrT(\scrM[{1\over p}])$ such that $\scrH_{\dbQ_p}$ is the subgroup of $\pmb{\GL}_{\scrM[{1\over p}]}$ that fixes $w_{\alpha}$ for all $\alpha\in\scrJ_{\scrH}$ (cf. [De3, Prop. 3.1 c)]). The existence of $g(\scrM)$ and $h_0$ allow us to speak about the torsor of $\scrH$ which parametrizes isomorphisms between $(\scrM,(w_{\alpha})_{\alpha\in\scrJ_{\scrH}})$ and $(g(\scrM),(w_{\alpha})_{\alpha\in\scrJ_{\scrH}})$. As $\scrH_{\dbF_p}$ is connected and $\dbZ_p$ is henselian, this torsor is trivial and thus there exists $h\in\scrH(\dbQ_p)$ such that $h(\scrM)=g(\scrM)$. Therefore $\scrT$ is a maximal torus of $g\scrH g^{-1}=h\scrH h^{-1}$.\endproof

\bigskip\noindent
{\bf 2.4. Basic definitions.}
{\bf (a)} We say $\scrC$ has a {\it lift of quasi $\text{CM}$ type} if there exists a maximal torus $\scrT$ of $\scrG$ such that we have $\phi(\Lie(\scrT))=\Lie(\scrT)$.

\smallskip
{\bf (b)} We say $\scrC$ is {\it semisimple} (resp. {\it unramified semisimple}) if the $B(k)$-linear automorphism $\phi^r$ of $M[{1\over p}]$ is a semisimple element of $\scrG(B(k))$ (resp. is a semisimple element of $\scrG(B(k))$ such that an integral power of it has all its eigenvalues belonging to $B(\bar k)$). 

\smallskip
{\bf  (c)} By a {\it torus of $\scrG_{B(k)}$ of $\dbQ_p$-endomorphisms of $\scrC$} we mean a torus $\scrT_{1,B(k)}$ of $\scrG_{B(k)}$ whose Lie algebra is $B(k)$-generated by elements fixed by $\phi$. Let $\scrT_{1,\dbQ_p}$ be the $\dbQ_p$-form of $\scrT_{1,B(k)}$  with respect to $(M[{1\over p}],\phi)$, cf. Example 2.2.1. Let $K$ be the smallest Galois extension of $\dbQ_p$ over which $\scrT_{1,\dbQ_p}$ splits. Let $K_1$ be the smallest unramified extension of $K$ which is unramified over a totally ramified extension $K_{1r}$ of $\dbQ_p$. Let $K_2$ be the composite field of $K_1$ and $B(k)$. Let $K_{2u}$ be the maximal unramified extension of $\dbQ_p$ included in $K_2$. It is easy to see that $K_2$ is the composite field of $K_{1r}$ and $K_{2u}$.

\smallskip
{\bf (d)} By an {\it $E$-pair} of $\scrC$ we mean a pair $(\scrT_{1,B(k)},\mu_1)$, where $\scrT_{1,B(k)}$ is a maximal torus of $\scrG_{B(k)}$ of $\dbQ_p$-endomorphisms of $\scrC$ and $\mu_1:\dbG_m\to \scrT_{1,K_1}$ is a cocharacter such that $\mu_{1,K_2}$, when viewed as a cocharacter of $\scrG_{K_2}$, is $\scrG(K_2)$-conjugate to $\mu_{K_2}$. If $\mu_1$ is definable over an unramified extension of $\dbQ_p$, then we refer to $(\scrT_{1,B(k)},\mu_1)$ as an {\it unramified $E$-pair}. By an {\it $E$-triple} of $\scrC$ we mean a triple $(\scrT_{1,B(k)},\mu_1,\tau)$, where $(\scrT_{1,B(k)},\mu_1)$ is an $E$-pair and where $\tau=(\tau_1,\ldots,\tau_l)$ is an $l$-tuple of elements of $\Gal(K_2/\dbQ_p)$ whose restrictions to $K_{2u}$ are all equal to the Frobenius automorphism $\grF_{2u}$ of $K_{2u}$ whose fixed field is $\dbQ_p$. Here $l\in\dbN^*$ while ``$E$" stands for endomorphisms. For $s\in\dbN^*$ and $j\in \{1,\ldots,l\}$ let $\tau_{sl+j}:=\tau_j$.

\smallskip
{\bf (e)} We say that an $E$-pair $(\scrT_{1,B(k)},\mu_1)$ of $\scrC$ satisfies the {\it $\grC$ condition} if there exists an $E$-triple $(\scrT_{1,B(k)},\mu_1,\tau)$ of $\scrC$ such that the following condition holds:

\medskip
{\bf (e1)} the product of the cocharacters of $\scrT_{1,K_2}$ of the form $\tau_{dl}\tau_{dl-1}\cdots\tau_j(\mu_{1,K_2})$ with $j\in\{1,\ldots,dl\}$, factors through the subtorus $Z^0(\scrG_{K_2})$ of $\scrT_{1,K_2}$; here $d\in\dbN^*$ is the smallest number such that $\mu_{1,K_2}$ is fixed by each element of $\Gal(K_2/\dbQ_p)$ that can be obtained from the product $\tau_{dl}\tau_{dl-1}\cdots\tau_1$ via a circular rearrangement of it.

\medskip
If moreover $l=1$ we say $(\scrT_{1,B(k)},\mu_1)$ satisfies the {\it cyclic $\grC$ condition}. 

\smallskip
{\bf (f)} We assume that $\scrC$ is basic (as defined in Subsection 2.1). We say $\grR$ (resp. $\grU$) holds for $\scrC$ if there exists an $E$-pair (resp. an unramified $E$-pair) of $\scrC$ that satisfies the $\grC$ condition. We say $T\grR$ (resp. $T\grU$) holds for $\scrC$ if each maximal torus $\scrT_{1,B(k)}$ of $\scrG_{B(k)}$ of $\dbQ_p$-endomorphisms of $\scrC$ (resp. each maximal torus $\scrT_{1,B(k)}$ of $\scrG_{B(k)}$ of $\dbQ_p$-endomorphisms of $\scrC$ which splits over $B(\bar k)$) is part of an $E$-pair (resp. of an unramified $E$-pair) of $\scrC$ that satisfies the $\grC$ condition. We say $Q\grR$ (resp. $Q\grU$) holds for $\scrC$ if there exists a $k_1$ and an $E$-pair (resp. an unramified $E$-pair) of $\scrC\otimes k_1$ that satisfies the $\grC$ condition. We say $TT\grR$ (resp. $TT\grU$) holds for $\scrC$ if for each $k_1$, $T\grR$ (resp. $T\grU$) holds for $\scrC\otimes k_1$.

\smallskip
{\bf (g)} We do not assume that $\scrC$ is basic. We say $\grR$ (resp. $\grU$, $T\grR$, $T\grU$, $Q\grR$, $Q\grU$, $TT\grR$, $TT\grU$) holds for $\scrC$ if there exists $h\in \grI(\scrC)$ such that the triple $(h(M),\phi,L_{\scrG(h)}^0(\phi))$ is a basic Shimura $F$-crystal over $k$ and $\grR$ (resp. $\grU$, $T\grR$, $T\grU$, $Q\grR$, $Q\grU$, $TT\grR$, $TT\grU$) holds for it. 

{\bf (h)} We say that an $E$-pair $(\scrT_{1,B(k)},\mu_1)$ of $\scrC$ is {\it admissible} if the filtered module 
$$(M[{1\over p}],\phi,F^1_{K_2})$$ 
over $K_2$ is admissible. Here $F^1_{K_2}$ is the maximal direct summand of $M\otimes_{W(k)} K_2$ on which $\mu_{1,K_2}$ acts via the inverse of the identical character of $\dbG_m$ (i.e., via weight $-1$).

{\bf (i)}  We say $\grA$ holds for $\scrC$ if there exists an $E$-pair $(\scrT_{1,B(k)},\mu_1)$ of $\scrC$ which is admissible. We say $T\grA$ holds for $\scrC$ if each maximal torus $\scrT_{1,B(k)}$ of $\scrG_{B(k)}$ of $\dbQ_p$-endomorphisms of $\scrC$ can be extended to an $E$-pair $(\scrT_{1,B(k)},\mu_1)$ of $\scrC$ which is admissible. As in (f), we speak also about $Q\grA$ or $TT\grA$ holding for $\scrC$.

{\bf (j)} We say $\scrC$ is {\it Uni-ordinary} if it has a lift $F^1$ such that $L^0_{\scrG}(\phi)_{B(k)}$ normalizes $F^1[{1\over p}]$ (i.e., we have $L^0_{\scrG}(\phi)_{B(k)}\leqslant\scrP_{B(k)}$). We say $\scrC$ is {\it Isouni-ordinary} if there exist elements $g\in \scrG(W(k))$ and $h\in \scrG(B(k))$ such that $\scrC_g$ is Uni-ordinary and we have $h\phi=g\phi h$ (i.e., $\scrC$ is inner isomorphic to a Uni-ordinary $\scrC_g$). 

{\bf (k)} A {\it principal bilinear quasi-polarization} of $\scrC$ is a perfect bilinear form $\lambda_M:M\times M\to W(k)$ whose $W(k)$-span is normalized by $\scrG$ and for which we have $\lambda_M(\phi(x),\phi(y))=p\sigma(\lambda_M(x,y))$ for all elements $x,y\in M.$ 

\medskip\noindent
{\bf 2.4.1. Example.} Let $(\scrT_{1,B(k)},\mu_1)$ be an $E$-pair of $\scrC$ such that the product of the cocharacters of $\scrT_{1,K_1}$ which belong to the $\Gal(K_1/\dbQ_p)$-orbit of $\mu_1$, factors through the subtorus $Z^0(\scrG_{K_1})$ of $\scrT_{1,K_1}$. We choose an element $\tau_0\in\Gal(K_2/\dbQ_p)$ whose restriction to $K_{2u}$ is $\grF_{2u}$ and whose order $o$ is the same as of $\grF_{2u}$ (its existence is implied by the fact that $K_2$ is the composite field of $K_{1r}$ and $K_{2u}$). Let $\{e_1,\ldots,e_s\}$ be the elements of $\Gal(K_2/K_{2u})$ listed in such a way that $e_s=1_K$. We have $\Gal(K_2/\dbQ_p)=\{e_a\tau_0^b|1\Le a\Le s,\,1\Le b\Le o\}$. Let $l:=os$ and $e_0:=e_s$. We define $\tau=(\tau_1,\ldots,\tau_l)$ as follows. For $i\in\{1,\ldots,l\}$ let $\tau_i:=\tau_0$ if $o$ does not divide $i$ and let $\tau_i:=e_{j-1}^{-1}e_j\tau_0$ if $i=o(s+1-j)$ with $j\in\{1,\ldots,s\}$. As $\tau_l\tau_{l-1}\cdots\tau_1=1_K$, let $d:=1$. As we have $\Gal(K_2/\dbQ_p)=\{\tau_l\tau_{l-1}\cdots\tau_j|1\Le j\Le l\}$, the condition 2.4 (e1) holds (cf. our hypothesis on the $E$-pair $(\scrT_{1,B(k)},\mu_1)$ of $\scrC$). Thus the $E$-pair $(\scrT_{1,B(k)},\mu_1)$ of $\scrC$ satisfies the $\grC$ condition. 

If moreover the $E$-pair $(\scrT_{1,B(k)},\mu_1)$ is unramified, then there exists a finite field extension $k_1=\dbF_{p^o}$ of $k$ such that $K_2=K_{2u}=B(k_1)$ and we can take $l=1$, $\tau=(\tau_1)$, and $d$ to be the smallest divisor of $o$ such that $\tau_0^d$ fixes $\mu_{1,K_2}$; thus $\scrC$ satisfies the cyclic $\grC$ condition.

\medskip\noindent
{\bf 2.4.2. Example.} Let $m\in\dbN^*$. We assume that the rank of $M$ is $2m$, that $\scrG$ is a product $\scrG_1\times_{W(k)} \cdots \times_{W(k)}  \scrG_m$ of $m$ copies of $\pmb{\GL}_2$, that $\phi$ permutes transitively the $\Lie(\scrG_i)[{1\over p}]$'s with $i\in\{1,\ldots,m\}$, that for each $i\in\{1,\ldots,m\}$ the image of $\mu$ in $\scrG_i$ does not factor through $Z(\scrG_i)$, and that the representation of $\scrG$ on $M$ is the direct sum of the standard rank 2 representations of the $m$ copies. This last assumption implies that $\scrG$ is the subgroup scheme of $\pmb{\GL}_M$ that fixes a semisimple $\dbZ_p$-subalgebra of $\End(M)$ formed by elements fixed by $\phi$ (this is also a particular case of property 2.5 (a) below). The rank of $F^1$ is $m$ and $\scrP=\scrP_1\times_{W(k)} \cdots\times_{W(k)} \scrP_m$ is a Borel subgroup scheme of $\scrG$. We also assume that there exists a maximal torus $\scrT=\scrT_1\times_{W(k)}  \cdots \times_{W(k)}  \scrT_m$ of $\scrP$ such that we have $\phi(\Lie(\scrT))=\Lie(\scrT)$ and $\phi(\Lie(\scrP))\subseteq\Lie(\scrP)$; thus the Dieudonn\'e module $(M,\phi)$ is ordinary. 

Let $g\in \scrG(W(k))$ be such that $\scrC_g$ is not basic. Thus $P_{\scrG}^+(g\phi)=\prod_{i=1}^m P_{\scrG}^+(g\phi)\cap \scrG_i$ is a Borel subgroup scheme of $\scrG$ and therefore $L_{\scrG}^0(\phi)_{B(k)}$ is a split maximal torus of $\scrG_{B(k)}$. We check that $\scrC_g$ is Isouni-ordinary. Based on [Va7, Thm. 3.1.2 (b) and (c)], up to a replacement of $g\phi$ by $hg\phi h^{-1}$ with $h\in \scrG(B(k))$, we can assume that $L_{\scrG}^0(\phi)$ is a maximal torus of $\scrG$ through which $\mu$ factors. Thus $L_{\scrG}^0(\phi)$ commutes with $\mu$ and therefore it is a maximal torus of $\scrP$.  Therefore $\scrC_g$ is Uni-ordinary. 

We now take $m=3$. Let $w:=(w_1,w_2,1_M)\in \scrG_1(W(k))\times \scrG_2(W(k))\times \scrG_3(W(k))$ be an element that normalizes $\scrT$ and such that for $i\in\{1,2\}$ the element $w_i$ takes $\scrP_i$ to its opposite $\scrP_i^{\text{opp}}$ with respect to $\scrT_i$. The Newton polygon slopes of $(M,w\phi)$ are ${1\over 3}$ and ${2\over 3}$ with multiplicities $3$. We have $L_{\scrG}^0(w\phi)=\scrT\leqslant\scrP$. Thus $\scrC_w$ is Uni-ordinary. Let $\scrU_1$ be the unipotent radical of $\scrP_1^{\text{opp}}$; it is a subgroup scheme of $U_{\scrG}^+(w\phi)$. Let $g_1\in \scrU_1(W(k))$ be such that modulo $p$ it is not the identity element. As $g_1\in U_{\scrG}^+(w\phi)(W(k))$, it is easy to see that we have $P_{\scrG}^+(g_1w\phi)=P_{\scrG}^+(w\phi)$. Thus $\scrC_{g_1w}$ is not basic and therefore (cf. previous paragraph) it is Isouni-ordinary. We show that the assumption that $\scrC_{g_1w}$ is Uni-ordinary leads to a contradiction. It is easy to see that this assumption implies that $L_{\scrG}^0(g_1w\phi)$ is a maximal torus of the Borel group scheme $P_{\scrG}^+(g_1w\phi)=P_{\scrG}^+(w\phi)$ (for instance, see Fact 2.5.7 below for a general result) which normalizes $F^1/pF^1$ (cf. very definitions). Let $b\in P_{\scrG}^+(w\phi)(W(k))$ be such that it normalizes $F^1/pF^1$ and we have $b(L_{\scrG}^0(g_1w\phi))b^{-1}=\scrT$, cf. [Bo, Ch. V, Thm. 19.2] and [DG, Vol. II, Exp. IX, Thms. 3.6 and 7.1]. Thus $bg_1w\phi b^{-1}=g_2w\phi$, where $g_2\in P_{\scrG}^+(w\phi)(W(k))$ normalizes $\scrT$. Therefore $g_2\in \scrT(W(k))$. It is easy to see that the images of $g_1$ and $b^{-1}g_2=g_1w\phi(b^{-1})w^{-1}$ in $\scrT_1(k)\backslash \scrG_1(k)/\scrT_1(k)$ are equal. As $g_1$ modulo $p$ is a non-identity element of $\scrU_1(k)$ and $b^{-1}g_2$ modulo $p$ belongs to $\scrP(k)$, we reached a contradiction. We conclude that $\scrC_{g_1w}$ is Isouni-ordinary without being Uni-ordinary.

\bigskip\noindent
{\bf 2.5. Some $\dbZ_p$ structures.} We begin by presenting two ways to construct $\dbZ_p$ structures.

\smallskip
{\bf (a)} Let $O$ be a $W(k)$-submodule of $\scrT(M)$ of finite rank such that $\phi(O)=O$ and $\phi^r$ fixes $O$. As $\phi(O)=O$, $O\otimes_{W(k)} W(\bar k)$ has a $\dbZ_p$-basis formed by elements fixed by $\phi\otimes\sigma_{\bar k}$. As $\phi^r$ fixes $O$, the $W(k)$-span of this $\dbZ_p$-basis contains $O$ and by reasons of ranks we get that it is $O$. Therefore $O$ is $W(k)$-generated by elements fixed by $\phi$. Similarly, a finite dimensional $B(k)$-vector subspace of $\scrT(M)[{1\over p}]$ fixed by $\phi^r$ and normalized by $\phi$, is $B(k)$-generated by elements fixed by $\phi$.

\smallskip
{\bf (b)} Let $\sigma_{\phi}:=\phi\circ\mu(p)$; it is a $\sigma$-linear automorphism of $M$. Let
$$M_{\dbZ_p}:=\{m\in M\otimes_{W(k)} W(\bar k)|(\sigma_{\phi}\otimes\sigma_{\bar k})(m)=m\}\;\;\;\text{and}\;\;\;\text{let}\;\;\;M_{\dbQ_p}:=M_{\dbZ_p}[{1\over p}].$$
We have $M\otimes_{W(k)} W(\bar k)=M_{\dbZ_p}\otimes_{\dbZ_p} W(\bar k)$. 

Let $\square$ be a closed subgroup scheme of $\pmb{\GL}_M$ which is an integral scheme. We assume that both $\mu$ and $\phi$ normalize $\Lie(\square_{B(k)})$. This implies that $\sigma_{\phi}$ normalizes $\Lie(\square)$. Thus $\square_{B(\bar k)}$ is the extension to $B(\bar k)$ of a connected subgroup of $\pmb{\GL}_{M_{\dbZ_p}[{1\over p}]}$, cf. Example 2.2.1. If $\square_{\dbZ_p}$ is the schematic closure of $\square_{\dbQ_p}$ in $\pmb{\GL}_{M_{\dbZ_p}}$, then its extension to $W(\bar k)$ is $\square_{W(\bar k)}$. 

As $\mu$ and $\phi$ normalize $\Lie(\scrG_{B(k)})$, $\sigma_{\phi}$ normalizes $\Lie(\scrG)$. Thus from the previous paragraph we get the existence of a unique closed subgroup scheme $\scrG_{\dbZ_p}$ of $\pmb{\GL}_{M_{\dbZ_p}}$ whose extension to $W(\bar k)$ is $\scrG_{W(\bar k)}$; it is a reductive group scheme over $\dbZ_p$. If $\square$ is a subgroup scheme of $\scrG$, then $\square_{\dbZ_p}$ is a subgroup scheme of $\scrG_{\dbZ_p}$.

\medskip\noindent
{\bf 2.5.1. Applications of the standard axioms.} Until the end of the paper we will assume that the two standard axioms of Definition 1.2.4 (b) hold for $\scrC$. Axiom 1.2.4 (b.i) implies that $\phi^r\in \scrG(B(k))$ and that we have $t_{\alpha}\in\scrT(M_{\dbZ_p})$ for all $\alpha\in\scrJ$. Thus the pair $(M_{\dbZ_p},(t_{\alpha})_{\alpha\in\scrJ})$ is a $\dbZ_p$ structure of $(M\otimes_{W(k)} W(\bar k),(t_{\alpha})_{\alpha\in\scrJ})$. The difference between each two such $\dbZ_p$ structures of $(M\otimes_{W(k)} W(\bar k),(t_{\alpha})_{\alpha\in\scrJ})$ is measured by a class in $H^1(\Gal(\dbZ_p^{\text{un}}/\dbZ_p),\scrG_{\dbZ_p})$, where $\dbZ_p^{\text{un}}$ is the maximal unramified, profinite discrete valuation ring extension of $\dbZ_p$. From Lang theorem (see [Se2, p. 132] and [Bo, Ch. V, Subsects. 16.3 to 16.6]) we get that this class is trivial. Thus the isomorphism class of the triple $(M_{\dbZ_p},\scrG_{\dbZ_p},(t_{\alpha})_{\alpha\in\scrJ})$ does not depend on the choice of the Hodge cocharacter $\mu:\dbG_m\to\scrG$ of $\scrC$. Also by replacing $\phi$ with $g\phi$, where $g\in \scrG(W(k))$, the isomorphism class of $(M_{\dbZ_p},\scrG_{\dbZ_p},(t_{\alpha})_{\alpha\in\scrJ})$ remains the same. From Lang theorem we also get that each torsor of $\scrG$ is trivial. This implies that there exists an isomorphism
$$i_M:M_{\dbZ_p}\otimes_{\dbZ_p} W(k)\arrowsim M\leqno (2)$$
that takes $t_{\alpha}$ to $t_{\alpha}$ for all $\alpha\in\scrJ$. Thus $\scrG=\scrG_{W(k)}$ (i.e., our notations match) and we refer to the triple $(M_{\dbZ_p},\scrG_{\dbZ_p},(t_{\alpha})_{\alpha\in\scrJ})$ as the {\it $\dbZ_p$ structure} of $(M,\phi,\scrG,(t_{\alpha})_{\alpha\in\scrJ})$. 

Axiom 1.2.4 (b.ii) is inserted for practical reasons i.e., to exclude situations that are not related to Shimura pairs of Hodge type and to get the following result. 

\medskip\noindent
{\bf 2.5.2. Theorem.} {\it {\bf (a)} The closed subgroup scheme $C=C_{\pmb{\GL}_M}(\scrG)$ of $\pmb{\GL}_M$ is reductive (and thus $C_1:=C_{\pmb{\GL}_M}(C)$ is well defined). 

\smallskip
{\bf (b)} The Lie algebra $\Lie(C)$ is $W(k)$-generated by elements fixed by $\phi$.}

\medskip
\proof
To prove (a) we work only with the $\scrG$-module $M$ and apply the axiom 1.2.4 (b.ii); thus the below reduction steps do not pay attention to $\phi$. To prove (a) we can assume that $\scrG$ is split. By considering a direct sum decomposition of $M$ into $\scrG$-modules on which $Z^0(\scrG)$ acts via distinct characters, we can assume as well that $Z^0(\scrG)=Z(\pmb{\GL}_M)$. Let $M[{1\over p}]:=\oplus_{i=1}^n M_i[{1\over p}]$ be a direct sum decomposition into irreducible $\scrG_{B(k)}$-modules, cf. Weyl complete reducibility theorem. Let $M_i:=M\cap M_i[{1\over p}]$. Thus $\oplus_{i=1}^n M_i$ is a $\scrG$-submodule of $M$. Due to the axiom 1.2.4 (b.ii), each simple factor of $\Lie(\scrG^{\der}_{B(k)})$ is of classical Lie type and the representation of $\Lie(\scrG^{\der}_{B(k)})$ on each $M_i[{1\over p}]$ is a tensor product of irreducible representations which are either trivial or are associated to minuscule weights (see [Se1, Prop. 7 and Cor. 1 of p. 182]). Thus the $\scrG_k$-module $M_i/pM_i$ is absolutely irreducible and its isomorphism class depends only on the isomorphism class of the $\scrG_{B(k)}$-module $M_i[{1\over p}]$, cf. the below well known Fact  2.5.3 and [Ja, Part I, Sect. 10.9]. 

By induction on $n\in\dbN^*$ we show that we can choose the decomposition $M[{1\over p}]:=\oplus_{i=1}^n M_i[{1\over p}]$ such that we have $M=\oplus_{i=1}^n M_i$. The case $n=1$ is trivial. The passage from $n$ to $n+1$ goes as follows. We have a short exact sequence $0\to M_1\to M\to M/M_1\to 0$ of $\scrG$-modules. Using induction, it suffices to consider the case $n=2$; thus the $W(k)$-monomorphism $M_2\hookrightarrow M/M_1$ becomes an isomorphism after inverting $p$. If the $\scrG_{B(k)}$-modules $M_1[{1\over p}]$ and $M_2[{1\over p}]$ are not isomorphic, then the $\scrG_k$-modules $M_1/pM_1$ and $M_2/pM_2$ are not isomorphic and therefore the natural $k$-linear map $M_1/pM_1\oplus M_2/pM_2\to M/pM$ is injective; this implies that we have $M=M_1\oplus M_2$. We assume now that the $\scrG_{B(k)}$-modules $M_1[{1\over p}]$ and $M_2[{1\over p}]$ are isomorphic. Thus $M_1$ and $M_2$ are isomorphic $\scrG$-modules. If they are trivial $\scrG$-modules, then we can replace $M_2$ by any direct supplement of $M_1$ in $M$ and thus we have $M=M_1\oplus M_2$. We now consider the case when $M_1$ and $M_2$ are non-trivial $\scrG$-modules. Let $\tilde M$ be a $W(k)$-submodule of $M$ which is a $\scrG$-module isomorphic to $M_1\oplus M_2$ and for which the length of the torsion $W(k)$-module $M/\tilde M$ has the smallest possible value $l\in\dbN$. We show that the assumption $l\neq 0$ leads to a contradiction. Due the smallest property of $l$, $\tilde M$ is not contained in $pM$. From this and the fact that $l\neq 0$, we get that the $\scrG$-module $\tilde M^\prime:={1\over p}\tilde M\cap M$ satisfies $\tilde M\subsetneqq \tilde M^\prime\subsetneqq {1\over p}$. Therefore the quotient $\tilde M^\prime/\tilde M$ is an irreducible $\scrG_k$-submodule of ${1\over p}\tilde M/\tilde M$. As the $\scrG$-module ${1\over p}\tilde M$ is isomorphic to $M_1\oplus M_2$, there exists an irreducible $\scrG$-submodule of $\tilde M^\prime$ whose reduction modulo $p$ is $\tilde M^\prime/\tilde M$ and which is isomorphic to $M_1$. From this we easily get that the $\scrG$-module $\tilde M^\prime$ is isomorphic to $M_1\oplus M_2$. As the length of the torsion $W(k)$-module $M/\tilde M^\prime$ is less than $l$, we reached the desired contradiction. 

Thus to prove (a) we can assume that $M=\oplus_{i=1}^n M_i$. As the isomorphism class of $M_i$ is uniquely determined by the isomorphism class of the $\scrG_{B(k)}$-module $M_i[{1\over p}]$, we can write $M=\oplus_{j\in J} M_j^{n_j}$, where each $M_j$ is isomorphic to some $M_i$, where $n_j\in\dbN^*$, and where for two distinct elements $j_1, j_2\in J$ the $\scrG$-modules $M_{j_1}$ and $M_{j_2}$ are not isomorphic and the $\scrG_k$-modules $M_{j_1}/pM_{j_1}$ and $M_{j_2}/pM_{j_2}$ are not isomorphic. We easily get that the group scheme $C$ is isomorphic to a product $\prod_{j\in J} \pmb{\GL}_{n_j}$ and therefore it is a reductive, closed subgroup scheme of $\pmb{\GL}_M$. Thus (a) holds. As $\phi(\Lie(C))=\Lie(C)$ and $\phi^r\in \scrG(B(k))$ fixes $\Lie(C)$, part (b) follows from the property 2.5 (a).\endproof

\medskip\noindent
{\bf 2.5.3. Fact.} {\it Let $\scrH$ be a split, simply connected group scheme over $\dbZ$ whose adjoint is absolutely simple and of classical Lie type $\theta$. Let $\scrT$ be a maximal split torus of $\scrH$. Let $\rho_{\varpi}:\scrH\to \pmb{\GL}_{\scrZ}$ be the representation associated to a minuscule weight $\varpi$ of the root system of the inner conjugation action of $\scrT$ on $\Lie(\scrH)$ (thus $\scrZ$ is a free $\dbZ$-module of finite rank, cf. [Hu, Subsect. 27.1]). Then the special fibres of $\rho_{\varpi}$ are absolutely irreducible.} 

\medskip
\proof
We use the notations of [Bou1, planches I to IV]. The minuscule weights are: $\varpi_i$ with $i\in\{1,\ldots,n\}$ if $\theta=A_n$, $\varpi_n$ if $\theta=B_n$, $\varpi_1$ if $\theta=C_n$, $\varpi_1$, $\varpi_{n-1}$, and $\varpi_n$ if $\theta=D_n$ (see [Bou2, Ch. VIII, Subsect. 7.3, Rm.] and [Se1, pp. 185--186]). Let $\grW$ be the Weyl group of $\scrH$ with respect to $\scrT$. Let $\grW_{\varpi}$ be the subgroup of $\grW$ that fixes $\varpi$. We have $\dim_{\dbZ}(\scrZ)=[\grW:\grW_{\varpi}]$, cf. [Bou2, Ch. VIII, Subsect. 7.3, Prop. 6]. Thus for each prime $p$, the absolutely irreducible representation of $\scrH_{\dbF_p}$ associated to weight $\varpi$ has dimension at least $\dim_{\dbZ}(\scrZ)$. As it is isomorphic to the representation of $\scrH_{\dbF_p}$ on a factor of the composition series of the fibre of $\rho_{\varpi}$ over $\dbF_p$, by reasons of dimensions we get that this fibre is absolutely irreducible.\endproof

\medskip\noindent
{\bf 2.5.4. Classes associated to admissible $E$-pairs.} We consider an $E$-pair $(\scrT_{1,B(k)},\mu_1)$ of $\scrC$ which is admissible; let $F^1_{K_2}$ be as in Subsection 2.4 (h). Let $\rho:\Gal(K_2)\to \pmb{\GL}_\scrW(\dbQ_p)$ be the admissible Galois representation that corresponds to $(M,\phi,F^1_{K_2})$. Thus $\scrW$ is a finite dimensional $\dbQ_p$-vector space and we have a $\Gal(K_2)$-isomorphism
$$\scrW\otimes_{\dbQ_p} B_{\text{dR}}(K_2)\arrowsim M\otimes_{W(k)} K_2\otimes_{K_2} B_{\text{dR}}(K_2)\leqno (3)$$
which respects the tensor product filtrations (the filtration of $\scrW$ is trivial and the filtration of $M\otimes_{W(k)} K_2$ is defined by $F^1_{K_2}$). For $\alpha\in\scrJ$, let $v_{\alpha}\in\scrT(\scrW)$ be the tensor that corresponds to $t_{\alpha}$ via (3) and Fontaine comparison theory. Let
$\gamma\in H^1(\dbQ_p,\scrG_{\dbQ_p})$ be the class of the right torsor of $\scrG_{\dbQ_p}$ that parameterizes isomorphisms between $(M_{\dbQ_p},(t_{\alpha})_{\alpha\in\scrJ})$ and $(\scrW,(v_{\alpha})_{\alpha\in\scrJ})$ (such a torsor exists, cf. (2) and (3)). 

\medskip\noindent
{\bf 2.5.5. Lemma.} {\it We assume that $\scrG^{\der}$ is simply connected. Then the class $\gamma$ is trivial i.e., there exists a $\dbQ_p$-linear isomorphism $\grJ:M_{\dbQ_p}\arrowsim \scrW$
that takes $t_{\alpha}$ to $v_{\alpha}$ for all $\alpha\in\scrJ$.}

\medskip
\proof  
The pointed set $H^1(\dbQ_p,\scrG^{\der}_{\dbQ_p})$ has also only one class (cf. [Kn, Thm. 1]) and we have an exact complex $H^1(\dbQ_p,\scrG^{\der}_{\dbQ_p})\to H^1(\dbQ_p,\scrG_{\dbQ_p})\to H^1(\dbQ_p,\scrG^{\ab}_{\dbQ_p})$ of pointed sets. Thus it suffices to show that the image $\gamma^{\ab}$ of $\gamma$ in  $H^1(\dbQ_p,\scrG^{\ab}_{\dbQ_p})$ is the trivial class.

Let $\scrG^{\ab}_{\dbZ_p}\hookrightarrow \pmb{GL}_{M^{\ab}_{\dbZ_p}}$ be a closed embedding monomorphism; we can assume that $M^{\ab}_{\dbZ_p}$ is a $\dbZ_p$-submodule of $\scrT(M_{\dbZ_p})$. Let $\scrG^{\ab,\text{big}}_{\dbZ_p}$ be a maximal torus of $\pmb{GL}_{M^{\ab}_{\dbZ_p}}$ that contains $\scrG^{\ab}_{\dbZ_p}$. The cocharacter $\mu_{1}^{\ab}:\dbG_m\to \scrG^{\ab}_{K_1}$ induced by $\mu_1$ is defined over a subfield of $K_1$ unramified over $\dbQ_p$ and thus by performing the operation $\grO_1$ we can assume it is the extension to $K_1$ of a cocharacter $\mu_{1,B(k)}^{\ab}:\dbG_m\to \scrG^{\ab}_{B(k)}$ over $B(k)$. Thus $(M^{\ab}_{\dbZ_p}\otimes_{\dbZ_p} B(k),(F^{i})_{i\in\dbZ},(1_{M^{\ab}_{\dbZ_p}}\otimes\sigma)\circ \mu_{1}^{\ab}({1\over p}))$ is a filtered $F$-isocrystal over $k$; here each $F^{i}$ is the maximal $B(k)$-vector subspace of $M^{\ab}_{\dbZ_p}\otimes_{\dbZ_p} B(k)$ on which $\dbG_m$ acts via  $\mu_{1}^{\ab}$ through weights at the most equal to $-i$.

Let $\rho:\Gal(B(k))\to \pmb{\GL}_{\scrW^{\ab}}(\dbQ_p)$ be the admissible Galois representation that corresponds to $(M^{\ab}_{\dbZ_p}\otimes_{\dbZ_p} B(k),(F^{i})_{i\in\dbZ},(1_{M^{\ab}_{\dbZ_p}}\otimes\sigma)\circ \mu_{1}^{\ab}({1\over p}))$. Here $\scrW^{\ab}$ is a $\dbQ_p$-vector space of the same dimension as the rank of $M^{\ab}_{\dbZ_p}$. Let $\scrT^{\ab,\text{big}}_{\dbQ_p}$ (resp. $\scrT^{\ab}_{\dbQ_p}$) be the maximal torus (resp. be the torus) of $\pmb{\GL}_{\scrW^{\ab}}$ that corresponds naturally to  $\scrG^{\ab,\text{big}}_{\dbQ_p}$ (resp. to $\scrG^{\ab}_{\dbQ_p}$) via Fontaine comparison theory. The torus $\scrT^{\ab,\text{big}}_{\dbQ_p}$ is naturally identified with $\scrG^{\ab,\text{big}}_{\dbQ_p}$ (being the group of invertible elements of an \'etale $\dbQ_p$-algebra) and it contains $\scrT^{\ab}_{\dbQ_p}$ as a subtorus. This implies that $\scrT^{\ab}_{\dbQ_p}$ extends to a torus $\scrT^{\ab}_{\dbZ_p}$ over $\dbZ_p$ which is the schematic closure of $\scrT^{\ab}_{\dbQ_p}$ in $\scrG^{\ab,\text{big}}_{\dbZ_p}$. Let $\scrW_{\dbZ_p}^{\ab}$ be a $\dbZ_p$-lattice of $\scrW^{\ab}$ normalized by $\scrT^{\ab}_{\dbZ_p}$. Let $(v_{\alpha})_{\alpha\in\scrJ^{\ab}}$ be a family of tensors of $\scrT(\scrW_{\dbZ_p}^{\ab})$ such that the torus $\scrT^{\ab}_{\dbZ_p}$ is the subgroup scheme of $\pmb{GL}_{\scrW_{\dbZ_p}^{\ab}}$ that fixes $v_{\alpha}$ for all $\alpha\in\scrJ^{\ab}$. For $\alpha\in\scrJ^{\ab}$ let $t_{\alpha}$ be the tensor of $\scrT(M_{\dbZ_p}^{\ab}\otimes_{\dbZ_p} B(k))$ that corresponds to $v_{\alpha}$ via Fontaine comparison theory; we can assume that $\scrJ^{\ab}\subseteq \scrJ$ and that each $t_{\alpha}$ with $\alpha\in\scrJ^{\ab}$ is identified canonically with $t_{\alpha}$ of the axiom 1.2.4 (b.i) under the identification of $M^{\ab}_{\dbZ_p}$ with a $\dbZ_p$-submodule of $\scrT(M_{\dbZ_p})$. 

From [Ki2, Prop. (1.3.4) and Cor. (1.3.5)] we get that there exist a $W(k)$-lattice $M^{\ab,\prime}_{W(k)}$ of $M_{\dbZ_p}^{\ab}\otimes_{\dbZ_p} B(k)$ and an isomorphism $\scrW_{\dbZ_p}^{\ab}\otimes_{\dbZ_p} W(k)\arrowsim M_{W(k)}^{\ab,\prime}$ that maps $v_{\alpha}$ to $t_{\alpha}$ for all $\alpha\in\scrJ^{\ab}$. Moreover, $\mu_{1,B(k)}^{\ab}$ extends to a cocharacter of $\scrG^{\ab}_{W(k)}$ and thus of $\pmb{GL}_{M^{\ab,\prime}_{W(k)}}$ and $1_{M_{\dbZ_p}^{\ab}}[{1\over p}]\otimes \sigma$ normalizes $M^{\ab,\prime}_{W(k)}$. This allows us to define a $\dbZ_p$ structure $M^{\ab,\prime}_{\dbZ_p}:=\{x\in M^{\ab,\prime}_{W(k)}|(1_{M_{\dbZ_p}^{\ab}}[{1\over p}]\otimes \sigma)(x)=x\}$ of $M^{\ab,\prime}_{W(k)}$ which is a $\dbZ_p$-lattice of $M_{\dbZ_p}^{\ab}[{1\over p}]$. We have $t_{\alpha}\in\scrT(M^{\ab,\prime}_{\dbZ_p})$ for all $\alpha\in\scrJ^{\ab}$. As each right torsor of $\scrT^{\ab}_{\dbZ_p}$ is trivial, we can choose the last isomorphism to map $\scrW_{\dbZ_p}^{\ab}$ onto $M_{\dbZ_p}^{\ab,\prime}$. Therefore there exists an isomorphism $\scrW^{\ab}\arrowsim M^{\ab}_{\dbZ_p}[{1\over p}]$ which maps $v_{\alpha}$ to $t_{\alpha}$ for all $\alpha\in\scrJ^{\ab}$. This implies that the class $\gamma^{\ab}$ is trivial. We conclude that the class $\gamma$ itself is trivial.\endproof

\medskip\noindent
{\bf 2.5.6. An ``adjoint'' reduction step.} We assume that there exists a non-trivial product decomposition $\scrG_{\dbZ_p}^{\ad}=\scrV_1\times_{\dbZ_p} \scrV_2$. Let $\phi_0=i_M^{-1}\phi i_M:M_{\dbZ_p}\otimes_{\dbZ_p} W(k)\to M_{\dbZ_p}\otimes_{\dbZ_p} W(k)$ and $\mu_0=i_M^{-1}\mu i_M:\dbG_m\to \scrG_{W(k)}\leqslant \pmb{\GL}_{M_{\dbZ_p}\otimes_{\dbZ_p} W(k)}$. We have $\phi_0=g(1_{M_{\dbZ_p}}\otimes\sigma)\mu_0({1\over p})$, where $g\in \scrG_{\dbZ_p}(W(k))$. Let $\scrG^1_{\dbZ_p}$ be a reductive, closed subgroup scheme of $\scrG_{\dbZ_p}$ of the same rank and the same $\dbZ_p$-rank as $\scrG_{\dbZ_p}$ and whose adjoint group scheme is $\scrV_1$. By replacing $i_M$ with its composite with an automorphism of $M_{\dbZ_p}\otimes_{\dbZ_p} W(k)$ defined by an element of $\scrG_{\dbZ_p}(W(k))$, we can assume that the cocharacter $\mu_0$ factors through $\scrG^1_{W(k)}$, cf. [Bo, Ch. V, Thm. 19.2] and [DG, Vol. II, Exp. IX, Thms. 3.6 and 7.1]. Let $g_0\in \scrG_{\dbZ_p}^{\der}(W(k))$ be an element whose image in $\scrV_1(W(k))$ is trivial and such that we have $g_0g\in \scrG^1_{\dbZ_p}(W(k))$ (it exists as we have a short exact sequence $1\to \scrG^1_{\dbZ_p}\to\scrG_{\dbZ_p}\to \scrV_2/\scrT_2^0\to 1$ of flat sheaves on the category of $\Spec \dbZ_p$-schemes and as $G^{\der}/\scrT^0=\scrV_1/\scrT^0_1\times_{\dbZ_p} \scrV_2/\scrT^0_2$, where $\scrT^0$ is a maximal torus of $\scrG^{\der}_{\dbZ_p}$ contained in $\scrG^1_{\dbZ_p}$ and where $\scrT^0_i:=\text{Im}(\scrT^0\to \scrV_i)$). The triple $(M_{\dbZ_p}\otimes_{\dbZ_p} W(k),g_0\phi_0,\scrG^1_{W(k)})$ is a Shimura $F$-crystal over $k$, cf. the last two sentences. Both axioms 1.2.4 (b.i) and (b.ii) hold for $(M_{\dbZ_p}\otimes_{\dbZ_p} W(k),g_0\phi_0,\scrG^1_{W(k)})$. Argument: axiom 1.2.4 (b.i) holds as $\scrG^1_{\dbZ_p}$ is the closed subgroup scheme of $\scrG_{\dbZ_p}$ that fixes $Z^0(\scrG^1_{\dbZ_p})$ and axiom 1.2.4 (b.ii) holds as a maximal torus of $\scrG^1_{W(\bar k)}$ is a maximal torus of $\scrG_{W(\bar k)}$. 
Thus from many ``adjoint" points of view, (by taking $\scrV_1$ to be $\dbZ_p$-simple and by replacing $\scrC$ with $(M_{\dbZ_p}\otimes_{\dbZ_p} W(k),g_0\phi_0,\scrG^1_{W(k)})$) one can assume that $\scrG^{\ad}_{\dbZ_p}$ is simple. We will use this fact in Section 6. 

\medskip\noindent
{\bf 2.5.7. Fact.} {\it If $\scrC$ is Uni-ordinary, then $L_{\scrG}^0(\phi)$ is a reductive, closed subgroup scheme of $\scrG$.}

\medskip
\proof
By very definitions, $L_{\scrG}^0(\phi)_{B(k)}$ is a subgroup of $\scrP_{B(k)}$. The Lie algebra $\Lie(L_{\scrG}^0(\phi)_{B(\bar k)})$ is $B(\bar k)$-generated by elements which are fixed by $\phi\otimes\sigma_{\bar k}$ and which leave invariant $F^1[{1\over p}]$. Due to the functorial aspect of [Wi, p. 513], these elements and the $t_{\alpha}$'s are fixed by the inverse  $\mu_{\text{can}}$  of the canonical split cocharacter of $(M,F^1,\phi)$ defined in [Wi, p. 512]. Thus the generic fibre of $\mu_{\text{can}}$ factors through $Z^0(L_{\scrG}^0(\phi)_{B(k)})$ and therefore also through $\scrG_{B(k)}$. From this and the standard properties of canonical split cocharacters we get that $\mu_{\text{can}}$ is a Hodge cocharacter of $\scrC$. Thus we can assume that $\mu=\mu_{\text{can}}$; therefore $\mu_{B(k)}$ factors through $Z^0(L_{\scrG}^0(\phi)_{B(k)})$. Let $T^0$ be the image of $\mu$; it is a torus of the center of $L_{\scrG}^0(\phi)$. 

By induction on $i\in\dbN^*$ we get the existence of a unique torus $T^i$ of the center of $L_{\scrG}^0(\phi)$ such that we have $\Lie(T^i)=\phi^i(\Lie(T^0))$. Let $T_0$ be the torus of $\pmb{\GL}_M$ generated by the commuting tori $T^i$'s; it is a torus of the center of $L_{\scrG}^0(\phi)$ and thus also of $\scrG$. We claim that $L_{\scrG}^0(\phi)$ is the reductive group scheme $C_0:=C_{\scrG}(T_0)$. Obviously $L_{\scrG}^0(\phi)$ is a closed subgroup scheme of $C_0$. As $\mu$ factors through $Z^0(C_0)$ we have $\phi(\Lie(C_0))=\Lie(C_0)$. Thus $\Lie(C_{0,B(k)})\subseteq\Lie(L_{\scrG}^0(\phi)_{B(k)})$ i.e., $C_{0,B(k)}$ is a subgroup of $L_{\scrG}^0(\phi)_{B(k)}$ (cf. [Bo, Ch. II, Subsect. 7.1] and the fact that $C_{0,B(k)}$ and $L_{\scrG}^0(\phi)_{B(k)}$ are connected being reductive groups). Thus $L_{\scrG}^0(\phi)_{B(k)}=C_{0,B(k)}$ and therefore $L_{\scrG}^0=C_0$ is reductive.\endproof

\bigskip\noindent
{\bf 2.6. Reduction to the basic context.} Let $g\in \scrG(W(k))$ and $h\in \scrG(B(k))$ be such that $L_{\scrG}^0(g\phi)$ is a reductive, closed subgroup scheme of $\scrG$ through which $\mu:\dbG_m\to\scrG$ factors and we have an equality $hg\phi=\phi h$, cf. [Va7, Subsubsect. 3.1.1 and Thm. 3.1.2]. Thus by performing the operation $\grO_2$ (i.e., by replacing $\scrC$ with $(h(M),\phi,\scrG(h))$), in this subsubsection we will also assume that $L^0_{\scrG}(\phi)$ is a Levi subgroup scheme of $P^+_{\scrG}(\phi)$ and that $\mu$ factors through it (if $\scrC$ is Uni-ordinary, then this automatically holds, cf. Fact 2.5.7 and its proof). Thus $(M,F^1,\phi,L_{\scrG}^0(\phi))$ is a Shimura filtered $F$-crystal over $k$ and $(M,\phi,L_{\scrG}^0(\phi))$ is basic. Thus in connection to Conjecture 1.2.2 and to Subproblem 1.2.3, we can always replace $\scrG$ by a Levi subgroup scheme of $P_{\scrG}^+(\phi)$ (and accordingly we can assume that $\scrC$ is basic). However, often we will not perform this replacement in Sections 5 to 7, as by keeping track of $\scrG$ we often get extra information on $L^0_{\scrG}(\phi)$ as follows. Let $L_{\scrG}^0(\phi)_{\dbZ_p}$ be the reductive, closed subgroup scheme of $\scrG_{\dbZ_p}$ which is the $\dbZ_p$ structure of $L_{\scrG}^0(\phi)$ obtained as in the property 2.5 (b) (for $\square=L_{\scrG}^0(\phi)$). 

\medskip\noindent
{\bf 2.6.1. Fact.} {\it We recall $L^0_{\scrG}(\phi)$ is a Levi subgroup scheme of $P^+_{\scrG}(\phi)$ through which $\mu$ factors. 

\medskip
{\bf (a)} Then we have $L^0_{\scrG}(\phi)=C_{\scrG}(Z^0(L^0_{\scrG}(\phi)))$. Moreover, $L_{\scrG}^0(\phi)_{\dbZ_p}$ is the centralizer of a $\dbG_m^s$ subgroup scheme of $\scrG_{\dbZ_p}$ in $\scrG_{\dbZ_p}$ with $s\in\{0,1\}$ and $(M,\phi)$ is a direct sum of isoclinic Dieudonn\'e modules over $k$  (which have only one Newton polygon slope).

\smallskip
{\bf (b)} The standard axioms 1.2.4 (b.i) and (b.ii) hold for $(M,\phi,L_{\scrG}^0(\phi))$.}

\medskip
\proof
Both $L^0_{\scrG}(\phi)$ and $C_{\scrG}(Z^0(L^0_{\scrG}(\phi)))$ are reductive, closed subgroup schemes of $\scrG$ (cf. [DG, Vol. III, Exp. XIX, Subsect. 2.8] for $C_{\scrG}(Z^0(L^0_{\scrG}(\phi)))$). Thus they coincide if and only if their generic fibres coincide. But as $L^0_{\scrG}(\phi)\leqslant C_{\scrG}(Z^0(L^0_{\scrG}(\phi)))$, the fact that $L^0_{\scrG}(\phi)_{B(k)}=C_{\scrG}(Z^0(L^0_{\scrG}(\phi)))_{B(k)}$ follows from the fact that $L^0_{\scrG}(\phi)_{B(k)}$ is the centralizer of the cocharacter $\nu_{B(k)}$ of $\scrG_{B(k)}$ which factors through $Z^0(L^0_{\scrG}(\phi)_{B(k)})$ and which is the Newton cocharacter of $\scrC$ (see [Va7, Subsect. 2.3]). The cocharacter $\nu_{B(k)}$ is fixed by $\phi$ and $\mu$ and thus it is the extension to $B(k)$ of a cocharacter $\nu$ of $Z^0(L^0_{\scrG}(\phi)_{\dbQ_p})$. As $Z^0(L^0_{\scrG}(\phi))$ is a torus, $\nu$ extends to a cocharacter of $Z^0(L^0_{\scrG}(\phi)_{\dbZ_p})$. Its centralizer in $\scrG_{\dbZ_p}$ is a reductive, closed subgroup scheme of $\scrG_{\dbZ_p}$ (cf. [DG, Vol. III, Exp. XIX, Subsect. 2.8]) whose generic fibre is $L^0_{\scrG}(\phi)_{\dbQ_p}$ and therefore it is $L^0_{\scrG}(\phi)_{\dbZ_p}$ itself. The image of this cocharacter of $Z^0(L^0_{\scrG}(\phi)_{\dbZ_p})$ is isomorphic to $\dbG_m^s$  with $s\in\{0,1\}$. As $\nu$ extends to a cocharacter of $Z^0(L^0_{\scrG}(\phi)_{\dbZ_p})$, $\nu_{B(k)}$ extends also to a cocharacter of $Z^0(L^0_{\scrG}(\phi))$. This implies that $(M,\phi)$ is a direct sum of isoclinic Dieudonn\'e modules over $k$. Thus (a) holds.

Axiom 1.2.4 (b.i) (resp. (b.ii)) holds for $(M,\phi,L_{\scrG}^0(\phi))$ due to the second part of (a) (resp. due to the fact that a maximal torus of $L_{\scrG}^0(\phi)$ is as well a maximal torus of $\scrG$).\endproof

\bigskip\noindent
{\bf 2.7. Proof of a general variant form of Theorem 1.2.5 (a) in the case with many endomorphisms.} We assume that $D$ is the $p$-divisible group of an abelian variety $A$ over $k$ and that there exists a polarization $\lambda_A$ of $A$ whose crystalline realization has a $W(k)$-span normalized by $\scrG$. We also assume that we have an equality
$$\{x\in\Lie(\scrG)[{1\over p}]|\phi(x)=x\}=\{x\in\Lie(\scrG)\otimes_{W(k)} B(\dbF)|\phi\otimes\sigma_{\dbF}(x)=x\}$$
which we express by saying that $\scrC$ is with many endomorphisms.
We check that, up to operations $\grO_1$ and $\grO_2$, there exists an abelian scheme $A_{W(k)}$ over $W(k)$ which is a lift of $A$ with respect to $\scrG$ (i.e., conditions 1.2 (b) to 1.2 (d) hold with $(V,\scrG_V^\prime)=(W(k),\scrG)$) and to which the Frobenius endomorphism of $A$ lifts.

To prove this, we can assume that $\scrG=L^0_{\scrG}(\phi)$ and that $\scrC$ is basic (cf. Subsection 2.6). As $\scrC$ is with many endomorphisms, $\gre:=\{x\in\Lie(\scrG)[{1\over p}]|\phi(x)=x\}\otimes_{\dbQ_p} B(k)$ is $\Lie(\scrG_{B(k)})$. Thus $\phi^r\in\scrG(B(k))$ fixes $\Lie(\scrG_{B(k)})$ and therefore we have $\phi^r\in Z^0(\scrG)(B(k))$. As $(M,\phi)$ is a direct sum of isoclinic Dieudonn\'e modules over $k$ (cf. Fact 2.6.1 (a)) and as $\mu$ commutes with $Z^0(\scrG)$ and thus with the cocharacter of $Z^0(\scrG)$ which extends the Newton cocharacter $\nu_{B(k)}$ (see proof of Fact 2.6.1), $(M,F^1,\phi)$ is a direct sum of filtered Dieudonn\'e modules whose associated Dieudonn\'e modules are isoclinic. Thus even if $p=2$, there exists a $p$-divisible group $D_{W(k)}$ over $W(k)$ which lifts $D$, which is a direct sum of $p$-divisible groups over $W(k)$ whose special fibres are isoclinic, and whose filtered $F$-crystal over $k$ is $(M,F^1,\phi)$. Due to the existence of the polarization $\lambda_A$, it is easy to see that there exists an abelian scheme $A_{W(k)}$ over $W(k)$ which lifts $A$ in such a way that its $p$-divisible group is $D_{W(k)}$ (cf. Serre--Tate deformation theory and Grothendieck's algebraization theorem). As $\phi^r\in Z^0(\scrG)(B(k))$ maps $F^1$ to $F^1$ and is the crystalline realization of the Frobenius endomorphism of $A$ and due to the direct sum property of $D_{W(k)}$, (even if $p=2$) the Frobenius endomorphism of $A$ lifts to an endomorphism of $A_{W(k)}$ (cf. Serre--Tate and Messing--Grothendieck deformation theories; it is an endomorphism and not a $\dbZ_p$-endomorphism of $A_{W(k)}$ due to the existence of $\lambda_A$).

\bigskip\noindent
{\bf 2.8. Some $W(k)$-algebras.} 
Let $e\in\dbN^*$. Let $X$ be an independent variable. Let $R:=W(k)[[X]]$. Let $\Rtil e$ (resp. $Re$) be the $W(k)$-subalgebra of $B(k)[[X]]$ formed by formal power series $\sum_{n=0}^\infty a_nX^n$ for which we have $a_n{\fracwithdelims[]{n}{e}}!\in W(k)$ for all $n$ (resp. for which the sequence $b_n:=a_n{\fracwithdelims[]{n}{e}}!$ is formed by elements of $W(k)$ and converges to 0). Thus $Re$ is a $W(k)$-subalgebra of $\Rtil e$. Let $\Phi_R$, $\Phi_{Re}$, and $\Phi_{\tilde Re}$ be the Frobenius lifts of $R$, $Re$, and $\Rtil e$ (respectively) that are compatible with $\sigma$ and that take $X$ to $X^p$. For $m\in\dbN^*$ let $I(m)$ be the ideal of $\tilde Re$ formed by formal power series with $a_0=a_1=\cdots=a_{m-1}=0$. By mapping $X$ to $0$ we get $W(k)$-epimorphisms $R\twoheadrightarrow W(k)$, $Re\twoheadrightarrow W(k)$, and $\Rtil e\twoheadrightarrow W(k)$ that respect the Frobenius lifts (the kernel of $\Rtil e\twoheadrightarrow W(k)$ is $I(1)$). The proof of the following elementary fact is left as an exercise. 

\medskip\noindent
 {\bf 2.8.1. Fact.} {\it We assume that $p\ge 3$ (resp. $p=2$). Let $V$ be a finite, totally ramified discrete valuation ring extension of $W(k)$ of degree at most $e$. Let $\pi_V$ be a uniformizer of $V$. Then there exist $W(k)$-epimorphisms $R\twoheadrightarrow V$, $Re\twoheadrightarrow V$, and $\Rtil e\twoheadrightarrow V$ (resp. $R\twoheadrightarrow V$ and $Re\twoheadrightarrow V$) that map $X$ to $\pi_V$.}
\finishproclaim

\bigskip\smallskip
\noindent
{\boldsectionfont 3. Unramified and ramified $\text{CM}$-isogeny classifications}

\bigskip
Let $\scrF=\{\scrC_g|g\in\scrG(W(k))\}$ and $\scrY(\scrF):=\cup_{g\in \scrG(W(k))} \scrI(\scrC_g)$ be as in Subsection 2.1. By the {\it strong $\text{CM}$-isogeny} (resp. by the {\it $\text{CM}$-isogeny}) {\it classification} of $\scrF$ we mean the description of the subset $SZ(\scrY(\scrF))$ (resp. $Z(\scrY(\scrF))$) of $\scrY(\scrF)$ formed by the inner isomorphism classes of those $\scrC_g$ with $g\in \scrG(W(k))$ which, up to the operation $\grO_2$ (resp. up to the operations $\grO_1$ and $\grO_2$), have a lift of quasi $\text{CM}$ type. Fact 3.1 collects few simple properties. These $\text{CM}$-isogenies classifications are too restrictive and difficult to be accomplished and thus in Subsection 3.2 we also introduce ramified lifts of $\scrC$ (or of $D$ with respect to $\scrG$) and the {\it (strong) ramified $\text{CM}$-isogeny classification} of $\scrF$. 

In Subsections 3.3 to 3.7 we include different properties required in Sections 5 to 8 and some remarks. In particular, Corollary 3.6.3 checks that for $p\ge 3$  the ramified lifts of $D$ with respect to $\scrG$ (see Definition 3.6.2) are in natural bijection to abstract ramified lifts of $\scrC$ (see Definition 3.2.1 (d)). We will use the notations of Subsections 2.1 and 2.5. We recall that the axioms 1.2.4 (b.i) and (b.ii) hold. 

\bigskip\noindent
{\bf 3.1. Fact.} {\bf (a)} {\it If $\scrC$ has a lift of quasi $\text{CM}$ type, then it is unramified semisimple.

\smallskip
{\bf (b)} If there exists a maximal torus of $\scrG_{B(k)}$ of $\dbQ_p$-endomorphisms of $\scrC$, then $\scrC$ is semisimple.

\smallskip
{\bf (c)} If $\scrC$ is semisimple, then by performing operation $\grO_1$ we can assume that there exists a maximal torus of $\scrG_{B(k)}$ of $\dbQ_p$-endomorphisms of $\scrC$.}

\medskip
\proof
We prove (a). Let $\scrT$ be a maximal torus of $\scrG$ such that we have $\phi(\Lie(\scrT))=\Lie(\scrT)$. Thus the element $\phi^r\in \scrG(B(k))$ (see Subsubsection 2.5.1) normalizes $\scrT$. Therefore we have $\phi^r\in N_{\scrG}(\scrT)(B(k))$. Let $m\in\dbN^*$ be such that $\phi^{rm}\in \scrT(B(k))$. As the torus $\scrT_{W(\bar k)}$ is split, part (a) follows. The proof of (b) is very much the same and thus left as an exercise. 

We prove (c). By performing the operation $\grO_1$, we can assume that the centralizer of the semisimple element $\phi^r\in \scrG(B(k))$ in $\scrG_{B(k)}$ is a reductive subgroup $C_2$. The Lie algebra $\Lie(Z^0(C_2))$ (or $\Lie(C_2^{\ad})$) is normalized by $\phi$ and fixed by $\phi^r$ and thus it is $B(k)$-generated by elements fixed by $\phi$. Let $C^{\ad}_{2,\dbQ_p}$ be the adjoint group over $\dbQ_p$ whose Lie algebra is $\{x\in\Lie(C_2^{\ad})|\phi(x)=x\}$; its extension to $B(k)$ is $C_2^{\ad}$ and it is constructed as the identity component of the group of Lie automorphisms of $\{x\in\Lie(C_2^{\ad})|\phi(x)=x\}$.  Let $T_{2,\dbQ_p}$ be a maximal torus of $C^{\ad}_{2,\dbQ_p}$. If $\scrT_{1,B(k)}$ is the maximal torus of $C_2$ (and thus of $\scrG_{B(k)}$) whose image in $C_2^{\ad}$ is $T_{2,B(k)}$, then $\Lie(\scrT_{1,B(k)})=\Lie(Z^0(C_2))\oplus \Lie(T_{2,B(k)})$ is $B(k)$-generated by elements fixed by $\phi$ and thus (c) holds.\endproof

\bigskip\noindent
{\bf 3.2. The ramified context.} Let $V$ be a finite, totally ramified discrete valuation ring extension of $W(k)$. Let $K:=V[{1\over p}]$ and let $\pi_V$ be a uniformizer of $V$. We assume that $e:=[V:W(k)]\ge 2$. Let the $W(k)$-algebras $R$, $Re$, and $\tilde Re$ be as in Subsection 2.8. For $m\in\dbN^*$ let $\Phi_{R_m}$ be the Frobenius lift of $R_m:=W(k)[[X_1,\ldots,X_m]]$ which is compatible with $\sigma$ and which takes $X_i$ to $X_i^p$ for all $i\in\{1,\ldots,m\}$. If $m=1$ we drop it as an index; thus $R_1=R$. The $(X_1,\ldots,X_m)$-adic completion $\hat{\Omega}_{R_m}$ of $\Omega_{R_m}$ is a free $R_m$-module that has $\{dX_1,\ldots,dX_m\}$ as an $R_m$-basis. Let $\chi_e:R\twoheadrightarrow V$ be the $W(k)$-epimorphism which takes $X$ to $\pi_V$. If $p\Ge 2$ (resp. $p\Ge 3$), we denote also by $\chi_e$ the $W(k)$-epimorphism $Re\twoheadrightarrow V$ (resp. $\tilde Re\twoheadrightarrow V$) defined by $\chi_e$ (cf. Fact 2.8.1). 

\medskip\noindent
{\bf 3.2.1. Definitions.} {\bf (a)} By a {\it lift of $\scrC$ to $R_m$} we mean a quadruple
$$(M_{R_m},F^1_{R_m},\phi_{M_{R_m}},\scrG_{R_m}),$$
where $M_{R_m}$ is a free $R_m$-module of the same rank as $M$, $F^1_{R_m}$ is a direct summand of $M_{R_m}$, $\scrG_{R_m}$ is a reductive, closed subgroup scheme of $\pmb{\GL}_{M_{R_m}}$, and $\phi_{M_{R_m}}:M_{R_m}\to M_{R_m}$ is a $\Phi_{R_m}$-linear endomorphism, such that the following three axioms hold:

\medskip
{\bf (i)} the $R_m$-module $M_{R_m}$ is generated by $\phi_{M_{R_m}}(M_{R_m}+p^{-1}F^1_{R_m})$;

\smallskip
{\bf (ii)} there exists a family of tensors $(t^{R_m}_{\alpha})_{\alpha\in\scrJ}$ of the $F^0$-filtration of $\scrT(M_{R_m})$ defined by $F^1_{R_m}$ such that we have $\phi_{M_{R_m}}(t^{R_m}_{\alpha})=t^{R_m}_{\alpha}$ for all $\alpha\in\scrJ$ and $\scrG_{R_m}[{1\over p}]$ is the closed subscheme of $\pmb{\GL}_{M_{R_m}[{1\over p}]}$ that fixes $t^{R_m}_{\alpha}$ for all $\alpha\in\scrJ$;

\smallskip
{\bf (iii)} the extension of $(M_{R_m},\phi_{M_{R_m}},\scrG_{R_m})$ via the $W(k)$-epimorphism $\chi_0:R_m\twoheadrightarrow W(k)$ that maps each $X_i$ to $0$, is $\scrC$.

\medskip
{\bf (b)} Let $\chi_{m;e}:R_m\twoheadrightarrow V$ be a $W(k)$-epimorphism; if $m=1$ we take $\chi_{1;e}:=\chi_e$. Let $F^1_V:=F^1_{R_m}\otimes_{R_m} {}_{\chi_{m;e}}V$. We refer to $(M_{R_m},F^1_V,\phi_{M_{R_m}},\scrG_{R_m})$ as a {\it lift of $\scrC$ to $R_m$ with respect to $V$}. 

\smallskip
{\bf (c)} We say $(M_{R_m},F^1_V,\phi_{M_{R_m}},\scrG_{R_m})$ is a {\it lift of $\scrC$ to $R_m$ of quasi $\text{CM}$ (resp. of $\text{CM}$) type with respect to $V$}, if there exists a maximal torus $\scrT_{R_m[{1\over p}]}$ of $\scrG_{R_m[{1\over p}]}$ such that the following two axioms hold:

\medskip
{\bf (i)} $\Lie(\scrT_{R_m[{1\over p}]})$ is left invariant (resp. is $R_m[{1\over p}]$-generated by elements fixed) by $\phi_{M_{R_m}}$;

\smallskip
{\bf (ii)} $F^1_V[{1\over p}]$ is a $\Lie(\scrT_{R_m[{1\over p}]})/\Ker(\chi_{m;e})\Lie(\scrT_{R_m[{1\over p}]})$-module.

\medskip
{\bf (d)} We have variants of (a) to (c), where we replace $R_m$ by $Re$ if $p\ge 2$ or by $\Rtil e$  if $p\ge 3$ (the $W(k)$-epimorphisms from either $Re$ or $\Rtil e$ onto $V$ being $\chi_e$). A lift of $\scrC$ to $\Rtil e$ (which is of quasi $\text{CM}$ type or of $\text{CM}$ type) with respect to $V$ is also called a {\it ramified lift of $\scrC$ to $V$ (of quasi $\text{CM}$ type or of $\text{CM}$ type)}.

\smallskip
{\bf (e)} If we have a principal bilinear quasi-polarization $\lambda_M:M\times M\to W(k)$ of $\scrC$, then by a {\it lift of $(\scrC,\lambda_M)$ to $R_m$} we mean a quintuple
$$(M_{R_m},F^1_{R_m},\phi_{M_{R_m}},\scrG_{R_m},\lambda_{M_{R_m}}),$$ 
where $(M_{R_m},F^1_{R_m},\phi_{M_{R_m}},\scrG_{R_m})$ is as in (a) and
$\lambda_{M_{R_m}}$ is a perfect bilinear form on $M_{R_m}$ which lifts $\lambda_M$, whose $R_m$-span is normalized by $\scrG_{R_m}$, and for which we have an identity $\lambda_M(\phi_{M_{R_m}}(x),\phi_{M_{R_m}}(y))=p\Phi_{R_m}(\lambda_{M_{R_m}}(x,y))$ for all elements $x,y\in M_{R_m}$. Similarly, definitions (b) to (d) extend to the principal bilinear quasi-polarized context.

\medskip\noindent
{\bf 3.2.2. Remarks.} {\bf (a)} Let $(M_{R_m},F^1_V,\phi_{M_{R_m}},\scrG_{R_m})$ be a lift of $\scrC$ to $R_m$ of $\text{CM}$ type with respect to $V$. Let $\scrT_{R_m[{1\over p}]}$ be a maximal torus of $\scrG_{R_m[{1\over p}]}$ such that the two axioms of Definition 3.2.1 (c) hold for it. Let $\scrT_{1,B(k)}$ be the pull-back of $\scrT_{R_m[{1\over p}]}$ via the $B(k)$-epimorphism $R_m[{1\over p}]\twoheadrightarrow B(k)$ that takes each $X_i$ to $0$. As $\Lie(\scrT_{R_m[{1\over p}]})$ is $R_m[{1\over p}]$-generated by elements fixed by $\phi_{M_{R_m}}$, $\scrT_{1,B(k)}$ is a maximal torus of $\scrG_{B(k)}$ of $\dbQ_p$-endomorphisms of $\scrC$. 
Similarly, if $\scrC$ has a ramified lift to $V$ of $\text{CM}$ type, then there exist maximal tori of $\scrG_{B(k)}$ of $\dbQ_p$-endomorphisms of $\scrC$.

\smallskip
{\bf (b)} We refer to Definition 3.2.1 (a). The reductive group scheme $\scrG_{R_m}$ over $R_m$ lifts $\scrG$ (cf. axiom (iii) of Definition 3.2.1 (a)) and thus, as $R_m$ is complete in the $(X_1,\ldots,X_m)$-adic topology, it is isomorphic to $\scrG\times_{W(k)} \Spec R_m$ (i.e., our notations match).

\smallskip
{\bf (c)} The recent works [VZ] and [Lau1,2] (resp. the work [Ki1]) would allow variants of the Definitions 3.2.1 (a) to (e) (resp. Definition 3.2.1 (d)) in terms of {\it Breuil windows}. However, as [VZ] and [Lau1,2] are stated only in terms of $p$-divisible groups and not in terms of suitable crystalline representations, one would encounter some technical difficulties in connection to the axioms (ii) of Definition 3.2.1 (a) and (i) of Definition 3.2.1 (c).

\medskip\noindent
{\bf 3.2.3. Ramified $\text{CM}$-isogeny classifications.} By the {\it strong ramified} (resp. by the {\it ramified}) {\it $\text{CM}$-isogeny classification} of $\scrF$ we mean the description of the subset
$$SZ^{\text{ram}}(\scrY(\scrF))\,\, (\text{resp.}\,\,Z^{\text{ram}}(\scrY(\scrF)))$$
 of $\scrY(\scrF)$ formed by inner isomorphism classes of those $\scrC_g$ with $g\in \scrG(W(k))$ for which, up to the operation $\grO_2$ (resp. up to operations $\grO_1$ and $\grO_2$), there exists a discrete valuation ring $V$ as in Subsection 3.2 and a ramified lift of $\scrC_g$ to $V$ of quasi $\text{CM}$ type. We have $Z(\scrY(\scrF))\subseteq Z^{\text{ram}}(\scrY(\scrF))$ and $SZ(\scrY(\scrF))\subseteq SZ^{\text{ram}}(\scrY(\scrF))$.

\bigskip\noindent
{\bf 3.3. Connections.} We refer to Definition 3.2.1 (a). If $p=2$ we assume that there exists a $2$-divisible group over $W(k)$ which lifts $D$ and whose Hodge filtration is the direct summand $F^1_{R_m}/\Ker(\chi_0)F^1_{R_m}$ of $M_{R_m}/\Ker(\chi_0)F^1_{R_m}=M$ (this condition is not necessary but this is the context in which [Fa, Thm. 10] is stated). From either loc. cit.  or [Va8, Thms. 3.2 and 3.3.1] we get the existence of a unique connection $\nabla_m:M_{R_m}\to M_{R_m}\otimes_{R_m} \hat{\Omega}_{R_m}$ on $M_{R_m}$ such that $\phi_{M_{R_m}}$ is horizontal with respect to it (i.e., the identity $\nabla_m\circ \phi_{M_{R_m}}=(\phi_{M_{R_m}}\otimes d\Phi_{R_m})\circ\nabla_m$ holds); it is integrable and nilpotent modulo $p$. There exists a unique (up to a unique isomorphism) $p$-divisible group over $R_m/pR_m$ whose $F$-crystal is $(M_{R_m},\phi_{M_{R_m}},\nabla_m)$  (cf. loc. cit. for the existence and cf. [BM, Thm. 4.1.1] for the uniqueness).

\bigskip\noindent
{\bf 3.4. Theorem.} {\it We assume that $p\ge 3$ and $\scrG=\pmb{\GL}_M$. Then the ramified lifts of $\scrC$ to $V$ are in natural bijection to lifts of $D$ to $p$-divisible groups over $V$.}

\medskip
\proof 
To a $p$-divisible group $D_V$ over $V$ that lifts $D$ one associates uniquely a ramified lift of $\scrC$ to $V$ as follows. Let
$$(M_{\tilde Re},\phi_{M_{\tilde Re}},\nabla)$$ 
be the extension via the $W(k)$-monomorphism $Re\hookrightarrow\tilde Re$ of (the projective limit indexed by $n\in\dbN^*$ of the evaluation at the thickening naturally attached to the closed embedding $\Spec V/pV\hookrightarrow\Spec Re/p^nRe$ of) the Dieudonn\'e $F$-crystal over $V/pV$ of $D_V\times_V \Spec V/pV$ (see [Me], [BBM], [BM], and [dJ, Subsect. 2.3]). Thus $M_{\tilde Re}$ is a free $\tilde Re$-module of the same rank as $D$, $\nabla:M_{\tilde Re}\to M_{\tilde Re}\otimes_{\tilde Re} \tilde RedX$ is
an integrable and nilpotent modulo $p$ connection on $M_{\tilde Re}$, and
$\phi_{M_{\tilde Re}}$ is a $\Phi_{\tilde Re}$-linear endomorphism of $M_{\tilde Re}$ which is horizontal with respect to $\nabla$ (its Verschiebung map is uniquely determined and therefore it is not added to the notations). If $F^1(M_{\tilde Re})$ is the inverse image in $M_{\tilde Re}$ of the Hodge filtration $F^1_V$ of $M_{\tilde Re}/\Ker(\chi_e)M_{\tilde Re}=H_{\text{dR}}^1(D_V/V)$ defined by $D_V$, then the
restriction of $\phi_{M_{\tilde Re}}$ to $F^1(M_{\tilde Re})$ is divisible by $p$ and $M_{\tilde Re}$ is $\tilde Re$-generated by $\phi_{M_{\tilde Re}}(M_{\tilde Re}+{1\over p}F^1(M_{\tilde Re}))$.

The quadruple $(M_{\tilde Re},F^1_V,\phi_{M_{\tilde Re}},\pmb{\GL}_{M_{\tilde Re}})$ is the ramified lift of $\scrC$ to $V$ associated to $D_V$. Due to the fully faithfulness part of [Fa, Thm. 5], to prove the theorem it suffices to show that every ramified lift of $\scrC$ to $V$ is associated to a $p$-divisible group over $V$ which lifts $D$. As $p\Ge 3$, each lift of $(M,\phi,\pmb{\GL}_M)$ is associated to (i.e., it is the filtered $F$-crystal over $k$ of) a unique $p$-divisible group over $W(k)$ that lifts $D$ (cf. Grothendieck--Messing deformation theory of [Me, Chs. IV and V]). Thus each lift of $\scrC$ to $R_m$ is associated to a unique $p$-divisible group over $R_m$ that lifts $D$, cf. end of Subsection 3.3. Therefore the theorem follows from the following general result (applied with $p\Ge 3$, $\scrG=\pmb{\GL}_M$, and $\scrJ=\emptyset$).\endproof

\bigskip\noindent
{\bf 3.5. Theorem.} {\it We assume that $p\ge 3$ but do not assume that $\scrG$ is $\pmb{\GL}_M$. We take $m$ to be $\dim(\scrG_{B(k)}^{\der})$. Then each lift $(M_{\tilde Re},F^1_{\tilde Re},\phi_{M_{\tilde Re}},\scrG_{\tilde Re})$ of $\scrC$ to $\tilde Re$ is the extension via a $W(k)$-homomorphism $R_m\to \Rtil e$ of a lift of $\scrC$ to $R_m$ (this makes sense, cf. the uniqueness of the connection $\nabla_m$ of Subsection 3.3).}

\medskip
\proof
Let $(t_{\alpha}^{\tilde Re})_{\alpha\in\scrJ}$ be a family of tensors of the $F^0$-filtration of $\scrT(M_{\tilde Re})$ defined by $F^1_{\tilde Re}$ which has the analogue meaning of the family of tensors $(t_{\alpha}^{R_m})_{\alpha\in\scrJ}$ of Definition 3.2.1 (a). Not to introduce extra notations, we can assume that $(M,F^1,\phi,(t_{\alpha})_{\alpha\in\scrJ})$ is the extension of $\scrC^{\tilde Re}:=(M_{\tilde Re},F^1_{\tilde Re},\phi_{M_{\tilde Re}},(t_{\alpha}^{\tilde Re})_{\alpha\in\scrJ})$ via the $W(k)$-epimorphism $\chi_e:\tilde Re\twoheadrightarrow W(k)$. Let 
$$\scrC^{R_m}:=(M\otimes_{W(k)} R_m,F^1\otimes_{W(k)} R_m,g^{\der}_{\text{univ}}(\phi\otimes\Phi_{R_m}),(t_{\alpha})_{\alpha\in\scrJ}),$$ 
where $g^{\der}_{\text{univ}}:\Spec R_m\to \scrG^{\der}$ is a universal morphism that identifies $\Spf R_m$ with the formal completion of $\scrG^{\der}$ along its identity section. Let $\nabla_{\text{univ}}$ be the unique connection on $M\otimes_{W(k)} R_m$ such that $g^{\der}_{\text{univ}}(\phi\otimes\Phi_{R_m})$ is horizontal with respect to it, cf. Subsection 3.3. Let $\delta$ be the flat connection on $M\otimes_{W(k)} R_m$ that annihilates $M\otimes 1$. 

We have $\omega_{\text{univ}}:=\nabla_{\text{univ}}-\delta\in \Lie(\scrG)\otimes_{W(k)} \hat{\Omega}_{R_m}$, cf. either [Fa2, Sect. 7, Rm. ii)] or [Va8, Thms. 3.2 and 3.3.1]. For the sake of completeness,  we include a proof of this result. We view $\scrT(M)$ as a module over the Lie algebra (associated to) $\End(M)$ and we denote also by $\nabla_{\text{univ}}$ the connection on $\scrT(M\otimes_{W(k)} R_m[{1\over p}])$ which extends naturally the connection $\nabla_{\text{univ}}$ on $M\otimes_{W(k)} R_m$. Each tensor $t_{\alpha}\in\scrT(M\otimes_{W(k)} R_m[{1\over p}])$ is fixed under the natural action of $g^{\der}_{\text{univ}}(\phi\otimes\Phi_{R_m})$ on $\scrT(M\otimes_{W(k)} R_m[{1\over p}])$. Thus we have $\nabla_{\text{univ}}(t_{\alpha})=(g^{\der}_{\text{univ}}(\phi\otimes\Phi_{R_m})\otimes d\Phi_{R_m})(\nabla_{\text{univ}}(t_{\alpha}))$. As we have $d\Phi_{R_m}(X_i)=pX_i^{p-1}dX_i$ for each $i\in\{1,\ldots,m\}$, by induction on $s\in\dbN^*$ we get that $\nabla_{\text{univ}}(t_{\alpha})=\omega_{\text{univ}}(t_{\alpha})\in \scrT(M)\otimes_{W(k)} (X_1,\ldots,X_m)^s \hat{\Omega}_{R_m}[{1\over p}]$. As $R_m$ is complete with respect to the $(X_1,\ldots,X_m)$-topology, we get that $\nabla_{\text{univ}}(t_{\alpha})=\omega_{\text{univ}}(t_{\alpha})=0$. But $\Lie(\scrG_{B(k)})\cap\End(M)$ is the Lie subalgebra of $\End(M)$ that annihilates $t_{\alpha}$ for all $\alpha\in\scrJ$. From the last two sentences we get that $\omega_{\text{univ}}\in \Lie(\scrG)\otimes_{W(k)} \hat{\Omega}_{R_m}$. 

Next we list three basic properties of the $W(k)$-algebra $\tilde Re$:

\medskip
{\bf (i)} we have $\tilde Re=\text{proj.lim.}_{m\in\dbN^*} \tilde Re/I(m)$, the transition $W(k)$-epimorphisms being the logical ones (see Subsection 2.8 for the ideals $I(m)$);

\smallskip
{\bf (ii)} the $W(k)$-module $I(m)/I(m+1)$ is free of rank $1$ for all $m\in\dbN^*$;

\smallskip
{\bf (iii)} we have an inclusion $I(m)^2+\Phi_{\tilde Re}(I(m))\subseteq I(m+1)$ for all $m\in\dbN^*$.
\medskip

Due to these properties, the arguments of either [Fa, Thm. 10 and Rm. (iii) of p. 136] or [Va8, Thm. 5.3] apply entirely to give us that $\scrC^{\tilde Re}$ is the extension of $\scrC^{R_m}$ through a $W(k)$-homomorphism $R_m\to\tilde Re$ that maps the ideal $(X_1,\ldots,X_m)$ to $I(1)$ (this extension is well defined as the connection $\nabla_{\text{univ}}$ exists and is unique). More precisely, by induction on $m\in\dbN^*$ one constructs a $W(k)$-homomorphism $R_m\to \tilde Re/I(m)$ that maps the ideal $(X_1,\ldots,X_m)$ to $I(1)/I(m)$ and such that $\scrC^{\tilde Re}$ modulo $I(m)$ is the extension of $\scrC^{R_m}$ through it. Strictly speaking, loci citi are stated in terms of a universal element of either $\scrG$ or of some closed, unipotent, smooth subgroup scheme $U$ of $\scrG^{\der}$ and not of $\scrG^{\der}$. But the image $\grH$ of the Kodaira--Spencer map of $\nabla_{\text{univ}}$ is the same regardless if we work with $\scrG^{\der}$, $U$, or $\scrG$ and thus each one of the loci citi apply in our present context of $\scrG^{\der}$ as well. More precisely, if we write $\omega_{\text{univ}}=\sum_{i=1}^m \omega_i\otimes dX_i$ with $\omega_i\in \Lie(\scrG)\otimes_{W(k)} R_m$, then $\grH$ is the image of $\sum_{i=1}^m R_m\omega_i$ in $\Hom_{W(k)}(F^1,M/F^1)\otimes_{W(k)} R_m$ and it is identified naturally with $\tilde F^{-1}(\Lie(\scrG^{\der}))\otimes_{W(k)} R_m$, where $\tilde F^{-1}(\Lie(\scrG^{\der}))$ is the maximal direct summand of either $\Lie(\scrG^{\der})$ or $\Lie(\scrG)$ on which a fixed Hodge cocharacter of $\scrC$ acts via the weight $-1$ (if the fixed Hodge cocharacter of $\scrC$ is the inverse of the canonical split cocharacter of $(M,F^1,\phi)$ defined in [Wi], then the mentioned closed, unipotent, smooth subgroup scheme $U$ of $\scrG^{\der}$ is defined uniquely by the requirement that $\Lie(U)=\tilde F^{-1}(\Lie(\scrG^{\der}))$; to be compared with [Va8, Subsect. 2.6 and Subsubsect. 5.2.7]).\endproof

\medskip\noindent
{\bf 3.5.1. Corollary.} {\it We assume that $p\ge 3$. Let $(M_{\tilde Re},F^1_{\tilde Re},\phi_{M_{\tilde Re}},\scrG_{\tilde Re})$ be a lift of $\scrC$ to $\tilde Re$. Let $D_V$ be the $p$-divisible group over $V$ that lifts $D$ and that corresponds to $(M_{\tilde Re},F^1_{\tilde Re}\otimes_{\tilde Re} {}_{\chi_e}V,\phi_{M_{\tilde Re}})$ via the natural bijection of Theorem 3.4. Let $\scrT_{\tilde Re[{1\over p}]}$ be a maximal torus of $\scrG_{\tilde Re[{1\over p}]}$ such that $\phi_{M_{\tilde Re}}$ leaves invariant $\Lie(\scrT_{\tilde Re[{1\over p}]})$ and $F^1_K:=F^1_{\tilde Re}\otimes_{\tilde Re[{1\over p}]} {}_{\chi_e[{1\over p}]} K$ is a $\Lie(\scrT_{\tilde Re[{1\over p}]})\otimes_{\tilde Re[{1\over p}]} K$-module. Then the following two properties hold:

\medskip
{\bf (a)} By performing the operation $\grO_1$ we can assume that $\Lie(\scrT_{\tilde Re[{1\over p}]})$ is $\tilde Re[{1\over p}]$-generated by elements fixed by $\phi_{M_{\tilde Re}}$ and thus that $\scrC$ has ramified lifts to $V$ of $\text{CM}$ type. 

\smallskip
{\bf (b)} Up to the operation $\grO_1$, the $p$-divisible group $D_V$ is with complex multiplication.}

\medskip
\proof
There exists a canonical and functorial (in $D_V$) identification 
$$(M\otimes_{W(k)}\tilde Re[{1\over p}],\phi\otimes\Phi_{\tilde Re})=(M_{\tilde Re}[{1\over p}],\phi_{M_{\tilde Re}})\leqno (4)$$
 under which the pull-back of the natural $B(k)$-epimorphism $\tilde Re[{1\over p}]\twoheadrightarrow B(k)$ that takes $X$ to $0$ is the identity automorphism of $(M[{1\over p}],\phi)$ (see [Fa, Sect. 6] for the existence part; the uniqueness part follows from the fact that no element of $\End(M[{1\over p}]\otimes_{B(k)} I(1)$ is fixed by $\phi\otimes\Phi_{\tilde Re}$). Therefore via (4), we can identify $F^1_K$ with a direct summand of $M\otimes_{W(k)} K=H^1_{\text{dR}}(D_V/V)[{1\over p}]=[M_{\tilde Re}/\Ker(\chi_e)M_{\tilde Re}][{1\over p}]$. 

Let $\scrT_{1,B(k)}$ be the pull-back of $\scrT_{\tilde Re[{1\over p}]}$ to a maximal torus of $\scrG_{B(k)}$ (via the mentioned $B(k)$-epimorphism $\tilde Re[{1\over p}]\twoheadrightarrow B(k)$). Thus  $\phi$ normalizes $\Lie(\scrT_{1,B(k)})$ and under the identification (4), $\Lie(\scrT_{1,B(k)})\otimes_{B(k)} \tilde Re[{1\over p}]$ gets identified with $\Lie(\scrT_{\tilde Re[{1\over p}]})$. As $\phi^r\in\scrG(B(k))$ normalizes $\Lie(\scrT_{1,B(k)})$ and thus $\scrT_{1,B(k)}$, by performing the operation $\grO_1$ we can assume that $\phi^r$ fixes $\Lie(\scrT_{1,B(k)})$. Thus $\Lie(\scrT_{1,B(k)})$ is $B(k)$-generated by elements fixed by $\phi$, cf. 2.5 (a). 

Let $\tilde t\in\Lie(\scrT_{\tilde Re[{1\over p}]})$ be an element which lifts an element $t\in\Lie(\scrT_{1,B(k)})$ fixed by $\phi$. As $\Phi_{\tilde Re}(X)=X^p$, the sequence $(\phi_{M_{\tilde Re}}^s(t))_{s\in\dbN^*}$ converges in the topology of the $\tilde Re[{1\over p}]$-module $\End(M_{\tilde Re})[{1\over p}]$ defined by the sequence $(I(m)\End(M_{\tilde Re})[{1\over p}])_{m\in\dbN^*}$ of ${\tilde Re}[{1\over p}]$-submodules to an element $t_0\in\Lie(\scrT_{\tilde Re[{1\over p}]})$ which is fixed by $\phi_{M_{\tilde Re}}$ and which lifts $t$. It is easy to see that this implies that  $\Lie(\scrT_{\tilde Re[{1\over p}]})$ is $\tilde Re[{1\over p}]$-generated by elements fixed by $\phi_{M_{\tilde Re}}$. This proves (a). As $t_0$ modulo $\Ker(\chi_e)[{1\over p}]$ leaves invariant $F^1_K$, from Grothendieck--Messing deformation theory we get that an integral $p$-power of $t_0$ is the crystalline realization of an endomorphism of $D_V$. From this and (a) we get that (b) holds.\endproof

\bigskip\noindent
{\bf 3.6. Connection to the Main Problem.} We assume $p\ge 3$. Let $D_V$ be a $p$-divisible group over $k$ that lifts $D$. Let $(M_{\tilde Re},F^1_V,\phi_{M_{\tilde Re}})$ be associated to $D_V$ as in the proof of Theorem 3.4. Under the identification (4), we can naturally view  $\scrG_{\tilde Re[{1\over p}]}$ as the subgroup scheme of $\pmb{\GL}_{M_{\tilde Re}[{1\over p}]}$ that fixes $t_{\alpha}$ for all $\alpha\in\scrJ$. Let $\scrG^{\prime}_{\tilde Re}$ be the schematic closure of $\scrG_{\tilde Re[{1\over p}]}$ in $\pmb{\GL}_{M_{\tilde Re}}$ (in general it is neither a closed subgroup scheme of $\pmb{\GL}_{M_{\tilde Re}}$ nor isomorphic to $\scrG_{\tilde Re}$). We have the following consequence of Theorem 3.5.

\medskip\noindent
{\bf 3.6.1. Corollary.} {\it We assume that $p\ge 3$ and that the schematic closure $\tilde \scrG_V^\prime$ of $\chi_e^*(\scrG^{\prime}_{\tilde Re})_{K}$ in $\pmb{\GL}_{M_{\tilde Re}\otimes_{\tilde Re} {}_{\chi_e}V}$ is a reductive group scheme over $V$ whose special fibre, under the canonical identification $M_{\tilde Re}\otimes_{\tilde Re} k=M/pM$, is $\scrG_k$. We also assume that there exists a cocharacter $\tilde\mu_V:\dbG_m\to\tilde \scrG_V^\prime$ that acts on $F^1_V$ via weight $-1$ and that fixes $[M_{\tilde Re}/\Ker(\chi_e)M_{\tilde Re}]/F^1_V=H^1_{\text{dR}}(D_V/V)/F^1_V$. Then $\scrG^{\prime}_{\tilde Re}$ is a reductive, closed subgroup scheme of $\pmb{\GL}_{M_{\tilde Re}}$ isomorphic to $\scrG_{\tilde Re}$ and $(M_{\tilde Re},F^1_V,\phi_{M_{\tilde Re}},\scrG^{\prime}_{\tilde Re})$ is a ramified lift of $\scrC$ to $V$.}

\medskip
\proof
For  $s\in\{1,\ldots,e\}$ let $R_{1,s}:=\tilde Re/I(s)=W(k)[[X]]/(X^s)$. By induction on $s\in\{1,\ldots,e\}$ we show that the schematic closure $\tilde \scrG^\prime_{R_{1,s}}$ of $\scrG_{B(k)[[X]]/(X^s)}$ in $\pmb{\GL}_{M_{\tilde Re}/I(s)M_{\tilde Re}}$ (i.e., in $\pmb{\GL}_{M_{\tilde Re}}\times_{\tilde Re} \Spec R_{1,s}$) is a reductive, closed subgroup  scheme. The case $s\le p-1$ is obvious as the ideal $(X)/(X^s)$ of $R_{1,s}$ has a nilpotent divided power structure. More precisely, for $s\le p-1$ the reduction modulo $I(s)[{1\over p}]$ of the identification (4) induces a canonical identification
$(M\otimes_{W(k)} R_{1,s},\phi\otimes\Phi_{R_1,s})=(M_{\tilde Re},\phi_{M_{\tilde Re}})\otimes_{\tilde Re} R_{1,s}$,
where $\Phi_{1,s}$ is the Frobenius lift of $R_{1,s}=R_1/(X^s)$ which is compatible with $\sigma$ and which annihilates $X$ modulo $(X^s)$. 
If $p-1\Le s\Le e-1$, then the passage from $s$ to $s+1$ goes as follows. 

Let $\grm_V$ be the maximal ideal of $V$. Under the canonical identification $M_{\tilde Re}\otimes_{\tilde Re} k=M/pM$, we can also identify $F^1_V/\grm_VF^1_V=F^1/pF^1$ and $\tilde \scrG_V^\prime=\scrG_k$; therefore we can view both $\tilde\mu_V$ modulo $\grm_V$ and $\mu_k$ as cocharacters of the parabolic subgroup $\scrP_k$ of $\scrG_k$. By replacing $\mu$ with a $\scrP(W(k))$-conjugate of it, we can assume that $\tilde\mu_V$ modulo $\grm_V$ commutes with $\mu_k$. As $\tilde\mu_V$ modulo $\grm_V$ and $\mu_k$ are two commuting cocharacters of $\scrP_k$ that act in the same way on $F^1/pF^1$ and $(M/pM)/(F^1/pF^1)$, they coincide. Let $\tilde\mu_{R_{1,s}}:\dbG_m\to\tilde \scrG^\prime_{R_{1,s}}$ be a cocharacter that lifts both $\tilde\mu_{V}$ modulo $\grm_V$ and $\mu$, cf. [DG, Vol. II, Exp. IX, Thms. 3.6 and 7.1]. Let $F^1_{R_{1,s}}$ be the direct summand of $M_{\tilde Re}/I(s)M_{\tilde Re}$ which lifts $F^1_V/\grm_VF^1_V=F^1/pF^1$ and which is normalized by $\tilde\mu_{R_{1,s}}$. 

Let $m:=\dim(\scrG_{B(k)}^{\der})$. Similarly to either [Fa, proof of Thm. 10 and Rm. (iii) of p. 136] or [Va8, proof of Thm. 5.3] we get that the quadruple $(M_{\tilde Re}/I(s)M_{\tilde Re},F^1_{R_{1,s}},\phi_{M_{\tilde Re}},(t_{\alpha})_{\alpha\in\scrJ})$ is induced from the quadruple $\scrC^{R_m}$ of the proof of Theorem 3.5 via a $W(k)$-homomorphism $j_s:R_m\to R_{1,s}$ that maps the ideal $(X_1,\ldots,X_m)$ to the ideal $(X)/(X^s)$ (the arguments for this are the same as the ones of the last paragraph of Theorem 3.5). Here we denote also by $\phi_{M_{\tilde Re}}$ its reduction modulo $I(s)$. As the ideal $(X^s)/(X^{s+1})$ of $R_{1,s+1}$ has naturally a trivial divided power structure and as $j_s$ lifts to a $W(k)$-homomorphism $j_{s+1}:R_m\to R_{1,s+1}$ that maps the ideal $(X_1,\ldots,X_m)$ to the ideal $(X^s)/(X^{s+1})$, the triple $(M_{\tilde Re}/I(s+1)M_{\tilde Re},\phi_{M_{\tilde Re}},(t_{\alpha})_{\alpha\in\scrJ})$ is the extension of $(M_{R_m},\phi_{M_{R_m}},(t_{\alpha}^{R_m})_{\alpha\in\scrJ})$ via such a homomorphism $j_{s+1}$. Thus $\tilde \scrG^\prime_{R_{1,s+1}}$ is the pull-back of $\scrG_{R_m}$ via the morphism $j_{s+1}:\Spec R_{1,s+1}\to\Spec R_m$ defined by $j_{s+1}$ (and denoted in the same way) and it is therefore a reductive, closed subgroup scheme of $\pmb{\GL}_{M_{\tilde Re}/I(s+1)M_{\tilde Re}}$. This ends the induction.

 Let $\tilde\mu_{R_{1,e}}:\dbG_m\to\tilde \scrG^\prime_{R_{1,e}}$ be a cocharacter that lifts both $\mu$ and the reduction modulo $p$ of $\mu_V$, cf. [DG, Vol. II, Exp. IX, Thms. 3.6 and 7.1] and the smoothness of $\scrG^\prime_{R_{1,e}}$. Let $F^1_{R_{1,e}}$ be the direct summand of $M_{\tilde Re}/I(e)M_{\tilde Re}$ which lifts $F^1_V/\grm_VF^1_V=F^1/pF^1$ and which is normalized by $\tilde\mu_{R_{1,e}}$. Similarly to either [Fa, proof of Thm. 10 and Rm. (iii) of p. 136] or [Va8, proof of Thm. 5.3] we get that $(M_{\tilde Re}/I(e)M_{\tilde Re},F^1_{R_{1,e}},\phi_{M_{\tilde Re}},(t_{\alpha})_{\alpha\in\scrJ})$
 is induced from $\scrC^{R_m}$ via a $W(k)$-homomorphism $j_e:R_m\to R_{1,e}$ that maps the ideal $(X_1,\ldots,X_m)$ to the ideal $(X)/(X^e)$. Here we denote also by $\phi_{M_{\tilde Re}}$ its reduction modulo $I(e)$. 

Let $D_{\text{univ}}$ be the $p$-divisible group over $R_m/pR_m$ whose $F$-crystal is $(M\otimes_{W(k)} R_m,g^{\der}_{\text{univ}}(\phi\otimes\Phi_{R_m}),\nabla_{\text{univ}})$ (cf. end of Subsection 3.3). A second induction on $s\in\{1,\ldots,e\}$ shows (based on Grothendieck--Messing deformation theory) that there exists a $W(k)$-homomorphism $j_s:R_m\to R_{1,s}$ that maps $(X_1,\ldots,X_m)$ to $(X)/(X^s)$ and such that the pull-back of $D_{\text{univ}}$ via the morphism $\Spec k[[X]]/(X^s)\to\Spec R_m$ defined by $j_s$ is $D_V\times_{V} \Spec k[[X]]/(X^e)$. Taking $s=e$ we get that we can assume that $D_V\times_{V} \Spec V/pV$ is the extension via $j_e$ modulo $p$ of $D_{\text{univ}}$. Thus $(M_{\tilde Re},\phi_{M_{\tilde Re}},(t_{\alpha})_{\alpha\in\scrJ})$ is the extension of $(M\otimes_{W(k)} R_m,g^{\der}_{\text{univ}}(\phi\otimes\Phi_{R_m}),(t_{\alpha})_{\alpha\in\scrJ})$ via a (any) $W(k)$-homomorphism $R_m\to \tilde Re$ which lifts $j_e$ modulo $p$ (the fact that under such an extension and the identification (4), each $t_{\alpha}$ is mapped to $t_{\alpha}$ follows from the fact that no element of $\scrT(M[{1\over p}])\otimes_{B(k)} I(e)[{1\over p}]$ is fixed by $\phi\otimes\Phi_{\tilde Re}$). Thus the closed embedding $\scrG^{\prime}_{\tilde Re}\hookrightarrow \pmb{\GL}_{M_{\tilde Re}}$ is the pull-back of the closed embedding $\scrG_{R_m}\hookrightarrow  \pmb{\GL}_{M\otimes_{W(k)} R_m}$ via a morphism $\Spec \tilde Re\to\Spec R_m$ which lifts $j_e$ modulo $p$. As $\scrG_{R_m}$ is a reductive, closed subgroup scheme of $\pmb{\GL}_{M\otimes_{W(k)} R_m}$, we conclude that $\scrG^{\prime}_{\tilde Re}$ is also a reductive, closed subgroup scheme of $\pmb{\GL}_{M_{\tilde Re}}$. 

As $\scrG^{\prime}_{\tilde Re}$ is smooth over $\tilde Re$ and due to the property (i) of the proof of Theorem 3.5, there exists a cocharacter $\mu_{\tilde Re}$ of $\scrG^{\prime}_{\tilde Re}$ that lifts both $\mu$ and $\tilde\mu_V$ (to be compared with [Va1, Lem. 5.3.2] or [Ki2, Prop. (1.1.5) and Lem. (1.4.5)]). Let $F^1_{\tilde Re}$ be the direct summand of $M_{\tilde Re}$ that lifts $F^1_V$ and that is normalized by $\mu_{\tilde Re}$. The quadruple $(M_{\tilde Re},F^1_{\tilde Re},\phi_{M_{\tilde Re}},\scrG^\prime_{\tilde Re})$ is a lift of $\scrC$ to $\tilde Re$ (the analogue for $\tilde Re$ of the axiom (iii) of Definition 3.2.1 holds for this quadruple due to the very definition of $\scrG^\prime_{\tilde Re}$). Thus $(M_{\tilde Re},F^1_V,\phi_{M_{\tilde Re}},\scrG^\prime_{\tilde Re})$ is a ramified lift of $\scrC$ to $V$. \endproof

\medskip\noindent
{\bf 3.6.2. Definition.} We assume $p\ge 3$. Let $D_V$ be a $p$-divisible group over $V$ that lifts $D$. We say that $D_V$ is a ramified lift of $D$ to $V$ with respect to $\scrG$ if (to be compared with the Manin Problem 1.2) the following three axioms hold: 

\medskip
{\bf (a)} under the canonical identification $H^1_{\text{dR}}(D_V/V)[{1\over p}]=M\otimes_{W(k)} V[{1\over p}]$ (see proof of Corollary 3.5.1), the schematic closure $\scrG^\prime_V$ of $\scrG_K$ in $\pmb{\GL}_{H^1_{\text{dR}}(D_V/V)}$ is a reductive, closed subgroup scheme of $\pmb{\GL}_{H^1_{\text{dR}}(D_V/V)}$; 

\smallskip
{\bf (b)} under the canonical identification $M/pM=H^1_{\text{dR}}(D_V/V)/\grm_VH^1_{\text{dR}}(A/V)$, the group scheme $\scrG^\prime_V$ lifts $\scrG_k$ (we recall that $\grm_V$ is the maximal ideal of $V$);

\smallskip
{\bf (c)} there exists a cocharacter $\dbG_m\to \scrG^\prime_V$ that acts on $F^1_V$ via the inverse of the identical character of $\dbG_m$ and that fixes $H^1_{\text{dR}}(D_V/V)/F^1_V$, where $F^1_V$ is the direct summand of $H^1_{\text{dR}}(D_V/V)$ which is the Hodge filtration of $D_V$.
\medskip

\medskip\noindent
{\bf 3.6.3. Corollary.} {\it We assume that $p\Ge 3$. Then the ramified lifts  of $\scrC$ to $V$ are in natural bijection to the ramified lifts of $D$ to $V$ with respect to $\scrG$.}

\medskip
\proof
Let $(M_{\tilde Re},F^1_{\tilde Re},\phi_{M_{\tilde Re}},\scrG_{\tilde Re})$ be a ramified lift of $\scrC$ to $\tilde Re$. Let $F^1_V:=F^1_{\tilde Re}\otimes_{\tilde Re} {}_{\chi_e}V$. Let $D_V$ be the $p$-divisible group over $V$ that corresponds to $(M_{\tilde Re},F^1_V,\phi_{M_{\tilde Re}})$ via Theorem 3.4. From Definitions 3.2.1 (a) and (d) and 3.6.2 we get that $D_V$ is a ramified lift of $D$ to $V$ with respect to $\scrG$.  The rule $(M_{\tilde Re},F^1_V,\phi_{M_{\tilde Re}},\scrG_{\tilde Re})\mapsto D_V$ defines a map which is injective (cf. Theorem 3.4) and onto (cf. Corollary 3.6.1) and thus it is the searched for bijection.\endproof 
 
\bigskip\noindent
{\bf 3.7. Remarks.} {\bf (a)} Sections 2.8, 3.4, and 3.6 hold with $k$ replaced by an arbitrary perfect field of characteristic $p$. Theorem 3.4 was first obtained in [Br] and [Zi2] (strictly speaking these references worked with $Re$ instead of $\Rtil e$ but as $\Phi_{\tilde Re}(\tilde Re)\subseteq Re$, it is easy to see that for $p\ge 3$ there exists a natural bijection between lifts of $\scrC$ to $Re$ with respect to $V$ and lifts of $\scrC$ to $\Rtil e$ with respect to $V$). 

\smallskip
{\bf (b)} In this paper the principally quasi-polarized contexts are treated as by products (consequences) of non--polarized contexts. This is so due to the following two reasons. First, often the principally quasi-polarized context is handled by making small (if any at all) modifications to the contexts that involve only $\scrC$. There exists no element of $\scrT(M)\otimes_{W(k)} I(1)$ fixed by $\phi\otimes \Phi_{\tilde Re}$. Thus if we have a principal quasi-polarization $\lambda_M:M\times M\to W(k)$ of $\scrC$ and if in Theorem 3.5 we have a lift $(M_{\tilde Re},F^1_{\tilde Re},\phi_{M_{\tilde Re}},\scrG_{\tilde Re},\lambda_{M_{\tilde Re}})$ of $(\scrC,\lambda_M)$ to $\tilde Re$, then the $W(k)$-homomorphism $R_m\to \tilde Re$ of Theorem 3.5 takes automatically $\lambda_M$ to $\lambda_{M_{\tilde Re}}$. 

To explain the second reason we assume that $\scrG$ is generated by $Z(\pmb{\GL}_M)$ and a reductive, closed subgroup scheme $\scrG^0$ of $\pmb{\text{SL}}_M$ and that the intersection $Z(\pmb{\GL}_M)\cap \scrG^0$ is either $\pmb{\mu}_2$ or $\Spec W(k)$. Then for most applications we can replace $\scrC$ by the direct sum $\scrC\oplus\scrC^{\vee}(1):=(M\oplus M^{\vee}(1),\phi\oplus p1_{M^{\vee}}\circ\phi,\scrG)$, where $M^{\vee}(1):=M^{\vee}$ as $W(k)$-modules and where $\scrG$ is a reductive, closed subgroup scheme of $\pmb{\GL}_{M\oplus M^{\vee}(1)}$ such that $\scrG^0$ acts on $M^{\vee}(1)$ via its action on $M^{\vee}$ and $Z(\pmb{\GL}_M)$ is naturally identified with $Z(\pmb{\GL}_{M\oplus M^{\vee}(1)})$. Let $\lambda_{M\oplus M^{\vee}(1)}$ be the natural principal alternating quasi-polarization of $\scrC\oplus\scrC^{\vee}(1)$. Then each ramified lift of $\scrC$ to $V$ induces naturally a unique ramified lift of $(\scrC\oplus\scrC^{\vee}(1),\lambda_{M\oplus M^{\vee}(1)})$ to $V$.

\bigskip\smallskip
\noindent
{\boldsectionfont 4. The basic results}

\bigskip
In this section we state our basic results pertaining to the (ramified) $\text{CM}$-classifications of Section 3 (see Basic Theorems 4.1 and 4.2). Let $(M_{\dbZ_p},\scrG_{\dbZ_p},(t_{\alpha})_{\alpha\in\scrJ})$ be as in Subsubsection 2.5.1. For properties such as $Q\grU$, $\grR$, $T\grA$, etc., we refer to Definitions 2.4.

\bigskip\noindent
{\bf 4.1. Basic Theorem.} {\it We assume that $\scrC$ is semisimple and basic. We have: 

\medskip
{\bf (a)} If $[\scrC]\in Z(\scrY(\scrF))$, then $Q\grU$ holds for $\scrC$.

\smallskip
{\bf (b)} Condition $\grR$ (resp. $Q\grR$, $T\grR$, or $TT\grR$) holds for $\scrC$ if and only if $\grA$ (resp. $Q\grA$, $T\grA$, or $TT\grA$) holds for $\scrC$.

\smallskip
{\bf (c)} We assume that $p\Ge 3$ and that $Q\grA$ holds for $\scrC$. Then we have $[\scrC]\in Z^{\text{ram}}(\scrY(\scrF))$.}

\medskip
Directly from Theorem 4.1 (b) and the Definition 2.4 (g) we get:

\medskip\noindent
{\bf 4.1.1. Corollary.} {\it We assume that $\scrC$ is semisimple. Then the condition $\grR$ (resp. $Q\grR$, $T\grR$, or $TT\grR$) holds for $\scrC$ if and only if $\grA$ (resp. $Q\grA$, $T\grA$, or $TT\grA$) holds for $\scrC$.}

\bigskip\noindent
{\bf 4.2. Basic Theorem.} {\it {\bf (a)} We assume that $\scrC$ is basic and semisimple. We also assume that each simple factor of $\scrG^{\ad}_{W(\bar k)}$ is of $B_n$, $C_n$, or $D_n$ Lie type. Then $TT\grR$ holds for $\scrC$.

\smallskip
{\bf (b)} We assume that $\scrC$ is semisimple and of $B_n$ and $D_n$ type in the sense of Definition 1.2.4 (a). Then $TT\grR$ holds for $\scrC$.} 

\medskip
In Sections 5 and 6 we prove Basic Theorems 4.1 and 4.2 (respectively). The following corollary is an abstract extension of [Zi1, Thm. 4.4] for $p\ge 3$ and for contexts related to Shimura pairs of either $B_n$ or $D_n^{\dbR}$ type. It also represents the very first situations where complete ramified $\text{CM}$-classifications are accomplished. 

\bigskip\noindent
{\bf 4.3. Corollary.} {\it We assume that $p\ge 3$ and that $\scrC$ is of $B_n$ and $D_n^{\dbR}$ type. Then for $g\in \scrG(W(k))$ we have $[\scrC_g]\in Z^{\text{ram}}(\scrY(\scrF))$ if and only if $\scrC_g$ is semisimple.}

\medskip
\proof
If $[\scrC_g]\in Z^{\text{ram}}(\scrY(\scrF))$, then by performing $\grO_1$ we can assume that there exists a maximal torus $\scrT_{1,B(k)}$ of $\scrG_{B(k)}$ of $\dbQ_p$-endomorphisms of $\scrC_g$ (see Remark 3.2.2 (a) and Corollary 3.5.1 (a) and (b)). Thus $(g\phi)^r\in\scrG(B(k))$ centralizes $\Lie(\scrT_{1,B(k)})$ and thus also $\scrT_{1,B(k)}$. Therefore $(g\phi)^r\in \scrT_{1,B(k)}(B(k))$ is semisimple i.e., $\scrC_g$ is semisimple. 

Suppose $\scrC_9$ is semisimple. To prove that $[\scrC_g]\in Z^{\text{ram}}(\scrY(\scrF))$, we can assume that $g=1_M$ (in order not to complicate the notations), that $L_{\scrG}^0(\phi)$ is a reductive group scheme (cf. Subsection 2.6), and that there exists a maximal torus $\scrT_{1,B(k)}$ of $\dbQ_p$-endomorphisms of $(M,\phi,L_{\scrG}^0(\phi))$ (cf. Fact 3.1 (c)). Let $(\scrT_{1,B(k)},\mu_1)$  be an $E$-pair of $(M,\phi,L_{\scrG}^0(\phi))$ that satisfies the $\grC$ condition (cf. Theorem 4.2 (b)); it is admissible (cf. Theorem 4.1 (b)). Thus $[\scrC]\in Z^{\text{ram}}(\scrY(\scrF))$, cf. Theorem 4.1 (c).\endproof 

\bigskip\smallskip
\noindent
{\boldsectionfont 5. Proof of Basic Theorem 4.1}

\bigskip
In this section we assume that $\scrC$ is basic and semisimple. In Subsections 5.1, 5.2, and 5.3 we prove Theorems 4.1 (a), 4.1 (b), and  4.1 (c) (respectively). 

\bigskip\noindent
{\bf 5.1. Proof of 4.1 (a).} To prove Theorem 4.1 (a) we can assume that there exists a maximal split torus $\scrT$ of $\scrG$ such that we have $\phi(\Lie(\scrT))=\Lie(\scrT)$. By performing the operation $\grO_1$, we can assume that $\phi^r\in\scrT(B(k))$ fixes $\Lie(\scrT)$ (see proof of Fact 3.1 (a)). Thus $\Lie(\scrT)$ is generated by elements fixed by $\phi$, cf. property 2.5 (a). Let $C_\scrT:=C_{\pmb{\GL}_M}(\scrT)$; it is a reductive, closed subgroup scheme of $\pmb{\GL}_M$ such that we have $\phi(\Lie(C_\scrT))=\Lie(C_\scrT)$. This implies that $F^1/pF^1$ is a $\Lie(C_{\scrT,k})$-module. As $C_\scrT(W(\bar k))$ is naturally a subset of $\Lie(C_{\scrT,W(\bar k)})$, the group $C_{\scrT,k}$ normalizes $F^1/pF^1$. Thus $\scrT_k$ is a maximal torus of $\scrP_k$. Let $\scrT_0$ be a maximal torus of $\scrP$ through which the cocharacter $\mu:\dbG_m\to \scrG$ factors. From [DG, Vol. II, Exp. IX, Thms. 3.6 and 7.1] and [Bo, Ch. V, Thm. 15.14] we get the existence of an element $g\in \scrG(W(k))$ which modulo $p$ belongs to $\scrP(k)$ and for which we have $g(\scrT_0)g^{-1}=\scrT$. By replacing $\mu$ with its inner conjugate through $g$ we can assume that $\mu$ factors through $\scrT$. As $\scrC$ is basic, the product of cocharacters of $\scrT$ of the orbit of $\mu$ under integral powers of $\phi$ factors through $Z^0(\scrG)$. Thus $(\scrT_{B(k)},\mu_{B(k)})$ is an unramified $E$-pair of $\scrC$ that satisfies the cyclic $\grC$ condition, cf. Example 2.4.1. This proves Theorem 4.1 (a).\endproof

\bigskip\noindent
{\bf 5.2. Proof of 4.1 (b).} To prove Theorem 4.1 (b), it is enough to show that an $E$-pair $(\scrT_{1,B(k)},\mu_1)$ of the basic $\scrC$ is admissible if and only if $\grU$ holds for it. If an $E$-pair $(\scrT_{1,B(k)},\mu_1)$ is admissible, then (as $\scrC$ is basic) from [Ko1, Subsects. 2.4 and 2.5] and [RaZ, Prop. 1.21] (see also [RR, Thm. 1.15]) we get that the product of the cocharacters of $\scrT_{1,K_1}$ which belong to the $\Gal(K_1/\dbQ_p)$-orbit of $\mu_1$ factors through $Z^0(\scrG_{K_1})$. Thus the $E$-pair $(\scrT_{1,B(k)},\mu_1)$ of $\scrC$ satisfies the $\grC$ condition, cf. Example 2.4.1. Thus $\grR$ holds for $\scrC$.

We now show that the converse holds i.e., we prove that if
$$(\scrT_{1,B(k)},\mu_1,\tau=(\tau_1,\tau_2,\ldots,\tau_l))$$ 
is an $E$-triple of $\scrC$ such that the condition 2.4 (e1) holds and if $K_2$ and $F^1_{K_2}$ are as in Definitions 2.4 (d) and (h), then the filtered module $(M[{1\over p}],\phi,F^1_{K_2})$ over $K_2$ is admissible. It is enough to show that the filtered module $(M[{1\over p}],\phi,F^1_{K_2})$ over $K_2$ is weakly admissible, cf. [CF, Thm. A]. Let $\bar K_2$ be the subfield of $\overline{B(\bar k)}$ generated by $B(\bar k)$ and $K_2$. For $i\in\{1,\ldots,l\}$ let $M_i:=M$. Let 
$$O:=(\oplus_{i=1}^l M_i)\otimes_{W(k)} K_2.$$ 
Let $\bar O:=O\otimes_{K_2} \bar K_2$. Let $d\in\dbN^*$ be as in the condition 2.4 (e1). By replacing $(\tau_1,\tau_2,\ldots,\tau_l)$ with $(\tau_1,\ldots,\tau_{dl})$ (cf. notations of Definition 2.4 (d)), we can assume that $d=1$. 

We denote also by $\tau$ the $\sigma$-linear automorphism
$$\tau:O\arrowsim O$$
which takes $m_i\otimes v_2\in M_i\otimes_{W(k)} K_2$ to $\phi(m_i)\otimes\tau_i(v_2)\in M_{i+1}\otimes_{W(k)} K_2$, where $M_{l+1}:=M_1$. We view $\scrT_{1,K_2}^l:=\prod_{i=1}^l \scrT_{1,K_2}$ (resp. $\scrG^l_{K_2}:=\prod_{i=1}^l \scrG_{K_2}$) as a torus (resp. as a reductive subgroup) of $\pmb{\GL}_O$. We embed $\scrT_{1,K_2}$ (resp. $\scrG_{K_2}$) diagonally into $\scrT_{1,K_2}^l$ (resp. $\scrG^l_{K_2}$). Let $\mu_2$ be the cocharacter of $\scrT_{1,K_2}^l$ which normalizes each $M_i\otimes_{W(k)} K_2$ and which acts on $M_i\otimes_{W(k)} K_2$ identified with $M\otimes_{W(k)} K_2$ as $\mu_{1,K_2}$ does. We consider the $\sigma$-linear automorphism
$$\sigma_2:=\tau\mu_2(p):O\arrowsim O.$$
We denote also by $\tau$ and $\sigma_2$ their $\sigma_{\bar k}$-linear extensions to $\bar O$. 

\medskip\noindent
{\bf 5.2.1. Lemma.} {\it All Newton polygon slopes of $(O[{1\over p}],\sigma_2)$ are $0$ (equivalently, $\bar O$ is $\bar K_2$-generated by elements fixed by $\sigma_2$)}.

\medskip
\proof
The actions of $\tau$ and $\sigma_2$ on cocharacters of $\scrT_{1,K_2}^l$ are the same. As $\scrC$ is basic, the Newton quasi-cocharacter $\nu$ of $\scrC$ factors through $Z^0(\scrG_{B(k)})$ (see [Va7, Cor. 2.3.2] or proof of Fact 2.6.1). We consider the quasi-cocharacter $\nu_2$ of $\scrT_{1,K_2}^l$ which is the mean average of the orbit of $\mu_2$ under integral powers of $\tau$. It factors through $Z^0(\scrG_{K_2}^l)$, cf. property 2.4 (e1). Strictly speaking we get directly this only for the first factor $\scrG_{K_2}$ of $\scrG_{K_2}^l$. However, due to the circular aspect of $\tau$ and of the condition 2.4 (e1), this extends automatically to all the other $l-1$ factors $\scrG_{K_2}$ of $\scrG_{K_2}^l$. As each $\tau_i$ extends the Frobenius automorphism $\grF_{2u}$ of $K_{2u}$, the composite of $\nu_2$ with $\scrG^l_{K_2}\to (\scrG^l_{K_2})^{\ab}=\prod_{i=1}^l \scrG_{K_2}^{\ab}$ is such that its natural projections on the factors $\scrG_{K_2}^{\ab}$ are all the same and equal to the composite of $\nu$ with the natural epimorphism $\scrG_{K_2}\twoheadrightarrow \scrG_{K_2}^{\ab}$. As $\nu_2$ factors through $Z^0(\scrG_{K_2}^l)$, from the last sentence we get that in fact $\nu_2$ factors through the subtorus $Z^0(\scrG)_{K_2}$ of $\scrT_{1,K_2}$ and this factorization coincides with the factorization of $\nu$ through the subtorus $Z^0(\scrG)_{K_2}$ of $\scrT_{1,K_2}$. From this and the fact that for each $a\in\dbN^*$ we have $\sigma_2^{la}=[(\prod_{m=1}^{la}\tau^m(\mu_2))(p)]\tau^{la}$, we easily get that all Newton polygon slopes of $(O[{1\over p}],\sigma_2)$ are $0$ ($\tau^{la}$ acts on $\oplus_{i=1}^l M_i=M^l$ as $\phi^{la}$ does).\endproof

\medskip\noindent
{\bf 5.2.2. Proposition.} {\it Let $t:=\dim_{K_2}(O)=l\text{rk}_{W(k)}(M)$. There exists a $\bar K_2$-basis $\varUpsilon=\{e_1,\ldots,e_t\}$ of $\bar O$ and a permutation $\pi$ of $\{1,\ldots,t\}$ such that for all $i\in\{1,\ldots,t\}$ we have $\tau(e_i)=p^{n_i}e_{\pi(i)}$, where $n_i\in\{0,1\}$ is $1$ if and only if we have $e_i\in (\oplus_{i=1}^l F^1_{K_2})\otimes_{K_2} \bar K_2$. Moreover, $(\oplus_{i=1}^l F^1_{K_2})\otimes_{K_2} \bar K_2$ has a $\bar K_2$-basis formed by those elements $e_i\in \varUpsilon$ with $n_i=1$.}
\medskip
\proof
For $i\in \{1,2,\ldots,l\}$ and $j\in\{0,1\}$ let $^iF^j$ be the
$\bar K_2$-vector subspace of $\bar O$ on which $\tau^i(\mu_2)$ acts trivially if $j=0$
and via the inverse of the identical character of $\dbG_m$ if $j=1$. For a function $\bar f:\{1,2,\ldots,l\}\to\{0,1\}$ let
$$F_{\bar f}:=\bigcap_{i\in\{1,2,\ldots,l\}}\,^iF^{\bar f(i)}.$$ 
Let $\scrQ$ be the set of such functions $\bar f$ with $F_{\bar f}\ne 0$. As $\tau^i(\mu_2)$'s commute (being cocharacters of $\scrT_{1,K_2}^l$) we have a direct sum decomposition $\bar O=\oplus_{\bar f\in\scrQ}F_{\bar f}$. 
Let 
$$
\bar\tau:\scrQ\to\scrQ
$$
be the bijection defined by the rule: $\bar\tau(\bar f)(i)=\bar f(i-1)$, where $\bar f(0):=\bar f(l)$. 
 
Let $I^0_{\bar\tau}:=\{\bar f\in\scrQ|\bar {\tau}(\bar f)=\bar f\}$. Let $\bar\tau=\prod_{j\in I_{\bar\tau}}\bar\tau_j$ be written as a product of disjoint cyclic permutations. We allow trivial cyclic permutations i.e., we have a disjoint union 
$$I_{\bar\tau}=I^1_{\bar\tau}\cup I^0_{\bar\tau}$$ 
with the property that $j\in I_{\bar\tau}$ belongs to $I^1_{\bar\tau}$ if and only if $\bar {\tau}_j$ is a non-trivial permutation. If $j\in I^1_{\bar\tau}$, then each function $\bar f\in\scrQ$ such that we have $\bar {\tau}_j(\bar f)\ne\bar f$ is said to be associated to $\bar {\tau}_j$. Also $\bar f\in I^0_{\bar\tau}$ is said to be associated to $\bar {\tau}_{\bar f}$. As $\bar\tau^{l}=1_\scrQ$, the order $d_j$ of the cyclic permutation $\bar\tau_j$ divides $l$. For each $j\in I_{\bar\tau}$ we choose arbitrarily an element $\bar f_j$ of $\scrQ$ which is associated to $\bar\tau_j$. We have $\tau^{d_j}(F_{\bar f_j})=\sigma_2^{d_j}(F_{\bar f_j})=F_{\bar f_j}$. Let $_pF_{\bar f_j}:=\{x\in F_{\bar f_j}|\sigma_2^{d_j}(x)=x\}$.

For each $j\in I_{\bar\tau}$ we choose a $\bar K_2$-basis $\{e_s|s\in\varUpsilon_j\}$ for $F_{\bar f_j}$ formed by elements of ${}_pF_{\bar f_j}$, cf. Lemma 5.2.1. For each cyclic permutation $\bar\tau_j$ of length $\ge 2$ (i.e., for when we deal with a $j\in I^1_{\bar\tau}$) and for every element $\bar f\in\scrQ$ associated to $\bar\tau_j$ but different from $\bar f_j$,
let $u(\bar f)\in\dbN^*$ be the smallest number such that $\bar f=\bar\tau^{u(\bar f)}(\bar f_j)$ and let $n_{u(\bar f),j}:=\sum^{u(\bar f)}_{i=1}\bar f_j(i)$. We get a $\bar K_2$-basis $\{{1\over p^{n_{u(\bar f),j}}}\tau^{u(\bar f)}(e_s)|s\in\varUpsilon_j\}$ for $F_{\bar f}$. The expressions of $n_{u(\bar f),j}$'s are a consequence of the following iteration formula $\tau^{u(\bar f)}=[(\prod_{m=1}^{u(\bar f)} \sigma_2^m(\mu_2))({1\over p})] \sigma_2^{u(\bar f)}$.

Let $\varUpsilon=\{e_1,\ldots,e_t\}$ be the $\bar K_2$-basis for $\bar O$ obtained by putting together the chosen $\bar K_2$-bases for $F_{\bar f}$'s with $\bar f\in\scrQ$. Let $\pi$ be the unique permutation of $\{1,\ldots,t\}$ such that for all $i\in\{1,\ldots,t\}$ we have $\tau(e_i)\in K_0e_{\pi(i)}$. From constructions we get that $\tau(e_i)=p^{n_i}e_{\pi(i)}$, where $n_i\in\{0,1\}$ has the desired property, and that $(\oplus_{i=1}^l F^1_{K_2})\otimes_{K_2} \bar K_2$ has a $\bar K_2$-basis formed by those elements $e_i\in \varUpsilon$ with $n_i=1$.\endproof

\medskip\noindent 
{\bf 5.2.3. End of the proof of 4.1 (b).} Let $\bar V_2$ be the ring of integers of $\bar K_2$. Let $\bar L$ be the $\bar V_2$-lattice of $\bar O$ generated by the elements of the $\bar K_2$-basis $\varUpsilon$ of Proposition 5.2.2. Let $\bar F$ be the direct summand of $\bar L$ generated by the elements of $\varUpsilon\cap ((\oplus_{i=1}^l F^1_{K_2})\otimes_{K_2} \bar K_2)$. We consider an arbitrary $B(k)$-submodule $\sharp$ of $M[{1\over p}]$ which is normalized by $\phi$. Let $\bar L_{\sharp}:=\bar L\cap (\oplus_{i=1}^l \sharp\otimes_{B(k)} \bar K_2)$; it is a $\bar V_2$-module which is a direct summand of $\bar L$ left invariant by $\tau$. As in Mazur theorem of [Ka, Thm. 1.4.1] we get that the Newton polygon of $(\bar L_{\sharp},\tau)$ (i.e., the Newton polygon of $l$ copies of $(\sharp,\phi)$) is below the Hodge polygon of $(\bar L_{\sharp},\tau,\bar L_{\sharp}\cap \bar F)$ (i.e., the Hodge polygon of $l$ copies of $(\sharp,\phi,(\sharp\otimes_{B(k)} K_2)\cap F^1_{K_2})$). This identification of Newton (resp. Hodge) polygons follows from the very definition of $O$ and $\tau$ (resp. from the very definition of $\bar F$). We conclude that the Newton polygon of $(\sharp,\phi)$ is below the Hodge polygon of $(\sharp,\phi,(\sharp\otimes_{B(k)} K_2)\cap F^1_{K_2})$. Therefore the filtered module $(M[{1\over p}],\phi,F^1_{K_2})$ over $K_2$ is indeed weakly admissible. Thus Theorem 4.1 (b) holds.\endproof

\bigskip\noindent
{\bf 5.3. Proof of 4.1 (c).} 
We begin the proof of Theorem 4.1 (c) with some \'etale considerations. To prove Theorem 4.1 (c) we can assume that there exists an $E$-pair $(\scrT_{1,B(k)},\mu_1)$ of $\scrC$ which is admissible. Let $\scrT_{1,\dbQ_p}$, $K_2$, and $F^1_{K_2}$ be as in Definitions 2.4 (b), (d), and (h). The torus $Z^0(\scrG_{\dbQ_p})$ is a subtorus of $\scrT_{1,\dbQ_p}$. The triple $(M[{1\over p}],\phi,F^1_{K_2})$ is an admissible filtered module. Let $M_{\dbQ_p}$, $\rho$, $\scrW$, $\gamma\in H^1(\dbQ_p,\scrG_{\dbQ_p})$, and $(v_{\alpha})_{\alpha\in\scrJ}$ be as in Subsubsection 2.5.4. As $p\Ge 3$, the Galois representation
$\rho:\Gal(K_2)\to \pmb{\GL}_\scrW(\dbQ_p)$ is associated to an isogeny class of $p$-divisible groups over the ring of integers $V_2$ of $K_2$ (cf. [Br, Cor. 5.3.3]). Let $\scrJ_1:=\Lie(\scrT_{1,\dbQ_p})$. For $\alpha\in\scrJ_1$ let $t_{\alpha}:=\alpha$. To prove Theorem 4.1 (c) we can assume that 
$$\scrJ\cap\scrJ_1=\Lie(Z^0(\scrG_{\dbQ_p}))\subseteq \Lie(\scrT_{1,\dbQ_p})$$ 
and that for each $\alpha\in\scrJ\cap\scrJ_1$ the two definitions of $t_{\alpha}$ define the same tensor of $\scrT(M[{1\over p}])$. Let $\scrJ_2:=\scrJ\cup\scrJ_1$. Let $\scrG^{\prime}_{\dbQ_p}$ be the subgroup of $\pmb{\GL}_\scrW$ that fixes $v_{\alpha}$ for all $\alpha\in\scrJ$; it is an inner form of $\scrG_{\dbQ_p}$. Until Subsubsection 5.3.2 below we will assume that there exists a $\dbQ_p$-linear isomorphism $\grJ:M_{\dbQ_p}\arrowsim \scrW$
that takes $t_{\alpha}$ to $v_{\alpha}$ for all $\alpha\in\scrJ$ (for instance, this holds if $\scrG^{\der}$ is simply connected, cf. Lemma 2.5.5). We use $\grJ$ to identify naturally $\scrG_{\dbQ_p}=\scrG^{\prime}_{\dbQ_p}$. Let $\scrG^{\prime}_{\dbZ_p}$ be the schematic closure of $\scrG^{\prime}_{\dbQ_p}$ in $\scrW_{\dbZ_p}:=\grJ(M_{\dbZ_p})$; we identify it with $\scrG_{\dbZ_p}$. Let 
$$\grJ_{\dbZ_p}:M_{\dbZ_p}\arrowsim \scrW_{\dbZ_p}$$ 
be the resulting $\dbZ_p$-linear isomorphism.

\medskip\noindent
{\bf 5.3.1. Crystalline machinery.} We now apply a standard crystalline machinery in order to show that $[\scrC]\in Z^{\text{ram}}(\scrY(\scrF))$. Let $K_3$ be a finite field extension of $K_2$ such that $\rho(\Gal(K_3))$ normalizes $\scrW_{\dbZ_p}$ and its ramification index $e$ is at least 2. Let $\tilde D_{K_3}$ be the $p$-divisible group over $K_3$ defined by restricting $\rho$ as follows
$$\varrho:\Gal(K_3)\to \pmb{\GL}_{\scrW_{\dbZ_p}}(\dbZ_p).$$
It extends to a $p$-divisible group $\tilde D_{V_3}$ over the ring of integers $V_3$ of $K_3$ (cf. the isogeny class part of the first paragraph of Subsection 5.3) and this extension is unique (cf. [Ta1]).

We use the notations of Subsection 2.8. To avoid extra notations, (by performing the operation $\grO_1$) we can assume that the residue field of $K_3$ is $k$. We fix a uniformizer $\pi_3$ of $V_3$. Let $\chi_{e,3}:\Rtil e\twoheadrightarrow V_3$ be a $W(k)$-epimorphism defined by $\pi_3$, cf. Fact 2.8.1. Let $(M_{\tilde Re},\phi_{M_{\tilde Re}},\nabla)$
be the Dieudonn\'e $F$-crystal over $\tilde Re/p\tilde Re$ of $\tilde D_{V_3}\times_{V_3} \Spec V_3/pV_3$ (viewed without the Verschiebung map), cf. proof of Theorem 3.4. 

Let $B^+(V_3)$ be the crystalline Fontaine ring of $V_3$ as defined in [Fa]. We recall that $B^+(V_3)$ is an integral, local $W(k)$-algebra which is endowed with a decreasing, exhaustive, and separated filtration $(F^i(B^+(W(k)))_{i\in\dbN}$, with a Frobenius lift $\grF$, and with a natural Galois action by $\Gal(K_3)$. Moreover we have a natural $W(k)$-epimorphism compatible with the natural Galois actions by $\Gal(K_3)$ of the form
$$s_{V_3}: B^+(V_3)\twoheadrightarrow \overline{V_3}^\wedge,$$
where $\overline{V_3}^\wedge$ is the $p$-adic completion of the normalization $\overline{V_3}$ of $V_3$ in $\overline{B(k)}$. We refer to loc. cit. for the natural $W(k)$-monomorphism $\tilde Re\hookrightarrow B^+(V_3)$ which respects the Frobenius lifts (and which is associated to the uniformizer $\pi_3$); its composite with $s_{V_3}$ is the composite of $\chi_{e,3}$ with the inclusion $V_3\hookrightarrow \overline{V_3}^\wedge$. We apply Fontaine comparison theory  to $\tilde D_{V_3}$ (see loc. cit. and [Va1, Subsect. 5.2]). We get a $B^+(V_3)$-monomorphism
$$i_{\tilde D_{V_3}}:M_{\tilde Re}\otimes_{\tilde Re} B^+(V_3)\hookrightarrow \scrW_{\dbZ_p}\otimes_{\dbZ_p} B^+(V_3)$$
 which has the following two properties:

\medskip
{\bf (a)} It respects the tensor product filtrations (the filtration of $\scrW_{\dbZ_p}$ is defined by: $F^1(\scrW_{\dbZ_p})=0$ and $F^0(\scrW_{\dbZ_p})=\scrW_{\dbZ_p})$.

\smallskip
{\bf (b)} It respects the Galois actions (the Galois action on $M_{\tilde Re}\otimes_{\tilde Re} B^+(V_3)$ is defined naturally via $s_{V_3}$ and the fact that $\Ker(s_{V_3})$ has a natural divided power structure). 
\medskip

For $\alpha\in\scrJ_2$ let $u_{\alpha}\in\scrT(M_{\tilde Re}[{1\over p}])$ be the tensor that corresponds to $\grJ(t_{\alpha})$ via $i_{\tilde D_{V_3}}$, cf. Fontaine comparison theory. From [Ki2, Cor. (1.4.3)] (to be compared with [Va8, Subsect. 5.5]) we get the existence of an isomorphism
$$\grJ_{\tilde Re}:M_{\tilde Re}\arrowsim \scrW_{\dbZ_p}\otimes_{\dbZ_p} \tilde Re$$
which takes $u_{\alpha}$ to $\grJ(t_{\alpha})$ for all $\alpha\in\scrJ$. Thus the schematic closure $\tilde \scrG_{\tilde Re}$ in $\pmb{\GL}_{M_{\tilde Re}}$ of the closed subgroup scheme of $\pmb{\GL}_{M_{\tilde Re}[{1\over p}]}$ that fixes $u_{\alpha}$ for all $\alpha\in\scrJ$, is a reductive subgroup scheme isomorphic to $\scrG_{\tilde Re}$. 

The tensorization of $(M_{\tilde Re}[{1\over p}],\phi_{M_{\tilde Re}})$ with the natural epimorphism $\tilde Re[{1\over p}]\twoheadrightarrow B(k)$ that takes $X$ to $0$, is the $F$-isocrystal of $D$ (cf. the very definition of $\tilde D_{V_3}$) and thus it is canonically isomorphic to $(M[{1\over p}],\phi)$. Under the resulting identification $(M_{\tilde Re}[{1\over p}],\phi_{M_{\tilde Re}})\otimes_{\tilde Re} B(k)=(M[{1\over p}],\phi)$, $u_{\alpha}$ is identified with $t_{\alpha}$ for all $\alpha\in\scrJ_2$. This implies that under the identification (4), $u_{\alpha}$ gets identified with $t_{\alpha}$ for all $\alpha\in\scrJ_2$. In particular, there exists a maximal torus of $\tilde \scrG_{\tilde Re[{1\over p}]}$ whose Lie algebra is $\tilde Re[{1\over p}]$-generated by those $u_{\alpha}$ with $\alpha\in\scrJ_1$ (i.e., which corresponds to the maximal torus $\scrT_{1,B(k)}$ of $\scrG_{B(k)}$ via the identification (4)).  

Due to the existence of the cocharacter $\mu_{1,K_2}:\dbG_m\to \scrT_{1,K_2}$ that acts on $F^1_{K_2}$ via the inverse of the identical character of $\dbG_m$, as in [Va1, Subsubsect. 5.3.1 and Lem. 5.3.2] (to be compared also with [Ki, Prop. (1.1.5) and Lem. (1.4.5)]) we argue that there exists a cocharacter $\tilde\mu_{\tilde Re}:\dbG_m\to\tilde \scrG_{\tilde Re}$ such that the following two properties hold:

\medskip
{\bf (c)}  there exists a direct sum decomposition $M_{\tilde Re}=F_{\tilde Re}^1\oplus F_{\tilde Re}^0$ such that $F_{\tilde Re}^1$ lifts the Hodge filtration $F^1_{V_3}$ of $M_{\tilde Re}\otimes_{\tilde Re} {}_{\chi_{e,3}} V_3$ defined by $\tilde D_{V_3}$ and for each $i\in\{0,1\}$, every element $\beta\in\dbG_m(\tilde Re)$ acts  on $F^i_{\tilde Re}$ through $\tilde\mu_{\tilde Re}$ as the multiplication with $\beta^{-i}$ (one can identify naturally $F^1_{V_3}[{1\over p}]=F^1_{K_2}\otimes_{K_2} K_3$);

\smallskip
{\bf (d)}  $\mu_{1,K_3}$ and the pull-back of $\tilde\mu_{\tilde Re}$ to a cocharacter of $\tilde \scrG_{K_3}=\scrG_{K_3}$ are $\scrG_{K_3}(K_3)$-conjugate (the identification $\tilde \scrG_{K_3}=\scrG_{K_3}$ used here is the one defined naturally by the tensorization of (4) via $\chi_{e,3}[{1\over p}]:\tilde Re[{1\over p}]\twoheadrightarrow K_3$). 

\medskip
Let $(M_1,F^1_1,\scrG_1,\tilde\mu_1):=(M_{\tilde Re},F^1_{\tilde Re},\tilde \scrG_{\tilde Re},\tilde\mu_{\tilde Re})\otimes_{\tilde Re} W(k)$. We have two identifications $M_1[{1\over p}]=M[{1\over p}]$ and $\scrG_{1,B(k)}=\scrG_{B(k)}$ and the pair $(M_1,\phi)$ is the Dieudonn\'e module over $k$ of the special fibre of $\tilde D_{V_3}$. As $\phi_{M_{\tilde Re}}(M_{\tilde Re}+{1\over p}F^1_{\tilde Re})=M_{\tilde Re}$, we have $\phi(M_1+{1\over p}F^1_1)$ and thus the cocharacter $\tilde\mu_1$ is a Hodge cocharacter of $(M_1,\phi,\scrG_1)$. Thus $(M_1,F^1_1,\phi,\scrG_1)$ is a Shimura filtered $F$-crystal over $k$. Moreover the triple
$$(M_{\tilde Re},F^1_{V_3},\phi_{M_{\tilde Re}},\tilde \scrG_{\tilde Re})$$
is a ramified lift of $(M_1,\phi,\scrG_1)$ to $V_3$; it is of $\text{CM}$ type, cf. the existence of the maximal torus $\scrT_{1,\tilde Re[{1\over p}]}$ of $\tilde\scrG_{\tilde Re[{1\over p}]}$. Let $h\in\scrG(B(k))$ be such that we have $h(M)=M_1$; for instance, we can take $h$ to be the composite of $\grJ_{\dbZ_p}\otimes 1_{W(k)}$ with the inverse of $\grJ_{\tilde Re}$ modulo the ideal $I(1)$ of $\tilde Re$ introduced in Subsection 2.8. Due to the property (d), it is easy to see that there exists an element $\tilde h\in\scrG(B(k))$ such that we have $F^1_1=\tilde h(F^1[{1\over p}])\cap M_1$. Thus we have $h\in\grI(\scrC)$ and an identity $(M_1,\phi,\scrG_1)=(h(M),\phi,\scrG(h))$. As  $h\in\grI(\scrC)$ and as $(M_{\tilde Re},F^1_{V_3},\phi_{M_{\tilde Re}},\tilde \scrG_{\tilde Re})$
is a ramified lift of $(M_1,\phi,\scrG_1)=(h(M),\phi,\scrG(h))$ to $V_3$ of $\text{CM}$ type, we have $[\scrC]\in Z^{\text{ram}}(\scrY(\scrF))$. Thus Theorem 4.1 (c) holds provided the $\dbQ_p$-linear isomorphism $\grJ$ exists.

\medskip\noindent
{\bf 5.3.2. End of the proof of Theorem 4.1 (c).} In this susubsection we show that there exists a $\dbQ_p$-linear isomorphism $\grJ$ as in Subsection 5.3 even if $\scrG^{\der}$ is not simply connected. As in [Va7, Subsubsect. 2.6.5], we consider a reductive group scheme $\scrG^{\flat}_{\dbZ_p}$ over $\dbZ_p$ equipped with an epimorphism $e:\scrG^{\flat}_{\dbZ_p}\twoheadrightarrow\scrG_{\dbZ_p}$ which has the following two properties:

\medskip
{\bf (i)} we have $\scrG^{\flat,\der}_{\dbZ_p}=\scrG^{\sc}_{\dbZ_p}$;

\smallskip
{\bf (ii)} the kernel $\Ker(e)$ is a finite product of tori of the form $\Res_{W(\dbF_{p^s})/\dbZ_p} \dbG_m$.

\medskip
We fix a closed embedding monomorphism $\scrG^{\flat}_{\dbZ_p}\hookrightarrow\pmb{GL}_{M^{\flat}_{\dbZ_p}}$. Let $\phi_0=i_M^{-1}\phi i_M:M_{\dbZ_p}\otimes_{\dbZ_p} W(k)\to M_{\dbZ_p}\otimes_{\dbZ_p} W(k)$ and $\mu_0=i_M^{-1}\mu i_M:\dbG_m\to \scrG_{W(k)}\leqslant \pmb{\GL}_{M_{\dbZ_p}\otimes_{\dbZ_p} W(k)}$. We have $\phi_0=g(1_{M_{\dbZ_p}}\otimes\sigma)\mu_0({1\over p})$, where $g\in \scrG_{\dbZ_p}(W(k))$. Due to (ii), there exist an element $g^{\flat}\in \scrG^{\flat}_{\dbZ_p}(W(k))$ that lifts $g$ and a cocharacter $\mu_0^{\flat}:\dbG_m\to \scrG^{\flat}_{W(k)}$ that lifts $\mu_0$. The triple
$$(M^{\flat}_{\dbZ_p}\otimes_{\dbZ_p} W(k),g^{\flat}(1_{M^{\flat}_{\dbZ_p}}\otimes\sigma)\mu_0^{\flat}({1\over p}),\scrG_{W(k)}^{\flat})$$
is a $p$-divisible object with a reductive group over $k$ in the sense of [Va7, Def. 1.2.1]. Let $\scrT_{1,B(k)}^{\flat}$ be the maximal torus of $\scrG_{B(k)}^{\flat}$ whose image in $\scrG_{B(k)}$ is $\scrT_{1,B(k)}$; its Lie algebra is $B(k)$-generated by elements fixed by $g^{\flat}(1_{M^{\flat}_{\dbZ_p}}\otimes\sigma)\mu_0^{\flat}({1\over p})$. Due to (ii), by performing the operation $\grO_1$ we can assume that $\scrG^{\flat}_{B(k)}$ is split and thus there exists also a cocharacter $\mu_1^{\flat}:\dbG_m\to \scrT_{1,K_1}^{\flat}$ that lifts $\mu_1$. 

One gets a filtered module over $K_2$
$$(M^{\flat}_{\dbZ_p}\otimes_{\dbZ_p} W(k),g^{\flat}(1_{M^{\flat}_{\dbZ_p}}\otimes\sigma)\mu_0^{\flat}({1\over p}),(F^{i,\flat}_{K_2})_{i\in\dbZ}),$$
where $F^{i,\flat}$ is the maximal $K_2$-vector subspace of $M^{\flat}_{\dbZ_p}\otimes_{\dbZ_p} K_2$ on which $\dbG_m$ acts via $\mu_{1,K_2}^{\flat}$ through weights at most $-i$. Subsubsections 5.2.1 to 5.2.3 can be entirely adapted to get that this filtered module over $K_2$ is admissible (the only difference is that one needs to consider functions $\bar f$ from $\{1,\ldots,l\}$ to a suitable finite subset of $\dbZ$). By assuming that $M_{\dbZ_p}$ is a direct summand of $M_{\dbZ_p}^{\flat}$ in such a way that $\scrG_{\dbZ_p}^{\flat}$ acts on it through  $\scrG_{\dbZ_p}$, one gets that the class $\gamma\in H^1(\dbQ_p,\scrG_{\dbQ_p})$ is the image of a class $\gamma^{\flat}\in H^1(\dbQ_p,\scrG^{\flat}_{\dbQ_p})$ which is associated similarly to $\gamma$ but obtained working with our new admissible filtered module over $K_2$. But as $\scrG^{\flat,\der}_{\dbQ_p}$ is simply connected (cf. (i)), the proof of Lemma 2.5.5 applies to give us that  $\gamma^{\flat}$ is the trivial class. Thus $\gamma$ itself is the trivial class and therefore a $\dbQ_p$-linear isomorphism $\grJ$ as in Subsection 5.3 does exist even if $\scrG^{\der}$ is not simply connected. Thus  $[\scrC]\in Z^{\text{ram}}(\scrY(\scrF))$ (cf. Subsubsection 5.3.1) i.e., Theorem 4.1 (c) holds.\endproof
 
\medskip\noindent
{\bf 5.3.3. Endomorphisms of $\tilde D_{V_3}$.} 
We refer to Subsubsection 5.3.1. Each element of $\End(M[{1\over p}])$ fixed by $\phi$ and $\scrT_{1,B(k)}$ defines an endomorphism of $(M[{1\over p}],\phi,F^1_{K_2})$ and therefore a $\dbQ$--endomorphism of $\tilde D_{V_3}$. Thus the homomorphism $\rho:\Gal(K_2)\to \pmb{\GL}_{\scrW}(\dbQ_p)$ and the homomorphism $\varrho$ of Subsubsection 5.3.1 factor through the group of $\dbZ_p$-valued points of the schematic closure in $\pmb{\GL}_{\scrW_{\dbZ_p}}$ of the torus of $\scrG^\prime_{\dbQ_p}$ that fixes the $\dbQ_p$-\'etale realizations of the $\dbQ$--endomorphisms of $\tilde D_{V_3}$ which correspond to those $t_{\alpha}$ with $\alpha\in\scrJ_1$.  Therefore $\tilde D_{V_3}$ is with complex multiplication.  This torus of $\scrG_{\dbQ_p}^\prime$ is canonically isomorphic to the torus $\scrT_{1,\dbQ_p}$ of Definition 2.4 (c) and thus it will be denoted as $\scrT_{1,\dbQ_p}$. 

\medskip\noindent
{\bf 5.3.4. Unramified version.} Referring to Subsubsection 5.3.3, we assume that the pair $(\scrT_{1,B(k)},\mu_1)$ is an unramified $E$-pair i.e., $\scrT_{1,B(k_1)}$ splits over the subfield $K$ of $B(\dbF)$. This implies that $K_2$ is unramified over $\dbQ_p$ and that $\scrT_{1,\dbQ_p}$ is unramified. Let $\scrT_{1,\dbZ_p}$ be the torus over $\dbZ_p$ whose generic fibre is $\scrT_{1,\dbQ_p}$. The group $\scrT_{1,\dbZ_p}(\dbZ_p)$ is the maximal compact subgroup of $\scrT_{1,\dbQ_p}(\dbQ_p)$ and therefore $\text{Im}(\rho)\leqslant\scrT_{1,\dbZ_p}(\dbZ_p)$. Based on this, we easily get that $\gamma\in H^1(\dbQ_p,\scrG_{\dbQ_p})$ is the image of a class in $H^1(\dbQ_p,\scrT_{1,\dbQ_p})$ and moreover that the proof of Lemma 2.5.5 can be adapted to give us that this last class is the trivial class. Based on Lemma 2.3, we can choose $\grJ:M_{\dbQ_p}\arrowsim \scrW$ such that $\scrW_{\dbZ_p}:=\grJ(M_{\dbZ_p})$ is normalized by $\scrT_{1,\dbZ_p}(\dbZ_p)$. Thus Subsubsection 5.3.1 can be redone in the unramified context in which we have $K_3=K_2$, $V_3$ is unramified over $\dbZ_p$, and we have a lift $A_{V_3}$ of $A$ whose $p$-divisible group is with complex multiplication and (to be compared with Subsection 2.6) a direct sum of $p$-divisible groups over $V_3=W(k)$ whose special fibres are isoclinic. Due to this last property, all these hold even if $p=2$. 

\bigskip\smallskip
\noindent
{\boldsectionfont 6. Proof of Basic Theorem 4.2}
\bigskip
 
In this section we prove Theorem 4.2. Let $(M_{\dbZ_p},\scrG_{\dbZ_p},(t_{\alpha})_{\alpha\in\scrJ})$ be the $\dbZ_p$ structure of $(M,\phi,\scrG,(t_{\alpha})_{\alpha\in\scrJ})$, cf. Subsection 2.5. Each simple factor of $\scrG^{\ad}_{\dbF_p}$ is of the form $\Res_{k_0/\dbF_p} \scrG^0_{k_0}$, where $k_0$ is a finite field and $\scrG^0_{k_0}$ is an absolutely simple, adjoint group over $k_0$ (cf. [Ti1, Prop. 3.1.2]). Thus each simple factor of $\scrG^{\ad}_{\dbZ_p}$ is of the form $\Res_{W(k_0)/\dbZ_p} \scrG^0$, where $\scrG^0$ is an absolutely simple, adjoint group scheme over $W(k_0)$ whose special fibre is $\scrG^0_{k_0}$ (cf. [DG, Vol. III, Exp. XXIII, Prop. 1.21]). Until Section 7 we will assume that each such group scheme $\scrG^0$ is of $B_n$, $C_n$, or $D_n$ Dynkin type.

We assume that there exists a maximal torus $\scrT_{1,B(k)}$ of $\scrG_{B(k)}$ of $\dbQ_p$-endomorphisms of $\scrC$. Let $\scrT_{1,\dbQ_p}$, $K$, $K_1$, and $K_2$ be as in Definition 2.4 (c). To prove Theorem 4.2, we will show that there exists a cocharacter $\mu_1$ of $\scrT_{1,K_1}$ with the property that the product of the cocharacters of $\scrT_{1,K_1}$ which belong to the $\Gal(K_1/\dbQ_p)$-orbit of $\mu_1$, factors through the subtorus $Z^0(\scrG_{K_1})$ of $\scrT_{1,K_1}$. Based on Example 2.4.1, this will implies that the $E$-pair $(\scrT_{1,B(k)},\mu_1)$ satisfies the $\grC$ condition and thus that $TT\grR$ holds for $\scrC$. The condition on the $\Gal(K_1/\dbQ_p)$-orbit of $\mu_1$ does not change if we replace $K_1$ by a larger Galois extension of $\dbQ_p$ (such as $K_2$) and thus, by performing the operation $\grO_1$ we can assume that $\scrG$ is split. Therefore the field $k$ contains each such field $k_0$.

Until Section 7 we will also assume that $L_{\scrG}^0(\phi)$ is a Levi subgroup scheme of $P_{\scrG}^+(\phi)$ and that $\mu:\dbG_m\to\scrG$ factors through a maximal torus $\scrT$ of $L^0_{\scrG}(\phi)$ contained in a Borel subgroup scheme $\scrB$ of $\scrP$, cf. Subsection 2.6. Let $L_{\scrG}^0(\phi)_{\dbZ_p}$ be the $\dbZ_p$ structure of $L_{\scrG}^0(\phi)$ obtained as in Subsection 2.5. It is a reductive, closed subgroup scheme of $\scrG_{\dbZ_p}$ which is the centralizer of the at most rank $1$ split torus of $\scrG_{\dbZ_p}$ whose extension to $B(k)$ is the image of the Newton cocharacter of $\scrC$, cf. Fact 2.6.1 (a) and its proof. 

In Subsection 6.1 we include some reduction steps. 
In Subsection 6.2 we include few simple properties. In Subsection 6.3 (resp. Subsection 6.4) we deal with the cases related to Shimura pairs of $B_n$ and $D_n^{\dbR}$ (resp. $C_n$ and $D_n^{\dbH}$) type. The proof of Theorem 4.2 ends in Subsection 6.5.

\bigskip\noindent
{\bf 6.1. Some reductions and notations.} Let $\mu^{\ad}:\dbG_m\to\scrG^{\ad}$ be the composite of the cocharacter $\mu:\dbG_m\to\scrG$ with the natural epimorphism $\scrG\twoheadrightarrow \scrG^{\ad}$. As $\scrG$ is split, each cocharacter of $\scrG^{\ad}_{K_2}$ which is $\scrG^{\ad}(K_2)$-conjugate to $\mu_{K_2}^{\ad}$, lifts uniquely to a cocharacter of $\scrG_{K_2}$ which is $\scrG_{K_2}(K_2)$-conjugate to $\mu_{K_2}$. Below we will consider only $E$-pairs of $\scrC$ which are as in Example 2.4.1. Thus based on the last two sentences, on Subsubsection 2.5.6, and on the fact that the statements 4.2 (a) and (b) pertain only to images in $\scrG^{\ad}_{K_2}$ of suitable products of cocharacters of $\scrG_{K_2}$ that factor through $\scrT_{1,K_2}$, to prove Theorem 4.2 we can assume that the adjoint group scheme $\scrG_{\dbZ_p}^{\ad}$ is $\dbZ_p$-simple and that the cocharacter $\mu^{\ad}$ is non-trivial. Thus $\scrG_{\dbZ_p}^{\ad}=\Res_{W(k_0)/\dbZ_p} \scrG^0$. If $\scrG^0$ is of $D_n$ Dynkin type, then it splits over $W(k_{02})$, where $k_{02}$ is the quadratic extension of $k_0$ (cf. [Se1, Cor. 2 of p. 182]). 

Let $m\in\dbN^*$ be such that $k_0:=\dbF_{p^m}$. We write 
$$\scrG^{\ad}=\scrG_1\times_{W(k)}\scrG_2\times_{W(k)}\cdots \times_{W(k)}\scrG_m,$$  
where $\scrG_i$ is a split, absolutely simple, adjoint group scheme over $W(k)$ and the numbering of $\scrG_i$'s is such that we have $\phi(\Lie(\scrG_i[{1\over p}]))=\Lie(\scrG_{i+1}[{1\over p}])$ for all $i\in\{1,\ldots,m\}$. Here and in all that follows in this section, the left lower or upper index $m+1$ has the same role as $1$ (thus $\scrG_{m+1}:=\scrG_1$, etc.). Let $\scrG^i$ be the semisimple,  normal, closed subgroup scheme of $\scrG^{\der}$ which is naturally isogenous to $\scrG_i$. We view the isomorphism (2) of Subsubsection 2.5.1 as an identification and therefore we can write $\phi=g(1_{M_{\dbZ_p}}\otimes\sigma)\mu({1\over p})$, where $g\in \scrG(W(k))$ (to be compared with Subsubsection 2.5.6). Let $g^{\ad}\in \scrG^{\ad}(W(k))$ be the image of $g$. For $\alpha\in\dbQ$ let $\grD_{\alpha}$ be the central division algebra over $B(k_0)$ of invariant $\alpha$. 

\medskip\noindent
{\bf 6.1.1. Fact.} {\it To prove the Theorem 4.2 we can also assume that we have a direct sum decomposition $M=\oplus_{i=1}^m M_i$
into $\scrG$-modules such that the following two conditions hold:

\medskip
{\bf (i)} if $i,j\in\{1,\ldots,m\}$ with $j\neq i$, then $M_i$ has no trivial $\scrG^i$-submodule and $M_j$ is a trivial $\scrG^i$-module; 

\smallskip
{\bf (ii)}  we have an identity $Z(\scrG)=\prod_{i=1}^m Z(\pmb{GL}_{M_i})$.}

\medskip
\proof
The arguments for this fact are similar to the ones of the proof of [Va1, Thm. 6.5.1.1 or Subsubsect. 6.6.5] but simpler as we are over $\dbZ_p$ and not over $\dbZ_{(p)}$ and we do not have to mention polarized Hodge $\dbQ$--structures. We recall the essence of loc. cit. 

We first assume that $\scrG^0$ is of $B_n$ Lie type. Thus $\scrG^0$ is split. We consider the spin faithful representation $\scrG^{0\sc}\hookrightarrow \pmb{\GL}_{M_0}$ over $W(k_0)$. Let $\pmb{\text{GSpin}}$ be the closed subgroup scheme of $\pmb{\GL}_{M_0}$ generated by $\scrG^{0\sc}$ and $Z(\pmb{\GL}_{M_0})$. Let $\scrG_{\dbZ_p}^\prime:=\Res_{W(k_0)/\dbZ_p} \pmb{\text{GSpin}}$. We consider its faithful representation on $M_{\dbZ_p}^\prime$, where $M_{\dbZ_p}^\prime$ is $M_0$ but viewed as a $\dbZ_p$-module. We identify naturally $\scrG_{\dbZ_p}^{\ad}=\scrG^{\prime\ad}_{\dbZ_p}$. Let $M^\prime:=M_{\dbZ_p}^\prime\otimes_{\dbZ_p} W(k)$. Let $\scrG^\prime:=\scrG_{W(k)}^\prime$. We have a unique direct sum decomposition $M^\prime=\oplus_{i=1}^m M_i^\prime$ of $\scrG^\prime$-modules which are also $W(k_0)\otimes_{\dbZ_p} W(k)$-modules. Let $g^\prime\in \scrG^\prime (W(k))$ be such that its image in $\scrG^{\prime\ad}(W(k))$ is $g^{\ad}$. Let $\mu^\prime$ be a cocharacter of $\scrG^\prime$ such that the cocharacter of $\scrG^{\prime\ad}=\scrG^{\ad}$ it defines naturally is $\mu^{\ad}$ and the triple $\scrC^\prime:=(M^\prime,g^{\prime}(1_{M_{\dbZ_p}}^{\prime}\otimes\sigma)\mu^\prime({1\over p}),\scrG^\prime)$ is a Shimura $F$-crystal over $k$. Let $\scrT_{1,B(k)}^\prime$ be the maximal torus of $\scrG^{\prime}_{B(k)}$ whose image in $\scrG^{\prime\ad}_{B(k)}$ is the same as of $\scrT_{1,B(k)}$; it is a maximal torus of $\scrG^{\prime}_{B(k)}$ of $\dbQ_p$-endomorphisms of $\scrC^\prime$. Thus $\scrC^\prime$ is semisimple, cf. Fact 3.1 (b). Let $\mu_1^\prime$ be the cocharacter of $\scrT_{K_1}^\prime$ which over $K_2$ is $\scrG^\prime(K_2)$-conjugate to $\mu^\prime_{K_2}$ and such that it defines the same cocharacter of $\scrG^{\prime\ad}_{K_1}=\scrG^{\ad}_{K_1}$ as $\mu_1$. The $E$-pair $(\scrT_{1,B(k)}^\prime,\mu_1^\prime)$ of $\scrC^\prime$ satisfies the $\grC$ condition if and only if the $E$-pair $(\scrT_{1,B(k)},\mu_1)$ of $\scrC$ satisfies the $\grC$ condition. Similarly, $TT\grR$ holds for $\scrC^\prime$ if and only if $TT\grR$ holds for $\scrC$. Thus to prove Theorem 4.2 for the case when $\scrG^0$ is of $B_n$ Lie type, we can replace $\scrC$ by $\scrC^\prime$. As the two conditions (i) and (ii) obviously hold if $\scrC$ is $\scrC^\prime$, the fact holds if $\scrG^0$ is of $B_n$ Lie type.

If $\scrG^0$ is of $C_n$ or $D_n$ Dynkin type, we will only list the modifications required to be performed to the previous paragraph. If $\scrG^0$ is of $C_n$ Lie type, then the spin representation has to be replaced by the standard rank $2n$ faithful representation $\scrG^{0\sc}\hookrightarrow \pmb{\GL}_{M_0}$ over $W(k_0)$. If $\scrG^0$ is of $D_n$ Dynkin type, then we have two disjoint subcases related to Shimura pairs of $D_n^{\dbH}$ and $D_n^{\dbR}$ (respectively) type, with the $D_n^{\dbR}$ type as in Definition 1.2.4 (a). 

In the first subcase the spin representation has to be replaced by the standard rank $2n$ faithful representation $\scrG^{01}\hookrightarrow \pmb{\GL}_{M_0}$ over $W(k_0)$. Here $\scrG^{01}$ is an isogeny cover of $\scrG^0$ for which such a representation is possible; its existence is implied by the fact that $\scrG^0$ splits over $W(k_{02})$. If $n>4$, then $\scrG^{01}$ is unique. If $n=4$, then we choose $\scrG^{01}$ such that the construction of $\mu^{\prime}$ is possible (we have only one choice for $\scrG^{01}$, due to the fact that the two subcases are disjoint).  The second subcase is in essence the same as the previous paragraph (the only difference being that $\scrG^{0\sc}$ is not necessarily split; however, as it splits over $W(k_{02})$, its spin representation is well defined over $W(k_0)$).\endproof

\bigskip\noindent
{\bf 6.2. Simple properties.} We first consider the case when $L_{\scrG}^0(\phi)$ is a torus (i.e., we have $L_{\scrG}^0(\phi)=\scrT$). Thus we have $\scrT_{1,B(k)}=\scrT_{B(k)}$ and $K_1\subseteq B(k)$. Let $\tau_1\in \Gal(K_1/\dbQ_p)$ be the restriction of $\sigma$ to $K_1$. The $E$-triple $(\scrT_{1,B(k)},\mu_{B(k)},\tau=(\tau_1))$ satisfies the condition 2.4 (e1) and is obviously admissible. Thus Theorem 4.2 holds if $L_{\scrG}^0(\phi)$ is a torus. From now on until Section 7 we will assume that $L^0_{\scrG}(\phi)$ is not a torus (i.e., we have $L^0_{\scrG}(\phi)\neq \scrT$). Let $\scrL_0$ be the $\dbQ_p$-form of $L^0_{\scrG}(\phi)_{B(\bar k)}$ with respect to $(M[{1\over p}],\phi)$. The tori $Z^0(\scrG_{\dbQ_p})$ and $\scrT_{1,\dbQ_p}$ are subgroups of $\scrL_0$. We have a direct sum decomposition
$$\Lie(\scrL_0)\otimes_{\dbQ_p} B(k_0)=(\Lie(Z^0(\scrL_0))\otimes_{\dbQ_p} B(k_0))\bigoplus_{i=1}^m \grL_0^i,\leqno (5)$$
where $\grL_0^i:=(\Lie(\scrL_0^{\der})\otimes_{\dbQ_p} B(k_0))\cap\Lie(\scrG_{iB(\bar k)})$. Each $\grL_0^i$ is a semisimple Lie algebra and thus it is also the Lie algebra of a semisimple group $\scrL_0^i$ over $B(k_0)$. Moreover, we have $\phi(\grL_0^i)=\grL_0^{i+1}$. Based on this and (5) we get that each $p$-adic field over which $\scrL_0$ splits must contain $B(k_0)$. Therefore $B(k_0)\subseteq K_1$. 

\medskip\noindent
{\bf 6.2.1. Lemma.} {\it We assume that $p\Ge 3$. Let $H_1:=\Gal(K_1/\dbQ_p)$. Let $H_0$ be a subgroup of $H_1^0:=\Gal(K_1/B(k_0))$ of even index. If $m$ is odd, then there exists an element $\tau_1\in H_1$ such that the following two conditions hold:

\medskip
{\bf (i)} all orbits under $\tau_1^m$ of the left translation action of $H_1$ on $H_1/H_0$ have an even number of elements;

\smallskip
{\bf (ii)} the action of $\tau_1$ on the residue field $l_1$ of $K_1$ is the Frobenius automorphism of $l_1$ whose fixed field is $\dbF_p$.}

\medskip
\proof
For $s\in\dbN$ let $H_{1s}$ be the $s$-th ramification group of $H_1$. Thus $H_1=H_{10}$, $H_1/H_{11}$ is cyclic, and the subgroup $H_{12}$ of $H_1$ is normal and (as $p\Ge 3$) has odd order. By replacing $H_1$ with $H_1/H_{12}$,  we can assume that $H_{12}=\{1_{K_1}\}$. Thus $H_{11}$ is a subgroup of $\dbG_m(l_1)$ and therefore it is cyclic. By replacing $H_1$ with its quotient through a normal subgroup of $H_{11}$ of odd order, we can assume that $H_{11}$ is of order $2^t$ for some $t\in\dbN$. The case $t=0$ is trivial and therefore we can assume that $t\Ge 1$. Let $H_{01}$ be the image of $H_0$ in $H_1^0/H_{11}$ and let $a$ be its index in $H_1^0/H_{11}$. If $a$ is even,  then the condition (i) is implied by (ii) and therefore we can choose any element $\tau_1\in H_1$ for which the condition (ii) holds. If $a$ is odd, then by replacing $H_1$ with its quotient through the subgroup of $H_{11}$ of order $2^{t-1}$ we can assume that $t=1$. Thus $H_{11}$ has order 2. As $H_{11}$ is a normal subgroup of $H_1$ of order $2$, it is included in the center of $H_1$. Thus $H_1$ is either cyclic or isomorphic to $H_{11}\times H_1/H_{11}$. If $H_1$ is cyclic, then we can take $\tau_1\in H_1$ such that it generates $H_1$ and the condition (ii) holds. If $H_1$ is isomorphic to $H_{11}\times H_1/H_{11}$, then we can take $\tau_1=(\tau_{11},\tau_{12})$ such that $\tau_{11}\in H_{11}$ and $\tau_{12}\in H_1/H_{11}$ generate these groups and the condition (ii) holds. In both cases the condition (i) also holds. \endproof   

\medskip\noindent
{\bf 6.2.2. Factors.} Let $\grN$ be the set of those elements $i\in\{1,\ldots,m\}$ for which the image of $\mu^{\ad}:\dbG_m\to\scrG^{\ad}$ in $\scrG_i$ is non-trivial. Let $\grM:=\{1,\ldots,m\}\setminus \grN$. If $i\in\grN$ (resp. $i\in\grM$), then $\scrG_i$ is called a {\it non-compact} (resp. {\it compact}) {\it factor} of $\scrG^{\ad}$ with respect to $\mu^{\ad}$. As $\mu^{\ad}$ is non-trivial, the set $\grN$ is non-empty.  Let $v\in\dbN^*$ be the number of elements of $\grN$. To simplify notations we will assume that $1\in \grN$. We will always choose the cocharacter $\mu:\dbG_m\to\scrG$ such that $\dbG_m$ acts via $\mu$ trivially on $M_i$ for all $i\in \grM$.

\bigskip\noindent
{\bf 6.3. Case 1.} Until Subsection 6.4 we will assume that $\scrG^0$ is of either $B_n$ or $D_n$ Dynkin type and that $\scrG^{\der}$ is simply connected (under these assumptions, one can also assume that $\scrG^i\to\pmb{GL}_{M_i}$ is the spin representation). Let $\scrG^{01}\to \scrG^0$ be an isogeny such that $\scrG^{01}$ is the \pmb{\text{SO}} group scheme of a quadratic form on a free $W(k_0)$-module $O_0$ of rank $t$. Here $t=2n+1$ (resp. $t=2n$) if $\scrG^0$ is of $B_n$ (resp. $D_n$) Dynkin type. Let 
$$\pmb{\text{S}}:=(\Res_{W(k_0)/\dbZ_p} \scrG^{01})\times_{\dbZ_p} \Spec W(k);$$ 
it is a semisimple group scheme over $W(k)$ whose adjoint group scheme is $\scrG^{\ad}$. Let $\mu_0:\dbG_m\to \pmb{\text{S}}$ be the unique cocharacter that lifts $\mu^{\ad}$ (for instance, if $\scrG^0$ is of $B_n$ Dynkin type, then $\mu_0=\mu^{\ad}$ as $\pmb{\text{S}}=\scrG^{\ad}$). Let $\pmb{\text{SO}}(O_i,b_i):=\scrG^{01}\times_{W(k_0)} \Spec W(k)$, where the $\dbZ_p$-embedding $W(k_0)\hookrightarrow W(k)$ is the same as the one that defines $\scrG_i=\scrG^0\times_{W(k_0)} \Spec W(k)$ and where $b_i$ is a perfect symmetric bilinear form on the $W(k)$-module $O_i$. We have $\pmb{\text{S}}=\prod_{i=1}^m \pmb{\text{SO}}(O_i,b_i)$. We have an identification $O_0\otimes_{\dbZ_p} W(k)=\oplus_{i=1}^m O_i$ of $W(k_0)\otimes_{\dbZ_p} W(k)$-modules. For $i\in\{1,\ldots,m\}$, let $W_i:=O_i\otimes_{W(k)} K_2$.

Let $\varUpsilon_i:=\{e_1^i,\ldots,e_t^i\}$ be a $K_2$-basis for $W_i$ such that the following three things hold:

\medskip
{\bf (i)} if $a,b\in\{1,\ldots,t\}$ with $a<b$, then the value of $b_i(e_a^i,e_b^i)$ is $0$ or $1$ depending on the fact that the pair $(a,b)$ belongs or not to the set $\{(1,2),\ldots,(2n-1,2n)\}$;

\smallskip
{\bf (ii)} the split torus $\scrT_{1,K_2}$ normalizes each $K_2e_a^i$;

\smallskip
{\bf (iii)} if $\dbG_m$ acts through $\mu_{1,K_2}$ on the $K_2$-span of some element $e^i_a$ with weight $-1$ (resp. $1$), then $a$ is odd (resp. is even).  

\medskip
The natural action of $\Gal(K_2/\dbQ_p)$ on cocharacters of $\scrT_{1,K_2}$ defines naturally an action of $\Gal(K_2/\dbQ_p)$ on $\varUpsilon:=\cup_{i=1}^m \varUpsilon_i$. For $\star\in\Gal(K_2/\dbQ_p)$, let $\pi_{\star}$ be the permutation of $\varUpsilon$ defined by $\star$. For each $i\in\{1,\ldots,m\}$, the set $\varUpsilon_i$ is normalized by $\Gal(K_2/B(k_0))$.

If $p>2$, then up to a replacement of $k$ by an extension of it of degree at most $2$, there exists an element $g_0\in \pmb{\text{S}}(W(k))$ whose image in $\scrG^{\ad}(W(k))$ is $g^{\ad}$. If $p=2$, then after a similar such replacement, the element $g_0$ exists provided we replace $g^{\ad}$ by the image in $\scrG^{\ad}(W(k))$ of a suitable element $hg\phi(h^{-1})$, where $h\in\scrG(W(k))$ normalizes $F^1/pF^1$ (to be compared with [Va7, Fact 2.6.3]). Thus by performing the operation $\grO_1$, we can assume that  there exists an element $g_0\in \pmb{\text{S}}(W(k))$ whose image in $\scrG^{\ad}(W(k))$ is $g^{\ad}$. Let $\phi_0:=g_0(1_{O_0}\otimes \sigma)\mu_0({1\over p})$; it is a $\sigma$-linear automorphism of $O_0\otimes_{\dbZ_p} B(k)$. 

Let $\scrT_{0,B(k)}$ be the maximal torus of $\pmb{\text{S}}_{B(k)}$ whose image in $\pmb{\text{S}}_{B(k)}^{\ad}=\scrG^{\ad}_{B(k)}$ is $\scrT_{0,B(k)}^\prime:=\text{Im}(\scrT_{1,B(k)}\to \scrG^{\ad}_{B(k)})$. Let $\scrT_{0,\dbQ_p}$ be the $\dbQ_p$-form of $\scrT_{0,B(k)}$ with respect to $(O_0\otimes_{\dbZ_p} B(k),\phi_0)$; its Lie algebra is $\{x\in \Lie(\scrT_{0,B(k)})=\Lie(\scrT_{0,B(k)}^\prime)|\phi_0(x)=x\}$. Let $\scrT_{0,\dbQ_p}^\prime$ be the $\dbQ_p$-form of $\scrT_{0,B(k)}^\prime$ which is the quotient of $\scrT_{0,\dbQ_p}$ by its finite subgroup whose extension to $B(k)$ is $\scrT_{0,B(k)}\cap Z(\pmb{\text{S}}_{B(k)})$. Equivalently, we have $\scrT_{0,\dbQ_p}^\prime=\scrT_{1,\dbQ_p}/Z(\scrG_{\dbQ_p})$.

Until Subsubsection 6.3.4 we will assume that $\scrC$ is basic. Thus $L_{\scrG}^0(\phi)=\scrG$ and therefore $\scrL_0$ is a $\dbQ_p$-form of $\scrG_{B(\bar k)}$. To show that there exists an $E$-pair $(\scrT_{1,B(k)},\mu_1)$ of $\scrC$ as in Example 2.4.1, we first prove the following lemma.

\medskip\noindent
{\bf 6.3.1. Lemma.} {\it  We recall that $\scrC$ is basic. The action of $\Gal(K_2/B(k_0))$ on $\varUpsilon_1$ has an orbit that contains $\{e_{2a_1-1}^1,e_{2a_1}^1\}$ for some element $a_1\in\{1,\ldots,n\}$.}

\medskip
\proof
In this proof by orbit we will mean an orbit of the action of $\Gal(K_2/B(k_0))$ on $\varUpsilon_1$. We consider an orbit $\tilde o_1$ whose  elements are pairwise perpendicular with respect to the bilinear form $b_1$. To fix the notations, we can assume that there exists $a\in\{1,...,n\}$ such that $e_{2a-1}^1\in\tilde o_1$. Let $\tilde o_2$ be the orbit that contains $e_{2a}^1$; we have $\tilde o_1\neq \tilde o_2$. The orbit decomposition of $\varUpsilon_1$ corresponds to a direct sum decomposition of $O_1[{1\over p}]$ in minimal $B(k)$-vector subspaces normalized by $\scrT_{1,B(k)}$ and $\phi_0^m$. Let $O_{1,1}$ and $O_{1,2}$ be the $B(k)$-vector subspaces of $O_1[{1\over p}]$ that correspond to $\tilde o_1$ and $\tilde o_2$ (respectively). The intersection $\grL_0^1\cap (\End(O_{1,1}\oplus O_{1,2})\otimes_{B(k)} B(\bar k))$ is the Lie algebra of a split semisimple group  over $B(k_0)$ of $D_s$ Lie type, where $s$ is the number of elements of $\tilde o_1$. Argument: this is implied by the fact that (as $\scrC$ is basic), for $j\in\{1,2\}$ all Newton polygon slopes of $(O_{1,j}[{1\over p}],\phi_0^m)$ are $0$. 

Thus if the lemma does not hold, then $\scrL_0^1$ has a split torus of rank equal to the sum of such $s$'s and thus of rank $n$ and therefore it is a split group over $B(k_0)$. Thus to prove the lemma we only have to show that the semisimple group $\scrL_0^1$ is in fact non-split. This is a rational statement. Thus to check it, based on [Va7, Thm. 1.3.3 and Subsect. 2.5] we can assume that $\phi(\Lie(\scrT))=\Lie(\scrT)$ and (in order to use the above notations on $\varUpsilon_i$'s) that $\scrT_{1,B(k)}=\scrT_{B(k)}$. Thus $K_2=B(k)$ and the actions of $\pi_{\sigma_{\dbF_p}}$ and $\phi$ on cocharacters of $\scrT$ coincide. We can also assume that $\mu_0:\dbG_m\to \pmb{\text{S}}$ fixes $e^i_a$ if either $i\in \grM$ or $a\Ge 3$. Let $g_1\in N_{\scrG}(\scrT)(W(k))$ be such that $g_1\phi(\Lie(\scrB))\subseteq\Lie(\scrB)$, cf. [Va7, Subsect. 2.5]. Let $w_1\in N_{\scrG}(\scrT)(W(k))$ be such that its image in $\scrG^{\ad}(W(k))$ belongs to $\scrG_1(W(k))$ and takes the image of $\scrB$ in $\scrG_1$ to its opposite with respect to the image of $\scrT$ in $\scrG_1$. As $\scrG^{\der}$ is simply connected of either $B_n$ or $D_n$ Dynkin type, $w_1$ takes the cocharacter of $\scrG_1$ defined by $\mu^{\ad}$ to its inverse. Thus the Shimura $F$-crystal $(M,w_1g_1\phi,\scrG)$ over $k$ is basic (to be compared with [Va7, Cases 1 and 2 of Subsubsect. 4.2.2]). Thus based on [Va7, Prop. 2.7] we can assume that $w_1g_1=1_M$ and therefore that:

\medskip
{\bf (i)} $\pi_{\sigma_{\dbF_p}}^m$ restricted to $\varUpsilon_1$ fixes $e_a^1$ for $a\Ge 3$ and permutes $e_1^1$ and $e_2^1$. 

\medskip
Thus the group $\scrL_0^1$ has a split torus $T^1_{1;n-1}$ of rank $n-1$: it is the torus of $\pmb{\text{S}}$ that fixes $e^1_1$ and $e^1_2$ and that normalizes $B(k_0)e^1_a$ for each $a\in\{3,\ldots,t\}$. The centralizer of $T^1_{1;n-1}$ in $\scrL_0^1$ is a non-split torus, cf. property (i). Thus the group $\scrL_0^1$ is non-split.\endproof

\medskip\noindent
{\bf 6.3.2. The choice of $\mu_1$ for the basic context.} Let $a_1$ be as in Lemma 6.3.1. Let $o$ be the orbit of $e_{2a_1-1}^1$ under $\Gal(K_2/\dbQ_p)$. For $i\in \grN\setminus\{1\}$ let $a_i\in\{1,\ldots,n\}$ be such that $\{e_{2a_i-1}^i,e_{2a_i}^i\}\subseteq o$. Let $\mu_2:\dbG_m\to \pmb{\text{S}}_{K_2}$ be the cocharacter that fixes all $e_a^i$'s except those of the form $e_{2a_i-1+j}^i$, where $i\in \grN$ and $j\in\{0,1\}$, and that acts as the identical (resp. as the inverse of the identical) character of $\dbG_m$ on each $K_2e_{2a_i}^i$ (resp. $K_2e_{2a_i-1}^i$) with $i\in \grN$; the cocharacter $\mu_2$ is $\pmb{\text{S}}(K_2)$-conjugate to $\mu_{0,K_2}$ (cf. also property 6.3 (iii) in the $D_n$ Dynkin type case). To define the cocharacter $\mu_1:\dbG_m\to\scrT_{1,K_1}$ it is enough to define the cocharacter $\mu_{1,K_2}:\dbG_m\to\scrT_{1,K_2}$. Let $\mu_{1,K_2}:\dbG_m\to \scrT_{1,K_2}$ be the unique cocharacter such that the cocharacter of $\scrG^{\ad}_{K_2}$ (resp. of $\scrG^{\ab}_{K_2}$) it defines naturally is the composite of $\mu_2$ with the isogeny $\pmb{\text{S}}_{K_2}\to \scrG^{\ad}_{K_2}$ (resp. is the one defined by $\mu_{K_2}$). As $\{e_{2a_i-1}^i,e_{2a_i}^i\}\subseteq o$, the product of the cocharacters of $\scrT_{1,K_2}$ that belong to the $\Gal(K_2/\dbQ_p)$-orbit of $\mu_{1,K_2}$ has a trivial image in $\scrG^{\ad}_{K_2}$ and thus it factors through $Z^0(\scrG_{K_2})$.

\medskip\noindent
{\bf 6.3.3. Remark.} If $m$ is odd and the action of $\Gal(K_2/B(k_0))$ on $\varUpsilon_1$ has only one orbit, then it is easy to check  based on Lemma 6.2.1 that there exists an $E$-pair $(\scrT_{1,B(k)},\mu_1)$ of $\scrC$ that satisfies the cyclic $\grC$ condition.

\medskip\noindent
{\bf 6.3.4. The non-basic context.} Until the Case 2 in Subsection 6.4 we will assume that $\scrC$ is non-basic. We use the previous notations of Subsection 6.3. Let $\scrT_{0,B(k)}^0$ be the subtorus of $\scrT_{0,B(k)}$ whose image in $\scrG^{\ad}_{B(k)}$ is the same as the image of $Z^0(L_{\scrG}^0(\phi))_{B(k)}$. As $L_{\scrG}^0(\phi)_{\dbZ_p}$ is the centralizer in $\scrG_{\dbZ_p}$ of a rank 1 split torus whose generic fibre is canonically identified with a subtorus of $\scrT_{1,\dbQ_p}$ (cf. beginning of Section 6) and whose extension to $B(k)$ is the image of the Newton cocharacter of $\scrC$, the group $C_{\pmb{\text{S}}_{B(k)}}(\scrT_{0,B(k)}^0)$ is a product 
$$\prod_{i=1}^m\prod_{j=1}^{m_f} \pmb{\text{SO}}(O_{i,j},b_{i,j}),\leqno (6)$$ 
where $O_i[{1\over p}]=\oplus_{j=1}^{m_f} O_{i,j}$ is the minimal direct sum decomposition normalized by $\scrT_{0,B(k)}^0$ and $b_{i,j}$ is the restriction of $b_i$ to $O_{i,j}$. We emphasize that $m_f\in\dbN^*$ does not depend on $i$ and that, in the case when $b_{i,j}=0$, we define $\pmb{\text{SO}}(O_{i,j},b_{i,j}):=\pmb{\GL}_{O_{i,j}}$. We choose the indexes such that we have $\phi(\End(O_{i,j}))=\End(O_{i+1,j})$. The direct sum decomposition 
$$(O_0[{1\over p}],\phi_0)=\oplus_{j=1}^{m_f} (\oplus_{i=1}^m O_{i,j}[{1\over p}],\phi_0)$$ 
is the Newton polygon slope decomposition. 

We will define a cocharacter $\mu_2:\dbG_m\to \pmb{\text{S}}_{K_2}$ that factors through $\scrT_{0,K_2}$. Let $j\in \{1,\ldots,m_f\}$. We first assume that $b_{1,j}\neq 0$. If $\dbG_m$ acts via $\mu_0$ trivially (resp. non-trivially) on $\oplus_{i=1}^m O_{i,j}$, then we define the action of $\dbG_m$ via $\mu_2$ on $O_{i,j}\otimes_{B(k)} K_2$ to be trivial (resp. to be obtained as in Subsubsection 6.3.2 but working with $(\oplus_{i=1}^m O_{i,j},\phi_0)$ instead of with $(\oplus_{i=1}^m O_i,\phi_0)$). 

Until the Case 2 we assume that $b_{1,j}=0$. Let $j_{\perp}\in\{1,\ldots,m_f\}\setminus\{j\}$ be the unique element such that $O_{i,j_{\perp}}$ is not perpendicular on $O_{i,j}$ with respect to $b_i$. Let $o_1,\ldots,o_s$ be the orbits of the action of $\Gal(K_2/\dbQ_p)$ on $\varUpsilon\cap [(\oplus_{i=1}^m O_{i,j})\otimes_{B(k)} K_2]$. Let $d_j:=\dim_{B(k)}(O_{1,j})$. Let $S_{+,j}$ (resp. $S_{-,j}$) be the set of those elements $i\in \grN$ with the property that $\mu_0$ acts via the inverse of the identical (resp. via the identical) character of $\dbG_m$ on a non-zero element of $O_{i,j}$. Let $c_{+,j}$ (resp. $c_{-,j}$) be the number of elements of $S_{+,j}$ (resp. of $S_{-,j}$). We have $S_{+,j}\cap S_{-,j}=\emptyset$ as otherwise $b_{1,j}\neq 0$. Thus $c_{+,j}+c_{-,j}\le v$ (see Subsubsection 6.2.2 for $v$). 

Let $\zeta_j\in\dbZ$ and $\xi_j\in\dbN^*$ be such that $g.c.d.(\zeta_j,\xi_j)=1$ and $(c_{+,j}-c_{-,j})\xi_j=d_j\zeta_j$. The only Newton polygon slope of $(O_{i,j},\phi_0^m)$ is ${{\zeta_j}\over {\xi_j}}\in [-1,1]$. Thus $d_j\in \xi_j\dbN^*$ and the $B(k_0)$-subalgebra
 $\{x\in\End(O_{i,j})\otimes_{B(k)} B(\bar k)|(\phi\otimes\sigma_{\bar k}^m(x)=x\}$ is (isomorphic to) $M_{d_j/\xi_j}(\grD_{{\zeta_j}\over {\xi_j}})$, cf. Dieudonn\'e--Manin classification of $F$-isocrystals over $\bar k$ (see [Ma, Sect. 2]). Thus:

\medskip
{\bf (i)} for each $l\in\{1,\ldots,s\}$ there exists $e_l\in\dbN^*$ such that the number of elements of $o_l\cap\varUpsilon_1$ is $\xi_je_l$; moreover, we have $\xi_j(\sum_{l=1}^s e_l)=d_j$. 

\medskip
For each $l\in\{1,\ldots,s\}$ we choose numbers $c_{l,j}^+$, $c_{l,j}^-\in\dbN$ such that the following two relations hold:

\medskip
{\bf (ii)} $c^+_{l,j}-c^-_{l,j}=\zeta_je_l$;

\smallskip
{\bf (iii)} $\sum_{l=1}^s c^+_{l,j}=c_{+,j}$ (and thus $\sum_{l=1}^s c^-_{l,j}=c_{-,j}$).

\medskip
For instance, if $\zeta_j\Ge 0$ we can choose $c^-_{1,j}=\cdots=c^-_{s-1,j}=0$ and $c^-_{s,j}=c_{-,j}$, the numbers $c^+_{l,j}$'s being now determined uniquely by the relation (ii). We define the action of $\dbG_m$ via $\mu_2$ on $(O_{i,j}\oplus O_{i,j_{\perp}})\otimes_{B(k)} K_2$ as follows. Let $l\in\{1,\ldots,s\}$, $i\in\{1,\ldots,m\}$, and $a\in\{1,\ldots,t\}$ be such that $e_a^i\in o_l$. The action of $\dbG_m$ via $\mu_2$ on $K_2e_a^i$ is:

\medskip
{\bf (iv)} via the identical character of $\dbG_m$ if $a$ is the smallest number in $\{1,\ldots,t\}$ such that $e_a^i\in o_l$ and $i$ is the $s_i$-th number in $S_{-,j}$ for some $s_i\in\{1+\sum_{x=1}^{l-1}c^-_{x,j},\ldots,\sum_{x=1}^{l}c^-_{x,j}\}$;

\smallskip
{\bf (v)} via the inverse of the identical character of $\dbG_m$ if $a$ is the smallest number in $\{1,\ldots,t\}$ such that $e_a^i\in o_l$ and $i$ is the $s_i$-th number in $S_{+,j}$ for some $s_i\in\{1+\sum_{x=1}^{l-1}c^+_{x,j},\ldots,\sum_{x=1}^{l}c^+_{x,j}\}$;

\smallskip
{\bf (vi)} trivial otherwise.

\medskip
The action of $\dbG_m$ via $\mu_2$ on a $K_2$-vector subspace $K_2e_a^i$ of $O_{i,j_{\perp}}$ is defined uniquely by the requirement that $\mu_2$ factors through the image of $\scrT_{0,K_2}$ in $\pmb{\GL}_{(O_{i,j}\oplus O_{i,j_{\perp}})\otimes_{B(k)} K_2}$. Due to the property (iii), $\mu_2$ is $\pmb{\text{S}}(K_2)$-conjugate to $\mu_{0,K_2}$.

Also, due to all these relations, the product of the cocharacters of $\scrT_{0,K_2}$ of the orbit under $\Gal(K_2/\dbQ_p)$ of the factorization of $\mu_2$ through $\scrT_{0,K_2}$ induces a cocharacter of $\pmb{\GL}_{(O_{i,j}\oplus O_{i,j_{\perp}})\otimes_{B(k)} K_2}$ that factors through $Z(\pmb{\GL}_{O_{i,j}\otimes_{B(k)} K_2})\times_{K_2} Z(\pmb{\GL}_{O_{i,j_{\perp}}\otimes_{B(k)} K_2})$. Thus choosing $\mu_1$ as in Subsubsection 6.3.2 we get that the product of the cocharacters of $\scrT_{1,K_2}$ which belong to the orbit under $\Gal(K_2/\dbQ_p)$ of $\mu_{1,K_2}$ factors through $L_{\scrG}^0(\phi)_{K_2}\cap (\prod_{i=1}^m\prod_{j=1}^{m_f} Z(\pmb{\GL}_{O_{i,j}\otimes_{B(k)} K_2}))$ and thus through $Z^0(L_{\scrG}^0(\phi)_{K_2})$, cf. (6).

\bigskip\noindent
{\bf 6.4. Case 2.} Until Subsection 6.5 we assume that $\scrG^0$ is of $C_n$ or $D_n$ Dynkin type, that $M_i$ has rank $2n$, and that $\scrC$ is basic. This Case 2 is very much the same as the Case 1 for $\scrC$ basic. We mention only the differences. The first three differences are: 

\medskip
{\bf (i)} we can assume that $O_i=M_i$ (thus $\dbG_m$ acts via $\mu$ trivially on $O_i$ for all $i\in \grM$);
\smallskip
{\bf (ii)} for the $C_n$ Dynkin type the form $b_i$ is alternating and not symmetric;

\smallskip
{\bf (iii)} we have $t=2n$ and $\pmb{\text{S}}=\scrG^{\der}$.  

\medskip
As $M_i=O_i$ let $\varUpsilon$, $\varUpsilon_i$, $e^i_a$, and $\pi_{\star}$ with $\star\in\Gal(K_2/\dbQ_p)$ be as in Subsection 6.3. Let $o_1,\ldots,o_s$ be the orbits of the action of $\Gal(K_2/\dbQ_p)$ on $\varUpsilon$ numbered in such a way that there exists $s_0\in\{0,\ldots,s\}$ with the property that for an element $l\in\{1,\ldots,s\}$ the orbit $o_l$ contains the set $\{e_{2a-1}^1,e_{2a}^1\}$ for some number $a\in\{1,\ldots,n\}$ if and only if $l\Le s_0$. The difference $s-s_0$ is an even number. We can also assume that if $s_1\in\{1,\ldots,{{s-s_0}\over 2}\}$, then the union $o_{s_0+2s_1-1}\cup o_{s_0+2s_1}$ contains the set $\{e_{2a-1}^1,e_{2a}^1\}$ for some number $a\in\{1,\ldots,n\}$. If $l\le s_0$ (resp. $l>s_0$) let $u_l\in\dbN^*$ be such that the number of elements of the set $\tilde o_l:=o_l\cap\varUpsilon_1$ is $2u_l$ (resp. is $u_l$). Lemma 6.3.1 gets replaced by the following one.

\medskip\noindent
{\bf 6.4.1. Lemma.} {\it {\bf (a)} If $v$ is odd, then $u_l$ is even for all $l\in\{s_0+1,\ldots,s\}$.

\smallskip 
{\bf (b)} We have $s_0>0$.}

\medskip
\proof
As $\scrC$ is basic, all Newton polygon slopes of $(M_1,\phi^m)$ are ${v\over 2}$. If $\tau_2\in\Gal(K_2/\dbQ_p)$ is an arbitrary element which extends the automorphism $\grF_{2u}^m$ of $K_{2u}$ (cf. the notations of Definition 2.4 (d)), then we have a direct sum decomposition 
$$(M_1\otimes_{W(k)} K_2,\phi^m\otimes\tau_2)=\oplus_{l=1}^s (\oplus_{e^1_a\in\tilde o_l}K_2e^1_a,\phi^m\otimes\tau_2).\leqno (7)$$
If $v$ is odd, then from the last two sentences we get that each $K_2$-vector space $\oplus_{e^1_a\in\tilde o_l} K_2e^1_a$ is even dimensional and thus $u_l$ is even if $l\in\{s_0+1,\ldots,s\}$. Therefore (a) holds.

The proof of (b) is the same as the proof of Lemma 6.3.1.\endproof 

\medskip
To define $\mu_1$ it is enough to define $\mu_{1,K_2}$. We consider two subcases as follows.

\medskip\noindent
{\bf 6.4.2. Choice of $\mu_1$ for $v$ odd.} If $l>s_0$ we consider a disjoint union $o_l=o_{l,1}\sqcup o_{l,2}$ with both $o_{l,1}$ and $o_{l,2}$ having ${u_l\over 2}$ elements, cf. Lemma 6.4.1 (a). Not to introduce extra notations we will assume that if $l-s_0\in 1+2\dbN$, then the sets $o_{l,1}\cap\varUpsilon_i$ and $o_{l+1,2}\cap\varUpsilon_i$ are perpendicular with respect to $b_i$ for all $i\in\{1,\ldots,m\}$. If $\scrG^0$ is of $C_n$ Dynkin type, then we choose $\mu_{1,K_2}$ such that $\dbG_m$ acts through it:

\medskip
{\bf (i)} trivially on $e^i_a$ if $i\in \grM$;

\smallskip
{\bf (ii)} trivially on $e^i_a\in o_{l,j}$ if $i\in \grN$, $l\ge s_0+1$, and $j\in\{1,2\}$ with $l-s_0-j$ even;

\smallskip
{\bf (iii)} trivially on $e_a^i\in o_l$ if $i\in \grN$, $l\Le s_0$, and $a$ is odd;

\smallskip
{\bf (iv)} via the weight $-1$ on the $K_2$-spans of all other elements of $\varUpsilon$. 

\medskip
For each $l\in\{1,\ldots,s\}$, $\mu_{1,K_2}$ acts non-trivially on half of the elements of $o_l$. Thus the product of the cocharacters of $\scrT_{1,K_2}$ which belong to the $\Gal(K_2/\dbQ_p)$-orbit of $\mu_{1,K_2}$ factors through $\prod_{i=1}^m Z(\pmb{\GL}_{M_i})=Z^0(\scrG_{K_2})$. As $\scrG^0$ is of $C_n$ Dynkin type, $\mu_{1,K_2}$ and $\mu_{K_2}$ are $\scrG(K_2)$-conjugate.

Suppose now $\scrG^0$ is of $D_n$ Dynkin type. In this case we might have to modify the above choice of $\mu_{1,K_2}$ as the cocharacters $\mu_{1,K_2}$ and $\mu_{K_2}$ might not be $\scrG(K_2)$-conjugate (they might be only $\scrG^+(K_2)$-conjugate, where $\scrG^+(K_2)$ is the co-fibre product of the monomorphisms $\scrG^{\der}\hookrightarrow \scrG$ and $\scrG^{\der}\hookrightarrow \scrG^{\der,+}$, the last one being a product of $m$ copies of the standard $\pmb{SO}_{2n}^+\hookrightarrow \pmb{O}_{2n}^+$ monomorphism). As $s_0>0$ (cf. Lemma 6.4.1 (b)), then by re-indexing for suitable elements $i\in\grN$ the unique pair $(e^i_{2a-1},e^i_{2a})$ of elements of $o_1$ on whose $K_2$-span $\mu_{1,K_2}$ acts non-trivially so that the pair becomes the pair $(e^i_{2a},e^i_{2a-1})$, we can assume that $\mu_{1,K_2}$ and $\mu_{K_2}$ are $\scrG(K_2)$-conjugate.

\medskip\noindent
{\bf 6.4.3. Choice of $\mu_1$ for $v$ even.} Let $\grN_0$ be a subset of $\grN$ that has ${v\over 2}$ elements. Let $\grN_1:=\grN\setminus \grN_0$. Let $e^i_a\in o_l$. If $l\le s_0$, then we define the action of $\dbG_m$ via $\mu_{1,K_2}$ on $K_2e^i_a$ as in  Subsubsection 6.4.2. If $l>s_0$, then we define the action of $\dbG_m$ via $\mu_{1,K_2}$ on $K_2e^i_a$ to be via the inverse of the identical character of $\dbG_m$ (resp. trivial) if and only if $i\in \grN_j$, where $j\in\{0,1\}$ is congruent to $l-s_0$ modulo $2$. Thus the number of elements of $o_l$ on whose $K_0$-spans $\dbG_m$ acts via $\mu_{1,K_2}$ as the inverse of the identical character of $\dbG_m$ is $u_l$ (resp. is ${u_l\over 2}$) if $l\Le s_0$ (resp. if $l>s_0$). Therefore the $\Gal(K_2/\dbQ_p)$-orbit of $\mu_{1,K_2}$ factors through $\prod_{i=1}^m Z(\pmb{\GL}_{M_i})=Z^0(\scrG_{K_2})$. 

If $\scrG^0$ is of $C_n$ Dynkin type, then $\mu_{1,K_2}$ and $\mu_{K_2}$ are $\scrG(K_2)$-conjugate. If  $\scrG^0$ is of $D_n$ Dynkin type, then as in the previous subcase we argue that $s_0>0$ and we can modify $\mu_{1,K_2}$ so that it is $\scrG(K_2)$-conjugate to $\mu_{K_2}$.

\medskip\noindent
{\bf 6.4.4. Remark.} If $\scrG^0$ is of $D_4$ Dynkin type, then Subsubsection 6.3.4 extends automatically to the context $\text{rk}_{W(k)}(M_i)=8$ of Case 2. The only difference: we have $M_i=O_i$.

\bigskip\noindent
{\bf 6.5. End of the proof of 4.2.} We recall that Subsection 6.1 achieved the reduction to Cases 1 and 2 of Subsections 6.3 and 6.4. Thus Theorem 4.2 (a) (resp. Theorem 4.2. (b)) follows from Subsubsections 6.3.2, 6.4.2, and 6.4.3  (resp. from Subsubsection 6.3.4).\endproof

\medskip\noindent
{\bf 6.5.1. Remark.} The approach of Subsubsection 6.3.4 extends in many cases to the case when $\scrC$ is basic and $\scrG^0$ is of $A_n$ Dynkin type. However, one has to deal not with only two sets $S_{+,j}$ and $S_{-,j}$ but with $n$ analogue sets and therefore in general it is much harder to show the existence of corresponding numbers $c_{l,u,j}^+$, where $u\in\{1,\ldots,n\}$ is a third index. This is the reason why in Subsection 6.4 we considered only the basic context (and why in Theorem 8.4 below we will rely as well on [Zi1, Thm. 4.4]).

\bigskip\smallskip
\noindent
{\boldsectionfont 7. First applications to abelian varieties}

\bigskip
Pairs of the form  $(\ddag,\lambda_\ddag)$ will denote polarized abelian schemes. By abuse of notations, we also denote by $\lambda_\ddag$ the different forms on the cohomologies (or homologies) of $\ddag$ induced by $\lambda_\ddag$. We now apply Theorems 4.1 and 4.2 to the geometric context of Subsubsections 1.1.1 and 1.2. Applications to Conjecture 1.2.2 (i) and to Subproblem 1.2.3 are included in Subsections 7.3 to 7.5. If $\triangle$ is an algebra, let $\triangle^{\text{opp}}$ be its opposite algebra. 

\bigskip\noindent
{\bf 7.1. Geometric setting.} Until the end we assume that $D$ is the $p$-divisible group of an abelian variety $A$ over $k$, that $\scrC=(M,\phi,\scrG)$ is a Shimura filtered $F$-crystal over $k$ such that axioms 1.2.4 (b.i) and (b.ii) hold, and that there exists a polarization $\lambda_A$ of $A$ whose crystalline realization $\lambda_A:M\times M\to W(k)$ has a $W(k)$-span normalized by $\scrG$. Let $F^1$ and $\mu$ be as in Subsection 2.1. Let $\scrG_{\dbZ_p}$ be as in Subsection 2.5. By performing the operation $\grO_1$ we can assume that the schematic closure $\scrT(\phi)$ of the group $\{\phi^{rm}|m\in\dbZ\}$ in $\scrG_{B(k)}$ is a torus over $B(k)$. This implies that we have an identity $\End(A)=\End(A_{\bar k})$. We identify 
$$\gre_A:=\End(A)^{\text{opp}}\otimes_{\dbZ}\dbQ$$ 
with a $\dbQ$--subalgebra of $\{x\in\End(M[{1\over p}])|\phi(x)=x\}$.

Let $\Pi:\End(M_{\dbQ_p})\to\End(M_{\dbQ_p})$
be the projector on $\Lie(\scrG_{\dbQ_p})$ along the perpendicular of $\Lie(\scrG_{\dbQ_p})$ with respect to the trace form $\grT$ on $\End(M_{\dbQ_p})$, cf. Lemma 2.2 (b). Let $\gre_{\dbQ_p}:=(\gre_A\otimes_{\dbQ}\dbQ_p)\cap\Im(\Pi)$; it is a Lie algebra over $\dbQ_p$. Let $\gre_{\dbQ_p}^{\perp}:=(\gre_A\otimes_{\dbQ} \dbQ_p)\cap\Ker(\Pi)$. As $\Pi$ is fixed by $\phi$, we have a direct sum decomposition
$$\gre_A\otimes_{\dbQ} \dbQ_p=\gre_{\dbQ_p}\oplus\gre_{\dbQ_p}^{\perp}$$
into $\dbQ_p$-vector spaces (or into modules over the Lie algebra $\gre_{\dbQ_p}$). Let $C_A$ be the reductive group over $\dbQ$ of invertible elements of $\gre_A$; thus $\Lie(C_A)$ is the Lie algebra associated to $\gre_A$. A classical theorem of Tate says that $C_{A,B(k)}$ is the centralizer of $\phi^r$ in $\End(M[{1\over p}])$. Thus $\gre_{B(k)}:=\gre_{\dbQ_p}\otimes_{\dbQ_p} B(k)$ is the Lie algebra of the centralizer of $\phi^r$ in $\scrG_{B(k)}$. Let $\gre^1_A$ be a semisimple $\dbQ$--subalgebra of $\gre_A$ which (inside $\End(M[{1\over p}])$) is formed by elements fixed by $\scrG_{B(k)}$ and which is stable under the involution of $\gre_A$ defined naturally by $\lambda_A$.

\bigskip\noindent
{\bf 7.2. Lemma.} {\it Let $\scrT_{1,B(k)}$ be a maximal torus of $\scrG_{B(k)}$ of $\dbQ_p$-endomorphisms of $\scrC$. Then there exists a maximal torus $\scrT^{\text{big}}_{1,B(k)}$ of $\pmb{\GL}_{M[{1\over p}]}$ of $\dbQ_p$-endomorphisms of $(M,\phi,\pmb{\GL}_M)$ and there exists an element $u\in C_A(\dbQ_p)$ such that the following four conditions hold:

\medskip
{\bf (i)} the element $u$ normalizes $M$ (i.e., we have $u(M)=M$) as well as any a priori given $W(k)$-lattice of $M[{1\over p}]$;

\smallskip
{\bf (ii)} the torus $u\scrT^{\text{big}}_{1,B(k)}u^{-1}$ is the extension to $B(k)$ of a maximal torus of $C_A$;

\smallskip
{\bf (iii)} the element $u$ fixes each element of $\gre^1_A$ and normalizes the $\dbQ_p$-span of $\lambda_A$;

\smallskip
{\bf (iv)} we have $\scrT_{1,B(k)}=\scrT^{\text{big}}_{1,B(k)}\cap \scrG_{B(k)}$.}

\medskip
\proof
Let $C^1_A$ be the identity component of the subgroup of $C_A$ that normalizes the $\dbQ$--span of $\lambda_A$ and that centralizes $\gre_A^1$. It is a reductive group over $\dbQ$. Let $\scrT^1_{\dbQ_p}$ be a maximal torus of $C^1_{A,\dbQ_p}$ that contains the $\dbQ_p$-form $\scrT_{1,\dbQ_p}$ of $\scrT_{1,B(k)}$ with respect to $(M[{1\over p}],\phi)$. From [Ha, Lem. 5.5.3] we deduce the existence of an element $u\in C^1_A(\dbQ_p)$ such that the condition (i) holds and $u\scrT^1_{\dbQ_p}u^{-1}$ is the extension to $\dbQ_p$ of a maximal torus $\scrT^{1,u}_{\dbQ}$ of $C^1_A$. Thus the condition (iii) holds. Let $\scrT^{1\text{big}}_{\dbQ}$ be a maximal torus of $C_A$ that contains $\scrT^{1u}_{\dbQ}$. The condition (ii) holds for $\scrT^{\text{big}}_{1,B(k)}:=u^{-1}\scrT^{1\text{big}}_{B(k)}u$. As $\scrT_{1,B(k)}$ is its own centralizer in $\scrG_{B(k)}$ and as the intersection $\scrT^{\text{big}}_{1,B(k)}\cap \scrG_{B(k)}$ contains and centralizes $\scrT_{1,B(k)}$, the condition (iv) holds. \endproof

\medskip
We have the following geometric consequence of Theorem 4.1 (c).

\bigskip\noindent
{\bf 7.3. Corollary.} {\it We assume that $p\Ge 3$ and that $Q\grA$ holds for $\scrC$. 

\medskip
{\bf (a)} Then by performing the operation  $\grO_1$ we can assume there exist an element $\tilde h\in \grP(\scrC)$ and an abelian variety $A(\tilde h)$ over $k$ such that the following two conditions hold:

\medskip
{\bf (i)} the abelian variety $A(\tilde h)$ is $\dbZ[{1\over p}]$-isogenous to $A$ and, under this $\dbZ[{1\over p}]$-isogeny, the Dieudonn\'e module of its $p$-divisible group is $(\tilde h(M),\phi)$ and is a direct sum of isoclinic Dieudonn\'e modules over $k$;

\smallskip
{\bf (ii)} there exists an abelian scheme $A(\tilde h)_{V_3}$ with complex multiplication over a finite, totally ramified discrete valuation ring extension $V_3$ of $W(k)$ which is a ramified lift of $A(\tilde h)$ to $V_3$ with respect to the schematic closure $\tilde \scrG(\tilde h)$ of $\tilde \scrG_{B(k)}$ in $\pmb{\GL}_{\tilde h(M)}$, where $\tilde \scrG$ is a $\pmb{\GL}_M(W(k))$-conjugate of $\scrG$ such that the triple $(M,\phi,\tilde \scrG)$ is a Shimura $F$-crystal over $k$. 

\medskip
{\bf (b)} Then by performing the operation $\grO_1$ we can assume that there exists an element $h\in \grI(\scrC)$ and an abelian variety $A(h)$ over $k$ such that the condition (i) and the following new condition hold:

\smallskip
{\bf (iii)} there exists an abelian scheme $A(h)_{V_3}$ over a finite, totally ramified discrete valuation ring extension $V_3$ of $W(k)$ which lifts $A(h)$ in such a way that the Frobenius endomorphism of $A(h)$ also lifts to it, which is a ramified lift of $A(h)$ to $V_3$ with respect to $\scrG(h)$,  and whose $p$-divisible group $D(h)_{V_3}$ is with complex multiplication.}

\medskip
\proof
We can assume that $\scrC$ is basic (cf. Subsection 2.6) and that there exists an $E$-pair $(\scrT_{1,B(k)},\mu_1)$ of $\scrC$ which is admissible (cf. hypotheses).  We first proof (a). Let $\scrT_{1,B(k)}^{\text{big}}$ and $u$ be as in Lemma 7.2; thus $u$ is fixed by $\phi$. Let $\tilde \scrG$, $\tilde \scrT_{1,B(k)}^{\text{big}}$, $\tilde \scrT_{1,B(k)}$, $\tilde\mu_1$, and $(\tilde t_{\alpha})_{\alpha\in\scrJ}$ be the inner conjugates of $\scrG$, $\scrT_{1,B(k)}^{\text{big}}$, $\scrT_{1,B(k)}$, $\mu_1$, and $(t_{\alpha})_{\alpha\in\scrJ}$ (respectively) through the element $u\in \pmb{\GL}_M(B(k))$. Let $\tilde\scrC:=(M,\phi,\tilde \scrG)$. The Lie algebra $\Lie(\tilde \scrT^{\text{big}}_{1,B(k)})$ is $B(k)$-generated by elements of $\gre_A$ and $(\tilde \scrT_{1,B(k)},\tilde\mu_1)$ is an $E$-pair of $\tilde\scrC$ which is admissible. 

We apply the proof of Theorem 4.1 (c) to $\tilde\scrC$ and $(\tilde \scrT_{1,B(k)},\tilde\mu_1)$ (see Subsection 5.3). Up to the operation $\grO_1$, we deduce the existence of an element $\tilde h\in \grI(\tilde\scrC)$ such that the $p$-divisible group $D(\tilde h)$ over $k$ whose Dieudonn\'e module is $(\tilde h(M),\phi)$ has a lift $D(\tilde h)_{V_3}$ to a finite, totally ramified discrete valuation ring extension $V_3$ of $W(k)$ with respect to $\tilde\scrG(\tilde h)$ such that each endomorphism of $D(\tilde h)$ whose crystalline realization is an element of $\Lie(\tilde \scrT^{\text{big}}_{1,B(k)})$ fixed by $\phi$ lifts to an endomorphism of $D(\tilde h)_{V_3}$ (cf. also Subsubsection 5.3.3). The $B(k)$-span of $\lambda_A$ is normalized by $\tilde G_{B(k)}$, cf. property 7.2 (iii). Let $A(\tilde h)$ be the abelian variety over $k$ defined by the condition (i). Let $A(\tilde h)_{V_3}$ be the abelian scheme over $V_3$ defined by $D(\tilde h)_{V_3}$, cf. Serre--Tate deformation theory. The fact that $A(\tilde h)_{V_3}$ is indeed an abelian scheme (and not only a formal abelian scheme over $\Spf V_3$) is implied by the fact that we are in a polarized context. As $D(\tilde h)_{V_3}$ is a ramified lift of $D(\tilde h)$ to $V_3$ with respect to $\tilde \scrG(\tilde h)$, the abelian scheme $A(\tilde h)_{V_3}$ is a ramified lift of $A(\tilde h)$ to $V_3$ with respect to $\tilde \scrG(\tilde h)$. As $\tilde\scrC$ is basic, the part of the condition (i) on Dieudonn\'e modules over $k$ holds (cf. Fact 2.6.1 (a)). Thus the condition (i) holds. The fact that $A(\tilde h)_{V_3}$ has complex multiplication follows from the condition 7.2 (ii). Thus (ii) holds as well. From the very definition of $\tilde\scrC$ we get that $\grI(\tilde\scrC)\subset\grP(\scrC)$; thus indeed $\tilde h\in\grP(\scrC)$. This proves (a). 

Part (b) follows from the proof of Theorem 4.1 (c) applied in the context of the $E$-pair $(\scrT_{1,B(k)},\mu_1)$ of $\scrC$ which is admissible. Up to the operation $\grO_1$, we get the existence of an element $h\in \grI(\scrC)$, of a finite, totally ramified discrete valuation ring extension $V_3$ of $W(k)$, and of a $p$-divisible group $D(h)_{V_3}$ over $V_3$ that is a ramified lift of $D(h)$ to $V_3$ with respect to $(h(M),\phi,\scrG(h))$ and that has the property that each endomorphism of $D(h)$ whose crystalline realization is an element of $\Lie(\scrT_{1,B(k)})\cap\End(M)$ fixed by $\phi$ lifts to an endomorphism of $D(h)_{V_3}$ (cf. Subsubsections 5.3.1 to 5.3.3). Let $A(h)$ and $A(h)_{V_3}$ be obtained as above but for $h$ instead of $\tilde h$. The $p$-divisible group $D(h)_{V_3}$ is with complex multiplication, cf. Subsubsection 5.3.3. As $\phi^r\in \scrT_{1,B(k)}(B(k))$ leaves invariant $h(M)$, the Frobenius endomorphism of $A(h)$ lifts to $A(h)_{V_3}$ (cf. also Subsubsection 5.3.3 and Serre--Tate deformation theory). As above we argue that the condition (i) holds. As $D(h)_{V_3}$ is a ramified lift of $D(h)$ to $V_3$ with respect to $(h(M),\phi,\scrG(h))$, the abelian scheme $A(h)_{V_3}$ is a ramified lift of $A(h)$ to $V_3$ with respect to $\scrG(h)$.\endproof 

\bigskip\noindent
{\bf 7.4. Proof of Theorem 1.2.5.} To prove Theorem 1.2.5, we can assume that $L_{\scrG}^0(\phi)$ is a Levi subgroup scheme of $P_{\scrG}^+(\phi)$, that $\mu:\dbG_m\to\scrG$ factors through $L^0_{\scrG}(\phi)$ (cf. Subsection 2.6), and that $TT\grR$ holds for $(M,\phi,L_{\scrG}^0(\phi))$ (cf. Theorem 4.2 (b)). Thus $TT\grU$ holds for $(M,\phi,L_{\scrG}^0(\phi))$ (cf. Theorem 4.1 (b)) and therefore for $\scrC$. Thus for $p\ge 3$ we get that all hypotheses of Corollary 7.3 hold. Therefore Theorem 1.2.5 (b) follows from Corollary 7.3 (b).

We now check that Theorem 1.2.5 (a) also follows from Corollary 7.3 (b) and Subsubsection 5.3.4. We consider a torus $\scrT_{1,B(k)}$ of $\scrG_{B(k)}$ of $\dbQ_p$-endomorphisms of $\scrC$ which splits over $B(k)$. Then in the proof of Corollary 7.3 (b), up to the operations $\grO_1$ and $\grO_2$, we can assume that we are in an unramified context (i.e., Subsubsection 5.3.4 implies that we can take $V_3=W(k)$ and we can allow $p=2$ as well) and thus Theorem 1.2.5 (a) holds.\endproof

\medskip
Corollary 7.3 (b) is our partial solution to Subproblem 1.2.3, cf. also Theorem 1.2.5 (b). We also have the following important variant of Corollary 7.3 (a). 

\bigskip\noindent
{\bf 7.5. Variant.} We assume that $p\Ge 3$ and that $Q\grA$ holds for $\scrC$. We also assume that $\Lie(\scrT_{1,B(k)})$ is $B(k)$-generated by elements of $\gre_A$; thus we can take $u$ to be $1_M$ and $\tilde\scrC=\scrC$. Thus the condition 7.3 (ii) holds with $\tilde h\in \grI(\scrC)=\grI(\tilde\scrC)$ and thus with $\tilde \scrG(\tilde h)=\scrG(\tilde h)$. This solves Conjecture 1.2.2 (i) under all our assumptions.

\bigskip\smallskip
\noindent
{\boldsectionfont 8. The context of standard Hodge situations}
\bigskip

If $(G,\scrX)$ is a Shimura pair, let $E(G,\scrX)$ be the subfield of $\dbC$ which is its reflex field, let $(G^{\ad},\scrX^{\ad})$ be its adjoint Shimura pair, and let $\Sh(G,\scrX)$ be the canonical model over $E(G,\scrX)$ of $\Sh(G,\scrX)$ (see [De1,2], [Mi3, Subsects. 1.1 to 1.8], and [Va1, Subsects. 2.2 to 2.8]). Let $\Sh(G,\scrX)/\scrK$ be the quotient of $\Sh(G,\scrX)$ by a compact subgroup $\scrK$ of $G(\dbA_f)$. See [Va1, Subsect. 2.4] for injective maps between Shimura pairs. For general properties of Shimura varieties of PEL type we refer to [Zi1], [LR], [Ko2, Ch. 5], [Mi3, p. 161], and [RaZ] (but in [Mi3, p. 161] one has to request that [De2, axiom 2.1.1.3] holds). The injective maps into Siegel modular varieties that define Shimura varieties of PEL type as used in these references, will be referred as {\it PEL type embeddings}. Let $O_{(w)}$ be the localization of the ring of integers of a number field with respect to a finite prime $w$ of it. 

In Subsection 8.1 we mainly introduce notations and a setting. Different ``properties" pertaining to the setting of Subsection 8.1 are introduced in Subsection 8.2. Theorems 8.3 and 8.4 are results that are needed in our approach to prove the Langlands--Rapoport conjecture for Shimura varieties of Hodge type (see already [Va10]). 

\bigskip\noindent
{\bf 8.1. Standard Hodge situation.} We recall part of the setting of [Va7, Sect. 5] pertaining to good reduction cases of Shimura varieties of Hodge type. Let
$$f:\Sh(G,\scrX)\hookrightarrow \Sh(\pmb{\text{GSp}}(W,\psi),\scrS)$$ 
be an injective map of Shimura pairs. Here the Shimura pair $(\pmb{\text{GSp}}(W,\psi),\scrS)$ defines a Siegel modular variety, cf. Subsubsection 1.2.1. We consider a
$\dbZ$-lattice $L$ of $W$ such that $\psi$ induces a perfect form 
$\psi:L\times L\to\dbZ$. 
Let $L_{(p)}:=L\otimes\dbZ_{(p)}$. We will assume that:

\medskip\noindent
{\it the schematic closure $G_{\dbZ_{(p)}}$ of $G$ in $\pmb{\text{GSp}}(L_{(p)},\psi)$ is a reductive group scheme over $\dbZ_{(p)}$.} 

\medskip\noindent 
It is easy to see that the group scheme $G^0_{\dbZ_{(p)}}:=G_{\dbZ_{(p)}}\cap \pmb{\text{Sp}}(L_{(p)},\psi)$ is reductive, cf. [Va7, Subsect. 5.1, Formula (12)]. Let $\scrK_p:=\pmb{\text{GSp}}(L_{(p)},\psi)(\dbZ_p)$; it is a hyperspecial subgroup of $\pmb{\text{GSp}}(W,\psi)_{\dbQ_p}(\dbQ_p)$. As $G_{\dbZ_{(p)}}$ is a reductive group scheme over $\dbZ_{(p)}$, the intersection $H:=G_{\dbQ_p}(\dbQ_p)\cap \scrK_p$ is a hyperspecial subgroup of $G_{\dbQ_p}(\dbQ_p)$. 

Let $v$ be a prime of the reflex field $E(G,\scrX)$ that divides $p$; it is unramified over $p$ (cf. [Mi4, Prop. 4.6 and Cor. 4.7]). Let $k(v)$ be the residue field of $v$.
Let $d:={{\dim_\dbQ(W)}\over 2}\in\dbN^*$. Let $\dbA_f$ (resp. $\dbA_f^{(p)}$) be the $\dbQ$--algebra of finite ad\`eles (resp. of finite ad\`eles with the $p$-component omitted). We have $\dbA_f=\dbA_f^{(p)}\times\dbQ_p$.

For integral canonical models of (suitable quotients of) Shimura varieties we refer to [Va1, Subsubsects. 3.2.3 to 3.2.6].  It is well known that the $\Spec \dbZ_{(p)}$-scheme $\scrM$ that parameterizing isomorphism classes
of principally polarized abelian schemes of relative dimension $d$ over $\Spec \dbZ_{(p)}$-schemes which have compatible level-$N$ symplectic similitude structures for all $N\in\dbN^*\setminus p\dbN^*$, together
with the natural action of $\pmb{\text{GSp}}(W,\psi)(\dbA_f^{(p)})$ on it, is an integral canonical model of
$\Sh(\pmb{\text{GSp}}(W,\psi),\scrS)/\scrK_p$ (for instance, see [De1, Thm. 4.21] and [Va1, Ex. 3.2.9 and Subsect. 4.1]). These structures and this action are defined naturally via the $\dbZ$-lattice $L$ of $W$ (see [Va1, Subsect. 4.1]). The set $\Sh(G,\scrX)/H(\dbC)$ is naturally identified with $G_{\dbZ_{(p)}}(\dbZ_{(p)})\backslash (\scrX\times G(\dbA_f^{(p)}))$, cf. [Mi4, Prop. 4.11 and Cor. 4.12]. From this and its analogue for $\Sh(\pmb{\text{GSp}}(W,\psi),\scrS)/\scrK_p(\dbC)$, one gets that the functorial morphism $\Sh(G,\scrX)/H\to\scrM_{E(G,\scrX)}=\Sh(\pmb{\text{GSp}}(W,\psi),\scrS)_{E(G,\scrX)}/\scrK_p$ is a closed embedding (to be compared with [Va1, Rm. 3.2.14]). 

Let $\scrN$ be the normalization of the schematic closure of $\Sh(G,\scrX)/H$ in $\scrM_{O_{(v)}}$. Let $(\scrA,\lambda_{\scrA})$ be the pull-back to $\scrN$ of the universal principally polarized abelian scheme over $\scrM$. Let
$(v_{\alpha})_{\alpha\in\scrJ}$ be a family of tensors of $\scrT(L_{(p)}^{\vee})$ such that $G$ is the subgroup of $\pmb{\GL}_W$ that fixes $v_{\alpha}$ for all $\alpha\in\scrJ$ (cf. [De3, Prop. 3.1 c)]). The choice of $L$ and $(v_\alpha)_{\alpha\in\scrJ}$ allows a moduli interpretation of $\Sh(G,\scrX)$ (see [De1,2], [Mi4], and [Va1, Subsect. 4.1 and Lem. 4.1.3]). For instance, $\Sh(G,\scrX)/H(\dbC)$ is the set of isomorphism classes of complex principally polarized abelian varieties of dimension $d$ that carry a family of Hodge cycles indexed by the set $\scrJ$, that have compatible level-$N$ symplectic similitude structures for all $N\in\dbN^*\setminus p\dbN^*$, and that satisfy certain axioms (see [Va1, Subsect. 4.1]). This moduli interpretation endows naturally the abelian scheme $\scrA_{E(G,\scrX)}$ with a family $(w_{\alpha}^{\scrA})_{\alpha\in\scrJ}$ of Hodge cycles (the Betti realizations of pull-backs of $w_{\alpha}^{\scrA}$ via $\dbC$-valued points of $\scrN_{E(G,\scrX)}$ correspond naturally to $v_{\alpha}$). Let $H_0$ be a compact, open subgroup of $G(\dbA_f^{(p)})$ that has the following three properties:

\medskip 
{\bf (a)} there exists $N_0\in\dbN^*$ such that $(N_0,p)=1$, $N_0\Ge 3$, and we have an inclusion
$$H_0\times H\subseteq \scrK(N_0):=\{g\in \pmb{\text{GSp}}(L,\psi)(\hat\dbZ)|g\equiv 1_{L\otimes_{\dbZ} \hat\dbZ} \, {\text{modulo}}\, N_0\hat\dbZ\};$$

\smallskip
{\bf (b)} the triple $\scrR:=(\scrA,\lambda_{\scrA},(w_{\alpha}^{\scrA})_{\alpha\in\scrJ})$ is the pull-back of an analogue triple $\scrR(H_0)=(\scrA_{H_0},\lambda_{\scrA_{H_0}},(w_{\alpha}^{\scrA_{H_0}})_{\alpha\in\scrJ})$ over $\scrN/H_0$, where $(\scrA_{H_0},\lambda_{\scrA_{H_0}})$ is the pull-back via the natural morphism $\scrN/H_0\to\scrM/\scrK^p(N_0)$ of the universal principally polarized abelian scheme over $\scrM/\scrK^p(N_0)$ (here $\scrK^p(N_0)\leqslant \pmb{\text{GSp}}(W,\psi)(\dbA_f^{(p)})$ is such that we have $\scrK(N_0)=\scrK^p(N_0)\times \scrK_p$);

\smallskip
{\bf (c)} the scheme $\scrN$ is a pro-\'etale cover of $\scrN/H_0$ (Serre Lemma implies that $\scrK^p(N_0)$ acts freely on $\scrM$ and thus $H_0$ acts freely on $\scrN$; to be compared with [Va1, Prop. 3.4.1]).

\medskip\noindent
{\bf 8.1.1. Some notations.} Let $k=\dbF_q$ be a finite field that contains $k(v)$. We consider a $W(k)$-morphism $z:\Spec W(k)\to\scrN/H_0$. Let 
$$(A_{W(k)},\lambda_{A_{W(k)}},(w_{\alpha})_{\alpha\in\scrJ})=z^*(\scrR(H_0)).$$
Let $y:\Spec k\to\scrN_{k(v)}/H_0$ and $(A,\lambda_A)$ be the special fibres of $z$ and $(A_{W(k)},\lambda_{A_{W(k)}})$ (respectively). Let $(M,\phi,\lambda_A)$ be the principally quasi-polarized Dieudonn\'e module of the principally quasi-polarized $p$-divisible group of $(A,\lambda_A)$. Let $F^1$ be the Hodge filtration of $M$ defined by $A_{W(k)}$. For $\alpha\in\scrJ$ let $t_{\alpha}\in\scrT(M[{1\over p}])$ be the de Rham component of the Hodge cycle $w_{\alpha}$ on $A_{W(k)}$; it belongs to the $F^0$-filtration of $\scrT(M[{1\over p}])$ defined by $F^1[{1\over p}]$. Let $\tilde \scrG$ be the schematic closure in $\pmb{\GL}_M$ of the subgroup of $\pmb{\GL}_{M[{1\over p}]}$ that fixes $t_{\alpha}$ for all $\alpha\in\scrJ$. We will assume that the triple $(f,L,v)$ is a {\it standard Hodge situation} in the sense of [Va7, Def. 5.1.2]. Thus the following two properties hold:

\medskip
{\bf (a)} the $\Spec O_{(v)}$-scheme $\scrN/H_0$ is smooth;

\smallskip
{\bf (b)} for each point $z\in\scrN/H_0(W(k))$, $\tilde \scrG$ is a reductive subgroup scheme of $\pmb{\GL}_M$.

\medskip
In [Va1,3] it is proved that (a) holds for $p\ge 5$ while in [Ki2] it is claimed that (a) holds for $p\ge 3$. Also it is known that (a) implies (b) if $p\ge 3$, cf. Remark 8.2.1 (c) below.

Let $\scrG:=G_{W(k)}$ and $\scrG^0:=G^0_{W(k)}$. By performing the operation $\grO_1$ we can assume that $\tilde \scrG$ is isomorphic to $\scrG$. By multiplying each $v_{\alpha}$ by a fixed integral power of $p$ we can assume that for all points $z\in\scrN/H_0(W(k))$ we have $t_{\alpha}\in\scrT(M)$ for all $\alpha\in\scrJ$. To match the notations with those of Sections 1 to 7, we will identify (non-canonically) $\tilde \scrG=\scrG$.

Each tensor $t_{\alpha}$ is fixed under the natural action of $\phi$ on $\scrT(M[{1\over p}])$ (cf. [Va7, Cor. 5.1.7]) and thus it is fixed by the inverse of the canonical split cocharacter of $(M,F^1,\phi)$ defined in [Wi, p. 512] (cf. the functorial aspects of [Wi, p. 513]). The resulting cocharacter $\mu:\dbG_m\to\scrG$ produces a direct sum decomposition $M=F^1\oplus F^0$ ($\dbG_m$ acts through $\mu$ on $F^1$ via weight $-1$ and it fixes $F^0$). We conclude that $\scrC:=(M,\phi,\tilde \scrG)$ is a Shimura $F$-crystal over $k$ having $\mu$ as a Hodge cocharacter. The axioms 1.2.4 (b.i) holds for $\scrC$, cf. the existence of the $t_{\alpha}$'s. The fact that the axioms 1.2.4 (b.ii) is implied by [De2, axiom 2.1.1.3] (cf. the existence of the Hodge cocharacters of $G_{\dbC}$ associated naturally to points $x\in X$).

The triple $\scrC$ depends only on $y$ and not on $z$ (cf. [Va7, paragraph before Subsubsect. 5.1.8] and thus we call it the {\it Shimura $F$-crystal attached} to the point $y\in\scrN_{k(v)}/H_0(k)$. Let $z_{\infty}:\Spec W(\bar k)\to\scrN$ be such that the resulting $W(\bar k)$-valued point of $\scrN/H_0$ factors through $z$. We refer to $\scrC\otimes\bar k$ as the Shimura $F$-crystal attached to the special fibre $y_{\infty}:\Spec \bar k\to\scrN_{k(v)}$ of $z_{\infty}$.  We also refer to $y_{\infty}$ (resp. $z_{\infty}$) as an {\it infinite lift} of $y$ (resp. of $z$). We also refer to $F^1$ as the lift of $\scrC$ defined by the point $z\in\scrN/H_0(W(k))$ that lifts $y\in\scrN_{k(v)}/H_0(k)$. 

For another point $y_j\in\scrN/H_0(k)$ let $(A_j,\lambda_{A_j})$, $\scrC_j=(M_j,\phi_j,\scrG_j)$, $y_{j,\infty}$, and $(t_{j,\alpha})_{\alpha\in\scrJ}$ be the analogues of $(A,\lambda_A)$, $\scrC$, $y_{\infty}$, and $(t_{\alpha})_{\alpha\in\scrJ}$ obtained by replacing $y$ with $y_j$. 

\medskip\noindent
{\bf 8.1.2. PEL type embeddings.} Let $C_{\dbQ}:=C_{\pmb{\GL}_W}(G)$. Let $G_1$ be the identity component of $C_{1,\dbQ}:=\pmb{\text{GSp}}(W,\psi)\cap C_{\pmb{\GL}_W}(C_{\dbQ})$; it contains $G$. Let $\scrX_1$ be the $G_1(\dbR)$-conjugacy class of homomorphisms $\Res_{\dbC/\dbR}\dbG_m\to G_{1,\dbR}$ that contain the composites of elements of $\scrX$ with the monomorphism $G_{\dbR}\hookrightarrow G_{1,\dbR}$. We get a PEL type embedding $f_1:(G_1,\scrX_1)\hookrightarrow (\pmb{\text{GSp}}(W,\psi),\scrS)$ through which $f$ factors. We call it the {\it PEL-envelope} of $f$, cf. [Va1, Rm. 4.3.12]. 

Let $G_{2,\dbZ_{(p)}}:=C_{\pmb{\text{GSp}}(L_{(p)},\psi)}(Z^0(G_{\dbZ_{(p)}}))$; it is a reductive group scheme over $\dbZ_{(p)}$ (cf. [DG Vol. III, Exp. XIX, Subsect. 2.8]). Let $G_2$ be the generic fibre of $G_{2,\dbZ_{(p)}}$; it contains $G_1$ and moreover we have $Z^0(G_1)=Z^0(G_2)$. As above we get an injective map $f_2:(G_2,\scrX_2)\hookrightarrow (\pmb{\text{GSp}}(W,\psi),\scrS)$ through which both $f$ and $f_1$ factor naturally. 

Let $i\in\{1,2\}$. Let $H_i:=G_{i,\dbQ_p}(\dbQ_p)\cap \scrK_p$. Let $v_i$ be the prime of the subfield $E(G_i,\scrX_i)$ of $E(G,\scrX)$ divided by $v$ and let $k(v_i)$ be its residue field. Let $\scrN_i$ be the normalization of the schematic closure of $\Sh(G_i,\scrX_i)/H_i$ in $\scrM_{O_{(v_i)}}$. If $H_{0i}:=G_i(\dbA_f^{(p)})\cap K^p(N_0)$, then $H_0\leqslant H_{01}\leqslant H_{02}$ and the morphism $\scrN_i\to\scrN_i/H_{0i}$ is a pro-\'etale cover.

The injective map $f_i$ is a PEL type embedding. The triple $(f_2,L_{(p)},v_2)$ is a standard Hodge situation (this well known fact follows from either [Zi1, Subsect. 3.5] or [LR]). Let $(\scrA_i,\lambda_{\scrA_i})$ be the pull-back to $\scrN_i$ of the universal principally polarized abelian scheme over $\scrM$. Let $\scrG_i$ be the integral, closed subgroup scheme of $\pmb{\GL}_M$ which has the analogue meaning of $\scrG=\tilde \scrG$ but obtained working with the $k$-valued point $y_i$ of $\scrN_{i,k(v_i)}/H_{0i}$ defined by $y$. The group scheme $\scrG_2$ is reductive.${}^1$ $\vfootnote{1} {If either $p>3$ or $p=2$ and $C_{1,\dbQ}$ is connected, then it is easy to see that Theorem 2.5.2 (a) implies that $\scrG_1$ is also a reductive group scheme (see also [LR] and [Ko2]; for $p>2$, cf. also Remark 8.2.1 (b) below).}$ 

We use the notations of Subsection 7.1 and (by performing the operation $\grO_1$) until the end we will assume that $\scrT(\phi)$ is a torus of $\scrG_{B(k)}$. Let $\gre_{1A}:=\gre_A\cap \Lie(\scrG_{1,B(k)})$. Identifying the opposite of the $\dbQ$--algebra that defines $\Lie(C_{\dbQ})$ with a semisimple $\dbQ$--subalgebra of $\gre_A$, we get that $\gre_{1A}$ is the maximal $\dbQ$--vector subspace of $\gre_A$ that centralizes $\Lie(C_{\dbQ})$ and that leaves invariant the $\dbQ$--span of $\lambda_A$. We can assume that $\grz_1:=\Lie(Z^0(G_{1,\dbZ_{(p)}}))$, when viewed as a set, is included in $\{v_{\alpha}|\alpha\in\scrJ\}$. 

\medskip\noindent
{\bf 8.1.3. Rational stratification.} Let $\grS_{\text{rat}}$ be the rational stratification of $\scrN_{k(v)}$ defined in [Va7, Subsect. 5.3]. We recall that if $y_1\in\scrN/H_0(k)$, then $y_{1,\infty}$ and $y_{\infty}$ are $\bar k$-valued points of the same (reduced) stratum of $\grS_{\text{rat}}$ if and only if there exists an isomorphism $(M_1\otimes_{W(k)} B(\bar k),\phi_1\otimes\sigma_{\bar k})\arrowsim (M\otimes_{W(k)} B(\bar k),\phi\otimes\sigma_{\bar k})$ that takes $t_{1,\alpha}$ to $t_{\alpha}$ for all $\alpha\in\scrJ$. The number of strata of $\grS_{\text{rat}}$ is finite, cf. [Va7, Rm. 5.3.2 (c)]. Let $\grs_0$ be the $G(\dbA_f^{(p)})$-invariant, reduced, closed subscheme of $\scrN_{k(v)}$ defined by the following property: the point $y_{\infty}$ factors through $\grs_0$ if and only if $\scrC$ is basic. Obviously $\grs_0$ is a union of strata of $\scrN_{k(v)}$.

\bigskip\noindent
{\bf 8.2. Properties.} Let $h\in \grI(\scrC)\subseteq \scrG(B(k))$. Let $A(h)$ be as in Subsubsection 1.1.1. We denote by $\lambda_{A(h)}$ the principal polarization of $A(h)$ defined naturally by $\lambda_A$; its crystalline realization is a rational multiple of $\lambda_A$. Let 
$$y(h):\Spec k\to\scrM_{k(v)}/\scrK^p(N_0)$$ 
be the morphism defined by $(A(h),\lambda_{A(h)})$ and its level-$N_0$ symplectic similitude structure induced naturally from the one of $(A,\lambda_A)$ defined by the point $y\in\scrN_{k(v)}/H_0(k)$. Let $y(h)_{\infty}:\Spec \bar k\to\scrM_{k(v)}$ be an infinite lift of $y(h)$. 

\medskip
{\bf (a)} For $p\Ge 3$ (resp. $p=2$) we say the {\it isogeny property} holds for the point $y\in\scrN_{k(v)}/H_0(k)$ if for every element $h\in \grI(\scrC)$, the (resp. up to the operation $\grO_1$ the) morphism $y(h)$ factors through $\scrN_{k(v)}/H_0$ and there exists a point $z(h)\in\scrN/H_0(W(k))$ which lifts this factorization (denoted in the same way) $y(h):\Spec k\to\scrN_{k(v)}/H_0$ and for which $t_{\alpha}$ is the de Rham realization of $z(h)^*(w_{\alpha}^{\scrA_{H_0}})$ for all $\alpha\in\scrJ$. We say the isogeny property holds for $(f,L,v)$ if the isogeny property holds for each point of $\scrN_{k(v)}/H_0$ with values in a finite field. 

\smallskip
{\bf (b)} We say the {\it weak isogeny property} holds for $(f,L,v)$ if $\grs_0$ is the only stratum of $\grS_{\text{rat}}$ that has a closed connected component.

\smallskip
{\bf (c)} We say the {\it Milne conjecture} holds for $(f,L,v)$ if for each point $y\in\scrN_{k(v)}/H_0(k)$ there exists a symplectic isomorphism $(M,\lambda_A)\arrowsim (L^{\vee}\otimes_{\dbZ} W(k),\psi^{\vee})$ that takes $t_{\alpha}$ to $v_{\alpha}$ for all $\alpha\in\scrJ$. Here $\psi^{\vee}$ is the alternating form on $L^{\vee}$ defined naturally by $\psi$.

\medskip\noindent
{\bf 8.2.1. Remarks.} {\bf (a)} There exists a standard but non-canonical identification $(L^{\vee}\otimes_{\dbZ} \dbZ_p,\psi^{\vee})=(H^1_{\acute{et}}(A_{\overline{B(k)}},\dbZ_p),\lambda_A)$ under which $v_{\alpha}$ is mapped to the $p$-component of the \'etale component of $w_{\alpha}$ for all $\alpha\in\scrJ$, cf. [Va7, Subsubsect. 5.1.5].

\smallskip
{\bf (b)} If $p\ge 3$, then the Milne conjecture holds for $(f,L,v)$ (cf. (a) and either [Va8, Cor. 1.4] or [Ki2, Cor. (1.4.3)]) and thus the property 8.1.1 (a) implies property 8.1.1 (b). 

\smallskip
{\bf (c)} The isogeny property for $(f,L,v)$ was announced in [Va1, Subsubsect. 1.7.1]. See [Va10] for a proof of it and of the weak isogeny property in most cases of interest. 

\smallskip
{\bf (d)} For each $\beta\in\dbG_m(W(k))$ there exists an element $\mu(\beta^{-1})\in \scrG(W(k))$ that acts on the $W(k)$-span of $\lambda_A$ via multiplication with $\beta$. Thus if there exists an isomorphism $(M,(t_{\alpha})_{\alpha\in\scrJ})\arrowsim (L^{\vee}\otimes_{\dbZ} W(k),(t_{\alpha})_{\alpha\in\scrJ})$, then there exists also an isomorphism of the form $(M,(t_{\alpha})_{\alpha\in\scrJ},\lambda_A)\arrowsim (L^{\vee}\otimes_{\dbZ} W(k),(t_{\alpha})_{\alpha\in\scrJ},\psi^{\vee})$. Moreover, as $\scrG^0$ is smooth and has a connected special fibre, such isomorphisms $(M,(t_{\alpha})_{\alpha\in\scrJ},\lambda_A)\arrowsim (L^{\vee}\otimes_{\dbZ} W(k),(t_{\alpha})_{\alpha\in\scrJ},\psi^{\vee})$ exist if and only if they exist in the flat topology of $W(k)$.

\medskip\noindent
{\bf 8.2.2. Lemma.} {\it We assume that $p\ge 3$ and that the isogeny property holds for the point $y\in \scrN_{k(v)}/H_0(k)$. We also assume that the condition $Q\grU$ (of Definition 2.4 (g)) holds for $\scrC$ (for instance this, holds if each simple factor of $(G^{\ad},\scrX^{\ad})$ is of $B_n$ or $D_n^{\dbR}$ type, cf. Theorem 1.2.5 (a)). Let $\pi$ be the Frobenius endomorphism of $A$. Then, up to operations $\grO_1$ and $\grO_2$ (i.e., up to a passage to a finite field extension of $k$ and up to the replacement of $y$ by $y(h)\in\scrN_{k(v)}/H_0(k)$ for some $h\in\grI(\scrC)$), there exists an element $\pi_{\dbQ}\in G(\dbQ)$ for which the following two properties hold:

\medskip
{\bf (i)} for each prime $l\in\dbN^*$ different from $p$, it is $G(\dbQ_l)$-conjugate to the $l$-adic realization of $\pi$ under the natural identification $H^1_{\acute et}(A,\dbQ_l)=L^{\vee}\otimes_{\dbZ} \dbQ_l$ of $\dbQ_l$-vector spaces;

\smallskip
{\bf (ii)} there exists an isomorphism $M[{1\over p}]\arrowsim L^{\vee}\otimes B(k)$ which takes $t_{\alpha}$  to $v_{\alpha}$ for all $\alpha\in\scrJ$ and which maps $\phi^r$ to $\pi_{\dbQ}\in G(B(k))$.}

\medskip
\proof
Based on the proof of Theorem 1.2.5 (a) (see Subsection 7.4), we can assume that there exists a lift $F^1$ of $\scrC$ such that the Frobenius endomorphism of $A$ lifts to an endomorphism $\grF$ of the abelian scheme $A_{W(k)}$ over $W(k)$ whose Hodge filtration is $F^1$. Let $z\in\scrN/H_0(W(k))$ be such that $A_{W(k)}=z^*(\scrA)$, cf. [Va9, Sect. 6.4]. We fix an $O_{(v)}$-monomorphism $W(k)\hookrightarrow \dbC$ and via it we consider the complex abelian variety $A_{\dbC}$. The Mumford--Tate group $G_{A_{\dbC}}$ of $A_{\dbC}$ is a reductive subgroup of $G$ and we have an injective map $(G_{A_{\dbC}},\scrX_{A_{\dbC}})\hookrightarrow (G,\scrX)$ of Shimura pairs (to be compared with [Va4, Subsect. 1.3]). The center of $G_{A_{\dbC}}$ has a $\dbQ$--valued point $\pi_{\dbQ}$ which is the Betti realization of $\grF_{B(k)}$; thus we also have $\pi_{\dbQ}\in G(\dbQ)$. The fact that the element $\pi_{\dbQ}\in G(\dbQ)$ satisfies the properties (i) and (ii) is well known (see [Pi], [Va4], etc.).\endproof

\bigskip\noindent
{\bf 8.3. Theorem.} {\it If $p=2$, we assume that the Milne conjecture holds for $(f,L,v)$. 

\medskip
{\bf (a)} Then $\grs_0$ is a stratum of $\grS_{\text{rat}}$ which is closed. 

\smallskip
{\bf (b)} We also assume that $G^{\der}$ is simply connected and that $y$ factors through $\grs_0/H_0$. Let $y_0\in\grs_0/H_0(k)$. Let $y_{0,\infty}:\Spec \bar k\to\scrN_{k(v)}$ be an infinite lift of $y_0$. Then up to the operation $\grO_1$, there exist elements $t\in G_2(\dbA_f^{(p)})$ and $h\in \grI(\scrC)$ such that we have an identity $y_{0,\infty}t=y(h)_{\infty}$ of $\bar k$-valued points of $\scrM_{k(v)}$.}

\medskip
\proof
The connected components of $\scrN$ are permuted transitively by $G(\dbA_f^{(p)})$, cf. [Va1, Lem. 3.3.2]. Thus to prove the proposition, we can assume that $y\in\grs_0/H_0(k)$ and that both  $y_{\infty}$ and $y_{0,\infty}$ factor through the special fibre of the same connected component $\scrN^0$ of $\scrN$. Let $\pi$ (resp. $\pi_0$) be the Frobenius endomorphism of $A$ (resp. of $A_0$). See [Ch, Subsect. 3.a] for the Frobenius tori $\scrT_{\pi}$ and $\scrT_{\pi_0}$ over $\dbQ$ of $\pi$ and $\pi_0$ (respectively). We recall that if $\dbQ[\pi]$ is the \'etale $\dbQ$--subalgebra of $\gre_A$ generated by $\pi$, then $\scrT_{\pi}$ is the smallest subtorus of $\Res_{\dbQ[\pi]/\dbQ} \dbG_m$ which has $\pi$ as a $\dbQ$--valued point. The crystalline realization of $\pi$ is $\phi^r\in \scrG(B(k))$ and therefore we have an identity $\scrT(\phi)=\scrT_{\pi, B(k)}$. Each element $b\in \grz_1$ defines naturally a $\dbZ_{(p)}$-endomorphism of any pull-back of $\scrA$, $\scrA_1$, or $\scrA_2$, to be denoted also by $b$. Thus we view $Z^0(G)$, $Z^0(G_1)$, $Z^0(G_2)$, and $\scrT_{\pi}$ as tori of $C_A$. 

We prove (a). As $y\in\grs_0/H_0(k)$, the Newton quasi-cocharacter of $(M,\phi,\scrG)$ factors through $Z^0(\scrG_{B(k)})$ (see [Va7, Cor. 2.3.2]) and thus it can be identified with a quasi-cocharacter $\mu_0$ of $Z^0(\scrG_{B(k)})$. This quasi-cocharacter depends only on the $\Gal(\dbQ_p)$-orbit of the composite $\mu^{\ab}:\dbG_m\to\scrG^{\ab}$ of the Hodge cocharacter $\mu$ of $\scrC$ with the canonical epimorphism $\scrG\twoheadrightarrow \scrG^{\ab}$. Moreover $\mu^{\ab}$ is uniquely attached to $\scrX$, cf. [Va7, Subsubsects. 5.1.1 and 5.1.9]. We conclude that, as the notation suggests, $\mu_0$ does not depend on the point $y\in\grs_0/H_0(k)$.

The torus $\scrT_{\pi}$ is the smallest torus of $C_A$ with the property that $\mu_0$ is a quasi-cocharacter of $\scrT_{\pi,B(k)}$, cf. Serre's result of [Pi, Prop. 3.5]. Thus $\scrT_{\pi}$ is naturally identified with a subtorus of $Z^0(G)$ uniquely determined by $\scrX$. Applying this also to $y_0$ we get that $\scrT_{\pi}=\scrT_{\pi_0}$. Thus $\pi_0\in\dbQ[\pi]$ is such that its image in each number field factor $S_0$ of $\dbQ[\pi]$ is non-trivial. Therefore from [Ta2] we get that the images of $\pi$ and $\pi_0$ in $S_0$ are both Weil $q$-integers. Thus the image of ${\pi\over \pi_0}$ in each $S_0$ is a root of unity and therefore by performing the operation $\grO_1$ we can assume that $\pi=\pi_0\in \scrT_{\pi}(\dbQ)=\scrT_{\pi_0}(\dbQ)\leqslant Z^0(G)(\dbQ)$. Let $\iota:A\to A_0$ be the $\dbQ$--isogeny defined by this equality, cf. [Ta2]. Let $(M_0[{1\over p}],\phi)\arrowsim (M[{1\over p}],\phi_0)$ be the isomorphism defined by $\iota$; we will view it as a natural identification. 

Let $\scrG^{0\prime}_{\dbQ_p}$ be the $\dbQ_p$-form of $\scrG^0_{B(k)}$ with respect to $(M[{1\over p}],\phi)$. We have $\scrG^{0\prime\ab}_{\dbQ_p}=\scrG^{0\ab}_{\dbQ_p}$ and therefore let $\scrG^{0\prime\ab}_{\dbZ_p}:=\scrG^{0\ab}_{\dbZ_p}$.

As the Milne conjecture holds for $(f,L,v)$, there exists an element $j\in \pmb{\GL}_M(B(k))$ such that $j(M)=M_0$ and $j$ takes $\lambda_A$ to $\lambda_{A_0}$ and takes $t_{\alpha}$ to $t_{0,\alpha}$ for all $\alpha\in\scrJ$. Thus $j$ commutes with $\phi^r=\phi^r_0\in Z^0(\scrG_{B(k)})(B(k))$. We can also assume that $j$ takes a Hodge cocharacter of $\scrC$ to a Hodge cocharacter of $\scrC_0$ (to be compared with [Va7, Lem. 5.1.9]). Thus we can identify $j^{-1}\phi_0 j=j^{-1}\phi j=g\phi$, where $g\in \scrG^0(W(k))$. From [Va7, Prop. 2.7 and Subsect. 4.7] we deduce the existence of an element $h\in \scrG^0(B(\bar k))$ such that we have $g(\phi\otimes\sigma_{\bar k})=h^{-1}(\phi\otimes\sigma_{\bar k})h$. In other words, there exists an isomorphism $(M_0\otimes_{W(k)} B(\bar k),\phi_0\otimes\sigma_{\bar k},(t_{0,\alpha})_{\alpha\in\scrJ},\lambda_{A_0})\arrowsim (h(M)\otimes_{W(k)} B(\bar k),\phi\otimes\sigma_{\bar k},(t_{\alpha})_{\alpha\in\scrJ},\lambda_A)$ defined by $hj^{-1}$. Thus $\grs_0$ is a stratum of $\grS_{\text{rat}}$ (cf. the definition of the rational stratification [Va7, Subsect. 5.3]); it is closed by its very definition. Thus (a) holds.

We prove (b); thus $G^{\der}$ is simply connected. As $j$ commutes with $\phi^r$ we get that $\phi^r=(j^{-1}\phi j)^r=(g\phi)^r$. But with respect to the $\dbQ_p$-form $\scrG^{0\prime}_{\dbQ_p}$ of $\scrG^0_{B(k)}$, the action of $\phi$ on $\scrG^0(B(k))=\scrG^{0\prime}_{\dbQ_p}(B(k))$ is the action of $\sigma$ on $\scrG^{0\prime}_{\dbQ_p}(B(k))=\scrG^0(B(k))$ and thus we have $\phi^r=(g\phi)^r=\prod_{i=0}^{r-1}\sigma^i(g)\phi^r$. Therefore $\prod_{i=0}^{r-1}\sigma^i(g)=1_M$ and thus $g$ defines naturally a class $\gamma_g\in H^1(\Gal(B(k)/\dbQ_p),\scrG^{0\prime}_{\dbQ_p})$. The image of $\gamma_g$ in $H^1(\Gal(B(k)/\dbQ_p),\scrG^{0\prime\ab}_{\dbQ_p})$ factors through $H^1(\Gal(W(k)/\dbZ_p),\scrG^{0\prime\ab}_{\dbZ_p})$. As $H^1(\Gal(W(k)/\dbZ_p),\scrG^{0\prime\ab}_{\dbZ_p})=0$ (cf. Lang theorem) and as the homomorphism $\scrG^0(W(k))\to \scrG^{0\ab}(W(k))$ is surjective, we can assume that we have $g\in \scrG^{\der}(W(k))$. Thus $\gamma_g$ is the image of some class $\gamma_g^{\der}\in H^1(\Gal(B(k)/\dbQ_p),\scrG^{\prime\der}_{\dbQ_p})$. As $G^{\der}$ is simply connected, the class $\gamma_g^{\der}$ is trivial (cf. [Kn, Thm. 1]) i.e., we can assume that we have $h\in \scrG^{\der}(B(k))$ and $g=h^{-1}\sigma(h)=h^{-1}\phi(h)$. Therefore $g\phi=h^{-1}\phi h$ and thus we have $h\in \grI(\scrC)$ and $j^{-1}\phi j=h^{-1}\phi h$. Let $\tilde h:=hj^{-1}\in \pmb{\GL}_M(B(k))$; it is fixed by $\phi$ and thus it is a $\dbQ_p$-valued point of $C_A$.

Let $Z_A$ be the reductive subgroup of $C_A$ that fixes $\lambda_A$ and $\grz_1[{1\over p}]$. We now check that we can assume that $\iota$ takes $b$ to $b$ for all $b\in \grz_1[{1\over p}]$ and takes $\lambda_A$ to $\lambda_{A_0}$. Let $\gamma_0\in H^1(\dbQ,Z_A)$ be the class that ``measures'' the existence of such a choice of $\iota$. Let $l$ be a rational prime. We check that the image of $\gamma_0$ in $H^1(\dbQ_l,Z_{A,\dbQ_l})$ is the trivial class. If $l=p$ (resp. if $l\neq p$), then this is so due to the previous paragraph (resp. due to the existence of all level-$l^m$ symplectic similitude structures of $\scrA$ with $m\in\dbN^*$ and on the fact that $\pi=\pi_0\in T_{\pi}(\dbQ)=T_{\pi_0}(\dbQ)\leqslant Z^0(G)(\dbQ)$). The triples $(A,\lambda_A,\grz_1)$ and $(A_0,\lambda_{A_0},\grz_1)$ lift to characteristic 0. But all pull-backs of $(\scrA,\lambda_{\scrA},\grz_1)$ via complex valued points of $\scrN^0$ are $\dbR$-isogenous (as each connected component of $\scrX$ is a $G^0(\dbR)$-conjugacy class). Thus $(A,\lambda_A,\grz_1)$ and $(A_0,\lambda_{A_0},\grz_1)$ are $\dbR$-isogenous i.e., the image of $\gamma_0$ in $H^1(\dbR,Z_{A,\dbR})$ is the trivial class. 

The group $Z_{A,\dbC}$ is isomorphic to the centralizer of a torus of $\pmb{\text{Sp}}(W,\psi)_{\dbC}$ in $\pmb{\text{Sp}}(W,\psi)_{\dbC}$. Thus it is the product of some $\pmb{\GL}_{n_1}$ groups with either a trivial group or with a $\pmb{\text{Sp}}_{2n_2}$ group (the ranks $n_1$ and $n_2$ do depend on the factors of such a product). Therefore we have a product decomposition $Z_A=Z_1\times_{\dbQ} Z_2$, where:

\medskip
{\bf (i)}  there exists a semisimple $\dbQ$--algebra $Z_{11}$ with involution $\vartheta_{11}$ such that $Z_1$ is the group scheme of invertible elements of $Z_{11}$ fixed by $*$;

\smallskip
{\bf (ii)}  $Z_2$ is either trivial or a simple connected semisimple group of $C_n$ Dynkin type. 

\medskip\noindent
The pair $(Z_{11},\vartheta_{11})$ is a product of semisimple $\dbQ$--algebras endowed with involutions which are either trivial or of second type. Thus $Z_1$ is a product of Weil restrictions of reductive groups whose derived groups are forms of $\pmb{\text{SL}}_n$ groups ($n\in\dbN^*$) and whose abelianizations are of rank 1. This implies that the Hasse principle holds for $Z_1$ (even if some $n$'s are even). It is well known that the Hasse principle holds for $Z_2$. We conclude that:

\medskip
{\bf (iii)} the Hasse principle holds for $Z_A$ and therefore the class $\gamma_0$ is trivial (cf. previous paragraph).

\medskip
It is well known that $Z_1(\dbQ)$ is dense in $Z_1(\dbQ_p)$. As $Z_{A,B(k)}$ is $C_{\pmb{\text{Sp}}(M[{1\over p}],\lambda_A)}(\scrT(\phi))$ and as $Z^0(\scrG)$ splits over a finite unramified extension of $W(k)$, the group $Z_{2,B(\bar k)}$ is split. Thus $Z_2(\dbQ)$ is dense in $Z_2(\dbQ_p)$, cf. [Mi4, Lem. 4.10]. Thus we get:

\medskip
{\bf (iv)} the group $Z_A(\dbQ)$ is dense in $Z_A(\dbQ_p)$.

\medskip
Due to the property (iii), we can assume that $j\in \scrG_2^0(B(k))$. Thus $\tilde h\in\scrG_2^0(B(k))$ is a $\dbZ_p$-isomorphism between the principally quasi-polarized Dieudonn\'e modules with endomorphisms associated to $(A_0,\lambda_{A_0},\grz_1)$ and $(A(h),\lambda_A,\grz_1)$. Let $s\in\dbN^*$. A theorem of Tate says that $\Hom_k(A(h),A_0)\otimes_{\dbZ} \dbZ_p$ is the set of $\dbZ_p$-homomorphisms between the Dieudonn\'e modules of the $p$-divisible groups of $A_0$ and $A(h)$ (see [Ta2, p. 99]; the passage from $\dbQ_p$ coefficients to $\dbZ_p$ coefficients is trivial). Based on this and the property (iv) we get that there exists a $\dbZ_{(p)}$-isomorphism $\tilde h_s$ between $(A(h),\lambda_A,\grz_1)$ and $(A_0,\lambda_{A_0},\grz_1)$ whose crystalline realization is congruent modulo $p^s$ to $\tilde h$.

Due to existence of the $\dbZ_{(p)}$-isomorphism $\tilde h_s$, there exists $t\in \pmb{\text{GSp}}(W,\psi)(\dbA_f^{(p)})$ such that we have an identity $y_{0,\infty}t=y(h)_{\infty}$ of $\bar k$-valued points of $\scrM_{k(v)}$  (cf. also [Mi2, Sect. 3]). The fact that we can take $t\in G_2^0(\dbA_f^{(p)})$ is checked easily by considering the level-$l^m$ symplectic similitude structures of $(y(h)_{\infty})^*(\scrA,\lambda_{\scrA})$ and $y_{0,\infty}^*(\scrA,\lambda_{\scrA})$  (here $l$ is a prime different from $p$ while $m\in\dbN^*$ is arbitrary).\endproof

\bigskip\noindent
{\bf 8.4. Theorem.} {\it  Let $n\in\dbN^*\setminus\{1\}$. We assume that all simple factors of $(G^{\ad},\scrX^{\ad})$ are of either $C_n$ or $D_n^{\dbH}$ type, that all simple factors of $(G_1^{\ad},\scrX_1^{\ad})$ are of $A_{2n-1}$ type, and that $C_{\dbQ}$ is indecomposable (equivalently, and that $C_{\dbQ}$ is the group scheme of invertible elements of a simple $\dbQ$--algebra). We also assume that the monomorphism $G^{\der}_{\dbC}\hookrightarrow G_{1,\dbC}^{\der}$ is a product of monomorphisms of one of the following two forms: $\pmb{\text{Sp}}_{2n}\hookrightarrow \pmb{\text{SL}}_{2n}$ and $\pmb{\text{SO}}_{2n}\hookrightarrow \pmb{\text{SL}}_{2n}$. 

\medskip
{\bf (a)} Then $TT\grA$ holds for $\scrC$.

\smallskip
{\bf (b)} Then up to the operations $\grO_1$ and $\grO_2$ (i.e., up to passage to a finite field extension of $k$ and up to the replacement of $y$ by $y(h)\in\scrM_{k(v)}/H_0$ for some $h\in\grI(\scrC)$), the following two properties hold:

\medskip\noindent
{\bf (b.i)} There exists an abelian scheme $A_{W(k)}$ over $W(k)$ which lifts $A$, to which the Frobenius endomorphism of $A$ lifts, and whose Hodge filtration $F^1$ is a lift of $\scrC$.  

\smallskip\noindent
{\bf (b.ii)} We assume that $p>2$. For each maximal torus $\scrT_{1,B(k)}$ of $\scrG_{B(k)}$ of $\dbQ_p$-endomorphisms of $\scrC$, there exist a discrete valuation ring $V$ which is a finite, totally ramified extension of $W(k)$ and an abelian scheme $A_V$ over $V$ which is a lift of $A$ with respect to $\scrG$, whose $p$-divisible group $D_V$ has complex multiplication, and for which each element of $\Lie(\scrT_{1,B(k)})$ fixed by $\phi$ is the crystalline realization of a $\dbQ$-endomorphism of $D_V$ (i.e., of a $\dbQ_p$-endomorphism of $A_V$).}

\medskip
\proof
We prove (a). We have $C_{1,\dbQ}=G_1$ (i.e., the group $C_{1,\dbQ}$ is connected), cf. the fact that $G_1^{\ad}$ is of $A_{2n-1}$ Dynkin type. The schematic closure $C_{\dbZ_{(p)}}$ of $C_{\dbQ}$ in $\pmb{\GL}_{L_{(p)}}$ is a reductive group scheme, cf. proof of Theorem 2.5.2 (a) applied over $\dbZ_p$. Thus the schematic closure of $G_1$ in $\pmb{\GL}_{L_{(p)}}$ is a reductive group scheme, the triple $(f_1,L,v_1)$ is a standard Hodge situation, and the $\Spec O_{(v_1)}$-scheme $\scrN_1/H_{01}$ is smooth (see [LR] and [Ko2]). Let $y_1$ and $\scrG_1$ be as in Subsubsection 8.1.2. Thus $\scrC_1:=(M,\phi,\scrG_1)$ is the Shimura $F$-crystal attached to the point $y_1\in\scrN_{1,k(v_1)}/H_{01}(k)$. 

Let $\scrT_{1,B(k)}$ be an arbitrary  maximal torus of $\scrG_{B(k)}$ of $\dbQ_p$-endomorphisms of $\scrC$. Let $K_1$ be as in Definition 2.4 (c). We need to find a cocharacter $\mu_1:\dbG_m\to \scrT_{1,K_1}$ such that, with $F^1_{K_2}$ is as in Definition 2.4 (h), the filtered module $(M,\phi,F^1_{K_2})$ is admissible. The existence of $\mu_1$ depends only on $\scrC$ up to the operations $\grO_1$ and $\grO_2$ and thus from now we will forget about $f$ and $\scrN$ and we will only keep in mind that the quadruple $(M,\phi,\Lie(C_{\dbZ_{(p)}}),\lambda_A)$ is the crystalline realization of a principally polarized abelian variety endowed with $\dbZ_{(p)}$-endomorphisms $(A,\lambda_A,\Lie(C_{\dbZ_{(p)}}))$ over $k$ and that our hypotheses get translated into properties of the group schemes $\scrG$, $\scrG_1$, etc. By performing the operation $\grO_1$, we can assume that $\scrG$ and $\scrG_1$ are split. Also, by enlarging $G$ we can assume that $Z^0(G)=Z^0(G_1)$. The main property required below is the following one (cf. also hypotheses):

\medskip
{\bf (i)} We have $Z^0(\scrG)=Z^0(\scrG_1)$ and the monomorphism $\scrG^{\der}\hookrightarrow \scrG_1^{\der}$ is a product of monomorphisms of one of the following two forms: $\pmb{\text{Sp}}_{2n}\hookrightarrow \pmb{\text{SL}}_{2n}$ and $\pmb{\text{SO}}^{\text{split}}_{2n}\hookrightarrow \pmb{\text{SL}}_{2n}$. 

\medskip
By performing the operation $\grO_2$, we can assume that $L_{\scrG}^0(\phi)$ is a reductive group scheme (cf. Subsection 2.6) and based on Fact 2.6.1 (a) it is easy to see that $L_{\scrG_1}^0(\phi)$ is a centralizer in $\scrG_1$ of a split rank torus of $Z^0(L_{\scrG}^0(\phi))$ and thus it is also a reductive group scheme. Thus we can work with $(M,\phi,L_{\scrG}^0(\phi))$ and $(M,\phi,L_{\scrG_1}^0(\phi))$ instead of $\scrC$ and $\scrC_1$. From Fact 2.6.1 (a) and the property (i) we get the existence of a direct sum decomposition 
$$(M,\phi)=\oplus_{j\in J} (M_j,\phi)$$ 
into Dieudonn\'e modules over $k$ that have only one Newton polygon slope such that for each $j\in J$ the following two properties hold:

\medskip
{\bf (ii)} the adjoint of the image $L^0(j)$ of $L_{\scrG}^0(\phi)$ in $\pmb{\GL}_{M_j}$ via the projection $\prod_{j\in J} \pmb{\GL}_{M_j}\twoheadrightarrow \pmb{\GL}_{M_j}$, has all simple factors of the same Lie type 
$$\theta(j)\in\{C_m,D_m,A_m|m\in\{1,\ldots,n\}\}\cup\{0\}$$
(here we have $\theta(j)=0$ if and only if this adjoint is trivial i.e., if and only if the image $L^0(j)$ of $L_{\scrG}^0(\phi)$ in $\pmb{\GL}_{M_j}$ is a torus);

\smallskip
{\bf (iii)} the image $L^0_1(j)$ of $L_{\scrG_1}^0(\phi)$ in $\pmb{\GL}_{M_j}$ via the same projection is either $L^0(j)$ or its adjoint has all simple factors of the same Lie type $A_{2m-1}$ and $\theta(j)\in\{C_m,D_m\}$. 

\medskip
The centralizer $\scrT_{1,B(k)}^\prime$ of $\scrT_{1,B(k)}$ in $\scrG_{1,B(k)}$ is a maximal torus of $\scrG_{1,B(k)}$ of $\dbQ_p$-endomorphisms of $\scrC_1$, cf. property (i) and the fact that $\scrT(\phi)$ is a torus (cf. Subsubsection 8.1.2). Let $\scrT_{1,1B(k)}$ be the centralizer of $\scrT_{1,B(k)}$ in $C_{\pmb{\GL}_{M[{1\over p}]}}(C_{B(k)})$; it is a torus over $B(k)$ which contains $\scrT_{1,B(k)}^\prime$. We use the notations of Definition 2.4 (c), an upper index $\prime$ being used in connection to $\scrT_{1,B(k)}^\prime$. We apply [Ha, Lem. 5.5.3] to the reductive subgroup of $C_A$ that normalizes the $\dbQ$--span of $\lambda_A$ and that fixes each $\dbQ$--endomorphism of $A$ defined by an element of $\Lie(C_{\dbQ})$. Thus up to a replacement of $(M,\phi,L_{\scrG}^0(\phi))$ by $(M,\phi,hL_{\scrG}^0(\phi)h^{-1})$, where $h\in \scrG_1(W(k))$ commutes with $\phi$, we can assume that $\Lie(\scrT_{1,B(k)}^\prime)$ is $B(k)$-generated by elements of $\gre_A$. From [Zi1, Thm. 4.4] we get that up to the operations $\grO_1$ and $\grO_2$, the isogeny class of principally polarized abelian varieties endowed with endomorphisms whose crystalline realization is $(M[{1\over p}],\phi,(\Lie(\scrT_{1,1B(k)})+\Lie(C_{B(k)}))\cap\End(A),\lambda_A)$ has a lift to the ring of fractions of a finite field extension $K_3^\prime$ of $K_2^\prime$. Thus the Hodge cocharacter of this lift, when viewed in the crystalline context, is the extension to $K_3^\prime$ of a cocharacter $\mu_1^\prime:\dbG_m\to \scrT_{1,K_1^\prime}^\prime$ such the $E$-pair $(\scrT_{1,B(k)}^\prime,\mu_1^\prime)$ of $\scrC_1$ is admissible. 

The key point is that the product $\grQ$ of the simple factors of $L_{\scrG}^0(\phi)^{\ad}$ which are of some $A_m$ Lie type, $m\Ge 2$, is the same as the similar product for $L_{\scrG_1}^0(\phi)^{\ad}$ (cf. property (i)). We choose a cocharacter $\mu_1:\dbG_m\to \scrT_{1,K_1}$ such that the following three properties hold (to be compared with Subsection 6.1):

\medskip
{\bf (iv.a)}  The cocharacters of $\grQ_{\overline{B(k)}}$ defined by $\mu_1$ and $\mu_1^\prime$ coincide.

\smallskip
{\bf (iv.b)} The cocharacter defined by $\mu_1$ of the product of the simple factors of $L_{\scrG}^0(\phi)^{\ad}_{K_1}$ which are not subgroups of $\grQ_{K_2}$ is constructed based on Subsection 6.4.

\smallskip
{\bf (iv.c)} A $L_{\scrG}^0(K_2)$-conjugate of $\mu_{1,K_2}$ is the extension to $K_2$ of a Hodge cocharacter $\mu_0:\dbG_m\to L_{\scrG}^0$ of $(M,\phi,L_{\scrG}^0(\phi))$.

\medskip
The admissible filtered modules over $K_2$ are stable under direct sums. Thus to check that $(\scrT_{1,B(k)},\mu_1)$ is admissible (i.e., to end the proof of (a)), we can work with a fixed $j_0\in J$ and we have to show that the filtered module $(M_{j_0}[{1\over p}],\phi,F^1_{K_2}\cap (M_{j_0}\otimes_{W(k)} K_2))$ over $K_2$ is admissible. Let $J_A:=\{j\in J|\theta(j)=A_m,\,\, m\in\dbN^*\setminus\{1\}\}$. If $j_0\in J_A$, then $(M_{j_0}[{1\over p}],\phi,F^1_{K_2}\cap M_{j_0}\otimes_{W(k)} K_2)$ is admissible as $(\scrT_{1,B(k)}^\prime,\mu_2)$ is admissible. If $j_0\in J\setminus J_A$ and $\theta(j_0)\neq 0$, then the fact that $(M_{j_0}[{1\over p}],\phi,F^1_{K_2}\cap M_{j_0}\otimes_{W(k)} K_2)$ is admissible follows from Theorems 4.1 (b) and 4.2 (b) (cf. Subsection 6.4). If $j_0\in J\setminus J_A$ and $\theta(j_0)=0$ and if $F^1_0$ is the maximal direct summand of $M$ on which $\dbG_m$ acts via $\mu_0$ through the weight $-1$, then from properties (i) and (iv.c) we get that the filtered module $(M_{j_0}[{1\over p}],\phi,F^1_{K_2}\cap M_{j_0}\otimes_{W(k)} K_2)$ is the extension to $K_2$ of the filtered Dieudonn\'e module $(M_{j_0},F^1_0,\phi)$ and thus it is admissible. Thus (a) holds.

Part (b.i) (resp. (b.ii)) follows from (a) and the proof of this is the same as the proof of Theorem 1.2.5 (a) (resp. Theorem 1.2.5 (b)) in Subsection 7.4.\endproof

\medskip\noindent
{\bf 8.4.1. Remark.} Referring to Theorem 8.4 (b), if the isogeny property holds for the point $y\in\scrN_{k(v)}/H_0(k)$, then it is easy to check that we have $A_{W(k)}=z^*(\scrA)$ and $A_V=z_V^*(\scrA)$ for suitable lifts $z\in \scrN/H_0(W(k))$ and $z_V\in \scrN/H_0(V)$ of $y\in\scrN_{k(v)}/H_0(k)$. 

\medskip\noindent
{\bf Acknowledgment.} We would like to thank University of Arizona, MPI--Bonn, Binghamton University, IAS--Princeton, and TIFR--Mumbai for good working conditions. We would also like to thank J. S. Milne for his encouragements to write this paper. This research was partially supported by the NSF grant DMS \#0900967.

\bigskip
\noindent
\references{37}
{\nspace{

\bigskip

\Ref[Bo]
A. Borel,
\sl Linear algebraic groups. Second edition,
\rm Grad. Texts in Math., Vol. {\bf 126}, Springer-Verlag, New York, 1991.

\Ref[Bou1]
N. Bourbaki,
\sl Lie groups and Lie algebras. Chapters 4--6,
\rm Elements of Mathematics (Berlin), Springer-Verlag, Berlin, 2002.

\Ref[Bou2]
N. Bourbaki,
\sl Lie groups and Lie algebras. Chapters 7--9,
\rm Elements of Mathematics (Berlin), Springer-Verlag, Berlin, 2005.

\Ref[Br]
C. Breuil,
\sl Groupes $p$-divisible, groupes finis et modules filtr\'es,
\rm Ann. of Math. {\bf 152} (2000), no. 2,  489--549.

\Ref[BBM]
P. Berthelot, L. Breen, and W. Messing, 
\sl Th\'eorie de Dieudonn\'e cristalline. II, 
\rm Lecture Notes in Math., Vol. {\bf 930}, Springer-Verlag, Berlin, 1982.

\Ref[BLR]
S. Bosch, W. L\"utkebohmert, and M. Raynaud,
\sl N\'eron models,
\rm Ergebnisse der Mathematik und ihrer Grenzgebiete (3), Vol. {\bf 21}, Springer-Verlag, Berlin, 1990.

\Ref[BM]  P. Berthelot and W. Messing, 
\sl Th\'eorie de Dieudonn\'e cristalline. III, 
\rm The Grothendieck Festschrift, Vol. I,  173--247, Progr. Math., Vol. {\bf 86}, Birkh\"auser Boston, Boston, MA, 1990.

\Ref[BO]
P. Berthelot and A. Ogus,
\sl F-crystals and de Rham cohomology. I, 
\rm Invent. Math. {\bf 72} (1983), no. 2,  159--199.

\Ref[BT]
F. Bruhat and J. Tits, 
\sl Groupes r\'eductifs sur un corps local. II. Sch\'emas en groupes. Existence d'une donn\'ee radicielle valu\'ee,
\rm Inst. Hautes \'Etudes Sci. Publ. Math., Vol. {\bf 60},  5--184, 1984.

\Ref[Ch]
W. C. Chi,
\sl $l$-adic and $\lambda$-adic representations associated to abelian varieties defined over number fields,
\rm Amer. J. Math. {\bf 114} (1992), no. 2,  315--353.

\Ref[CF] P. Colmez and J.-M. Fontaine, 
\sl Construction des repr\'esentations $p$-adiques semi-stables, 
\rm Invent. Math. {\bf 140} (2000), no. 1,  1--43.

\Ref[dJ] J. de Jong, 
\sl Crystalline Dieudonn\'e module theory via formal and rigid geometry, 
\rm Inst. Hautes \'Etudes Sci. Publ. Math., Vol. {\bf 82},  5--96, 1995.

\Ref[De1]
P. Deligne,
\sl Travaux de Shimura,
\rm S\'eminaire  Bourbaki, 23\`eme ann\'ee (1970/71), Exp. No. 389, Lecture Notes in Math., Vol. {\bf 244},  123--165, Springer-Verlag, Berlin, 1971.

\Ref[De2]
P. Deligne,
\sl Vari\'et\'es de Shimura: interpr\'etation modulaire, et
techniques de construction de mod\`eles canoniques,
\rm Automorphic forms, representations and $L$-functions (Oregon State Univ., Corvallis, OR, 1977), Part 2,   247--289, Proc. Sympos. Pure Math., {\bf 33}, Amer. Math. Soc., Providence, RI, 1979.

\Ref[De3]
P. Deligne,
\sl Hodge cycles on abelian varieties,
\rm Hodge cycles, motives, and Shimura varieties, Lecture Notes in Math., Vol.  {\bf 900},  9--100, Springer-Verlag, Berlin-New York, 1982.

\Ref[DG]
M. Demazure, A. Grothendieck, et al.,
\sl Sch\'emas en groupes. Vols. {\bf II--III},
\rm S\'eminaire de G\'eom\'etrie Alg\'ebrique du Bois Marie 1962/64 (SGA 3), Lecture Notes in Math., Vols. {\bf 152--153}, Springer-Verlag, Berlin-New York, 1970. 

\Ref[Fa]
G. Faltings,
\sl Integral crystalline cohomology over very ramified valuation rings, 
\rm J. Amer. Math. Soc. {\bf 12} (1999), no. 1,  117--144.

\Ref[Fo]
J.-M. Fontaine,
\sl Repr\'esentations $p$-adiques semi-stables,
\rm P\'eriodes $p$-adiques (Bures-sur-Yvette, 1988), Ast\'erisque {\bf 223},  113--184, Soc. Math. de France, Paris, 1994.

\Ref[FR]
J.-M. Fontaine and M. Rapoport,
\sl Existence de filtrations admissibles sur des isocristaux,
\rm Bull. Soc. Math. France {\bf 133}  (2005),  no. 1,  73--86. 

\Ref[Ha]
G. Harder,
\sl \"Uber die Galoiskohomologie halbeinfacher Matrizengruppen II,
\rm  Math. Z. {\bf 92} (1966),  396--415

\Ref[He]
S. Helgason,
\sl Differential geometry, Lie groups, and symmetric spaces,
\rm Pure and Applied Mathematics, Vol. {\bf 80}, Academic Press, Inc. [Harcourt Brace Jovanovich, Publishers], New York-London, 1978.

\Ref[Hu]
J. Humphreys, 
\sl Introduction to Lie algebras and representation theory. Second printing, revised,
\rm Grad. Texts in Math., Vol. {\bf 9}, Springer-Verlag, New York-Berlin, 1978.

\Ref[Ii1] 
Y. Ihara,
\sl The congruence monodromy problems,
\rm J. Math. Soc. Japan {\bf 20} (1968),  107--121.

\Ref[Ii2] 
Y. Ihara,
\sl On congruence monodromy problems, 
\rm Lecture Notes, Vols. {\bf 1} and {\bf 2}, Department of Mathematics, Univ. Tokyo, 1968 and 1969.

\Ref[Ii3] 
Y. Ihara,
\sl Some fundamental groups in the arithmetic of algebraic curves over finite fields,
\rm Proc. Nat. Acad. Sci. U.S.A. {\bf 72} (1975), no. 9,  3281--3284. 

\Ref[Ja]
J. C. Jantzen,
\sl Representations of algebraic groups. Second edition,
\rm  Math. Surv. and Monographs, Vol. {\bf 107}, Amer. Math. Soc., Providence, RI, 2003.

\Ref[Ka] N. Katz, 
\sl Slope filtration of $F$-crystals, 
\rm Journ\'ees de G\'eom\'etrie Alg\'ebrique de Rennes (Rennes 1978), Vol. I,   113--163, Ast\'erisque {\bf 63}, Soc. Math. France, Paris, 1979.

\Ref[Ki1] 
M. Kisin,
\sl Crystalline representations and $F$-crystals, 
\rm Algebraic geometry and number theory, 459--496, Progr. Math., Vol. {\bf 253}, Birkh\"auser Boston, Boston, MA, 2006.

\Ref[Ki2] M. Kisin,
\sl Integral canonical models of Shimura varieties of abelian type, 
\rm J. Amer. Math. Soc. {\bf 23} (2010), no. 4,  967--1012.

\Ref[Kn]
M. Kneser,
\sl Galois-Kohomologie halbeinfacher algebraischer Gruppen \"uber $p$-adischen K\"orpern. II., 
\rm Math. Zeit. {\bf 89} (1965),  250--272. 

\Ref[Ko1]
R. E. Kottwitz, 
\sl Isocrystals with additional structure,
\rm Compositio Math. {\bf 56} (1985), no. 2,  201--220.

\Ref[Ko2]
R. E. Kottwitz, 
\sl Points on some Shimura varieties over finite fields, 
\rm J. Amer. Math. Soc. {\bf 5} (1992), no. 2,  373--444.

\Ref[Lan]
R. Langlands,
\sl Some contemporary problems with origin in the Jugendtraum,
\rm Mathematical developments arising from Hilbert problems (Northern Illinois Univ., De Kalb, IL, 1974),  401--418, Proc. Sympos. Pure Math., Vol. {\bf 28}, Amer. Math. Soc., Providence, RI, 1976.

\Ref[Lau1]
E. Lau, 
\sl Frames and finite group schemes over complete regular
local rings, 
\rm Documenta Math. {\bf 15} (2010), 545--569.
 
\Ref[Lau2]
E. Lau, 
\sl A relation between Dieudonn\'e displays and crystalline Dieudonn\'e theory, \rm http://arxiv.org/abs/1006.2720.

\Ref[Le]
J. D. Lewis,
\sl A survey of the Hodge conjecture. Second edition, 
\rm CRM Monograph Series, Vol. {\bf 10}, Amer. Math. Soc., Providence, RI, 1999.

\Ref[LR]
R. Langlands and M. Rapoport,
\sl Shimuravariet\"aten und Gerben, 
\rm J. Reine Angew. Math. {\bf 378} (1987),  113--220.

\Ref[Ma] 
J. I. Manin, 
\sl The theory of formal commutative groups in finite characteristic, 
\rm Russian Math. Surv. {\bf 18} (1963), no. 6,  1--83. 

\Ref[Me]
W. Messing,
\sl The crystals associated to Barsotti-Tate groups:
with applicactions to abelian schemes,
\rm Lecture Notes in Math., Vol. {\bf 264}, Springer-Verlag, Berlin-New York, 1972.

\Ref[Mi1]
J. S. Milne,
\sl Points on Shimura varieties mod $p$,
\rm Automorphic forms, representations and $L$-functions (Oregon State Univ., Corvallis, OR, 1977), Part 2,   165--184, Proc. Sympos. Pure Math., Vol. {\bf 33}, Amer. Math. Soc., Providence, RI, 1979.

\Ref[Mi2]
J. S. Milne,
\sl The conjecture of Langlands and Rapoport for Siegel modular varieties,
\rm Bull. Amer. Math. Soc. (N.S.) {\bf 24} (1991),  no. 2,  335--341.

\Ref[Mi3]
J. S. Milne,
\sl The points on a Shimura variety modulo a prime of good
reduction,
\rm The Zeta functions of Picard modular surfaces,  153--255, Univ. Montr\'eal, Montreal, Quebec, 1992.

\Ref[Mi4]
J. S. Milne,
\sl Shimura varieties and motives,
\rm Motives (Seattle, WA, 1991), Part 2,  447--523, Proc. Sympos. Pure Math., Vol. {\bf 55}, Amer. Math. Soc., Providence, RI, 1994.

\Ref[Mi5]
J. S. Milne,
\sl Points on Shimura varieties over finite fields: the conjecture of Langlands and Rapoport,
\rm manuscript November 2009, http://arxiv.org/abs/0707.3173.

\Ref[Oo]
F. Oort, 
\sl Newton polygons and formal groups: conjectures by Manin and Grothendieck, 
\rm Ann. of Math. (2) {\bf 152} (2000), no. 1,  183--206.

\Ref[Pf]
M. Pfau,
\sl The reduction of connected Shimura varieties at a prime of good reduction,
\rm Ph. D. thesis (1993), Univ. of Michigan, U.S.A.

\Ref[Pi]
R. Pink,
\sl $l$-adic algebraic monodromy groups, cocharacters, and the Mumford--Tate conjecture,
\rm J. Reine Angew. Math. {\bf 495} (1998),  187--237.

\Ref[Re1]
H. Reimann,
\sl The semi-simple zeta function of quaternionic Shimura varieties,
\rm Lecture Notes in Math., Vol. {\bf 1657}, Springer-Verlag, Berlin, 1997.

\Ref[Re2]
H. Reimann,
\sl Reduction of Shimura varieties at parahoric levels,
\rm Manuscripta Math. {\bf 107}  (2002),  no. 3,  355--390. 

\Ref[RaZ]
M. Rapoport and T. Zink,
\sl Period spaces for $p$-divisible groups, 
\rm Annals of Mathematics Studies, Vol. {\bf 141}, Princeton University Press, Princeton, NJ, 1996.

\Ref[ReZ]
H. Reimann and T. Zink,
\sl Der Dieudonn\'emodul einer polarisierten abelschen Mannigfaltigkeit vom CM-typ,
\rm Ann. of Math. (2) {\bf 128} (1988), no. 3,  461--482.

\Ref[RR] 
M. Rapoport and M. Richartz, 
\sl On the classification and specialization of $F$-isocrystals with additional structure, 
\rm Compositio Math. {\bf 103} (1996), no. 2,  153--181.

\Ref[Sa1]
I. Satake,
\sl Holomorphic imbeddings of symmetric domains into a Siegel space,
\rm  Amer. J. Math. {\bf 87} (1965),  425--461.

\Ref[Sa2]
I. Satake,
\sl Symplectic representations of algebraic
groups satisfying a certain analyticity condition,
\rm Acta Math. {\bf 117} (1967),  215--279.

\Ref[Se1] J.-P. Serre, 
\sl Groupes alg\'ebriques associ\'es aux modules de Hodge--Tate, 
\rm Journ\'ees de G\'eom. Alg. de Rennes (Rennes, 1978),  Vol. III,   155--188, Ast\'erisque {\bf 65}, Soc. Math. France, Paris, 1979.

\Ref[Se2]
 J. -P. Serre, 
\sl Galois Cohomology, 
\rm Springer-Verlag, Berlin, 1997.

\Ref[Ta1] 
J. Tate,
\sl $p$-divisible groups,
\rm Proc. Conf. Local Fields (Driebergen, 1966), 158--183, Springer, Berlin, 1967.

\Ref[Ta2] 
J. Tate, 
\sl Classes d'isog\'enie des vari\'et\'es sur un corps fini (d'apr\`es J. Honda), 
\rm S\'eminaire Bourbaki 21\`eme ann\'ee (1968/69), , Exp. 352, Lecture Notes in Math., Vol. {\bf 179}, 95--110, Springer-Verlag, Berlin, 1971. 

\Ref[Ti1] 
J. Tits, 
\sl Classification of algebraic semisimple groups, 
\rm Algebraic Groups and Discontinuous Subgroups (Boulder, CO, 1965),   33--62, Proc. Sympos. Pure Math., Vol. {\bf 9}, Amer. Math. Soc., Providence, RI, 1966.

\Ref[Ti2]
J. Tits,
\sl Reductive groups over local fields, 
\rm Automorphic forms, representations and $L$-functions (Oregon State Univ., Corvallis, OR, 1977), Part 1,   29--69, Proc. Sympos. Pure Math., Vol. {\bf 33}, Amer. Math. Soc., Providence, RI, 1979.

\Ref[Va1]
A. Vasiu, 
\sl Integral canonical models of Shimura varieties of preabelian
 type, 
\rm Asian J. Math. {\bf 3} (1999), no. 2,  401--518.

\Ref[Va2]
A. Vasiu, 
\sl On two theorems for flat, affine groups schemes over a discrete valuation ring,
\rm Centr. Eur. J. Math. {\bf 3} (2005), no. 1,   14--25.

\Ref[Va3]
A. Vasiu,
\sl Integral canonical models of unitary Shimura varieties,
\rm Asian J. Math. {\bf 12} (2008), no. 2,  151--176.

\Ref[Va4]
A. Vasiu,
\sl Some cases of the Mumford--Tate conjecture and Shimura varieties,
\rm Indiana Univ. Math. J. {\bf 57} (2008), no. 1, 1--75.

\Ref[Va5]
A. Vasiu, 
\sl Projective integral models of Shimura varieties of Hodge type with compact factors,
\rm J. Reine Angew. Math. {\bf  618} (2008), 51--75. 

\Ref[Va6] 
A. Vasiu,
\sl Mod $p$ classification of Shimura $F$-crystals,
\rm Math. Nachr. {\bf 283} (2010), no. 8, 1068--1113.

\Ref[Va7]
A. Vasiu,
\sl Manin problems for Shimura varieties of Hodge type,
\rm J. Ramanujan Math. Soc. {\bf 26} (2011), no. 1, 31--84.

\Ref[Va8]
A. Vasiu,
\sl A motivic conjecture of Milne,
\rm 60 pages, to appear in J. Reine Angew. Math., http://arxiv.org/abs/math/0308202.

\Ref[Va9]
A. Vasiu,
\sl Generalized Serre--Tate ordinary theory,  
\rm 196 pages (including a 3 page index), accepted for publication in final from (i.e., to be published) by International Press, available at http://www.math.binghamton.edu/adrian/

\Ref[Va10]
A. Vasiu,
\sl On the Tate and Langlands--Rapoport conjectures for special fibres of integral canonical models of Shimura varieties of abelian type,
\rm manuscript 2008 available at http://www.math.binghamton.edu/adrian.

\Ref[VZ]
A. Vasiu and T. Zink,
\sl Breuil's classification of $p$-divisible groups over regular local rings of arbitrary dimension, 
\rm Algebraic and arithmetic structures of moduli spaces (Sapporo 2007),  461--479, Adv. Stud. Pure Math., Vol. {\bf 58}, Math. Soc. Japan, Tokyo, 2010. 

\Ref[Zi1]
T. Zink,
\sl Isogenieklassen von Punkten von Shimuramannigfaltigkeiten mit Werten in einem endlichen K\"orper,
\rm Math. Nachr. {\bf 112} (1983),  103--124.

\Ref[Zi2]
T. Zink,
\sl Windows for display of $p$-divisible groups,
\rm Moduli of abelian varieties (Texel Island, 1999),   491--518, Progr. Math. {\bf 195}, Birkh\"auser, Basel, 2001. 

\Ref[Wi]
J.-P. Wintenberger,
\sl Un scindage de la filtration de Hodge pour certaines vari\'et\'es alg\'ebriques sur les corps locaux,
\rm Ann. of Math. (2) {\bf 119} (1984), no. 3,  511--548.

}}

\bigskip
\hbox{Adrian Vasiu}
\hbox{Department of Mathematical Sciences, Binghamton University,}
\hbox{Binghamton,  P. O. Box 6000, New York 13902-6000, U.S.A.}
\hbox{e-mail: adrian\@math.binghamton.edu\;\;phone 1--607--777--6036.}

\enddocument